\documentclass[msom,nonblindrev]{informs3_mod}

\usepackage{amsmath,amssymb,mathtools}
\usepackage{threeparttable,booktabs,caption, graphicx, makecell, multicol, multirow, xcolor, lscape}
\usepackage{soul}
\usepackage{bm}
\usepackage{hyperref}
\usepackage{cleveref}
\usepackage{rotating}
\usepackage[shortlabels]{enumitem}
\usepackage{float}
\usepackage{datetime}
\usepackage{xcolor}
\usepackage{subfigure}
\usepackage{cases}

\OneAndAHalfSpacedXI

\newcommand{\Ascr}{\ensuremath{\mathcal A}}
\newcommand{\Bscr}{\ensuremath{\mathcal B}}

\newcommand{\Fscr}{\ensuremath{\mathcal F}}

\newcommand{\Iscr}{\ensuremath{\mathcal I}}

\newcommand{\Pscr}{\ensuremath{\mathcal P}}

\newcommand{\Rscr}{\ensuremath{\mathcal R}}
\newcommand{\Sscr}{\ensuremath{\mathcal S}}
\newcommand{\Tscr}{\ensuremath{\mathcal T}}

\newcommand{\Wscr}{\ensuremath{\mathcal W}}
\newcommand{\Xscr}{\ensuremath{\mathcal X}}
\newcommand{\Yscr}{\ensuremath{\mathcal Y}}
\newcommand{\Zscr}{\ensuremath{\mathcal Z}}

\def\E{\mathbb{E}}
\def\R{\mathbb{R}}
\def\P{\mathbb{P}}

\newcommand{\paran}[1]{\left\{#1\right\}}
\newcommand{\Paran}[1]{\left[#1\right]}
\newcommand{\Para}[1]{\left(#1\right)}

\newcommand{\I}[1]{\ensuremath{\mathbb{I}_{\left\{#1\right\}}}} 

\newcommand{\setaiw}{\ensuremath{\Ascr_{\textup{\textsf{IW}}}}}
\newcommand{\iwdro}{\textsf{IW-DRO}}
\newcommand{\ziw}{\ensuremath{z_{\textup{\textsf{IW}}}}}
\newcommand{\Jiw}{\ensuremath{J_{\textup{\textsf{IW}}}}}
\newcommand{\eapx}{\ensuremath{\epsilon_{\textup{\textsf{apx}}}}}
\newcommand{\varepsilonnp}{\ensuremath{\varepsilon_{\textup{\textsf{NP}}}}}
\newcommand{\varepsilonp}{\ensuremath{\varepsilon_{\textup{\textsf{P}}}}}
\newcommand{\rnp}{\ensuremath{r_{\textup{\textsf{NP}}}}}
\newcommand{\rp}{\ensuremath{r_{\textup{\textsf{P}}}}}
\newcommand{\rme}{\ensuremath{r_{\textup{\textsf{ME}}}}}
\newcommand{\setame}{\ensuremath{\Ascr_{\textup{\textsf{ME}}}}}
\newcommand{\varepsilonme}{\ensuremath{\varepsilon_{\textup{\textsf{ME}}}}}
\newcommand{\zme}{\ensuremath{z_{\textup{\textsf{ME}}}}}
\newcommand{\Cme}{\ensuremath{C_{\textup{\textsf{ME}}}}}

\newcommand{\munp}{\ensuremath{\hat{\mu}_{\textup{\textsf{NP}}}}}
\newcommand{\mup}{\ensuremath{\hat{\mu}_{\textup{\textsf{P}}}}}
\newcommand{\mume}{\ensuremath{\hat{\mu}_{\textup{\textsf{ME}}}}}

\newcommand{\iwdroapx}{\textsf{IW-DRO-Apx}}

\newcommand{\setag}{\ensuremath{\Ascr_{\textup{\textsf{G}}}}}
\newcommand{\JP}{\ensuremath{J_{\textup{\textsf{P}}}}}
\newcommand{\JD}{\ensuremath{J_{\textup{\textsf{D}}}}}
\newcommand{\mur}{\ensuremath{\hat{\mu}_{\textup{\textsf{R}}}}}
\newcommand{\muru}{\ensuremath{\hat{\mu}_{\textup{\textsf{RU}}}}}
\newcommand{\muapx}{\ensuremath{{\mu}_{\textup{\textsf{apx}}}}}

\usepackage{natbib}
 \bibpunct[, ]{(}{)}{,}{a}{}{,}%
 \def\bibfont{\small}%
\TheoremsNumberedThrough     
\ECRepeatTheorems

\EquationsNumberedThrough    

\MANUSCRIPTNO{}

\begin{document}


\RUNAUTHOR{Wang, Chen, and Wang}

\RUNTITLE{Contextual Optimization under Covariate Shift}

\TITLE{Contextual Optimization under Covariate Shift:\\
A Robust Approach by Intersecting Wasserstein Balls}

\ARTICLEAUTHORS{%
\AUTHOR{Tianyu Wang}
\AFF{Department of Industrial Engineering and Operations Research, Columbia University, New York, NY 10027,
\EMAIL{tianyu.wang@columbia.edu}} 
\AUTHOR{Ningyuan Chen}
\AFF{Department of Management, University of Toronto Mississauga,\\
Rotman School of Management, University of Toronto, Canada,
\EMAIL{ningyuan.chen@utoronto.ca}} 
\AUTHOR{Chun Wang}
\AFF{School of Economics and Management, Tsinghua University, Beijing, China,
\EMAIL{wangchun@sem.tsinghua.edu.cn}} 

} 

\ABSTRACT{%
In contextual optimization, a decision-maker leverages contextual information, often referred to as covariates, to better resolve uncertainty and make informed decisions.
In this paper, we examine the challenges of contextual decision-making under covariate shift, a phenomenon where the distribution of covariates differs between the training and test environments. Such shifts can lead to inaccurate upstream estimations for test covariates that lie far from the training data, ultimately resulting in suboptimal downstream decisions.
To tackle these challenges, we propose a novel approach called Intersection Wasserstein-balls DRO (\textsf{IW-DRO}), which integrates multiple estimation methods into the distributionally robust optimization (DRO) framework. At the core of our approach is an innovative ambiguity set defined as the intersection of two Wasserstein balls, with their centers constructed using appropriate nonparametric and parametric estimators. On the computational side, we reformulate the \textsf{IW-DRO} problem as a tractable convex program and develop an approximate algorithm tailored for large-scale problems to enhance computational efficiency.
From a theoretical perspective, we demonstrate that \textsf{IW-DRO} achieves superior performance compared to single Wasserstein-ball DRO models. We further establish performance guarantees by analyzing the coverage of the intersection ambiguity set and the measure concentration of both estimators under the Wasserstein distance. Notably, we derive a finite-sample concentration result for the Nadaraya-Watson kernel estimator under covariate shift.
The proposed \textsf{IW-DRO} framework offers practical value for decision-makers operating in uncertain environments affected by covariate shifts. Through extensive numerical experiments on income prediction and portfolio optimization problems, using both synthetic and real-world data, we demonstrate that \textsf{IW-DRO} outperforms various benchmark policies in these challenging settings.
}


\HISTORY{This draft is of \today.\footnote{The earlier draft \citep{wang2021distributionally}, titled ``Distributionally Robust Prescriptive Analytics with Wasserstein Distance'', focuses on the measure concentration properties of the kernel estimator under Wasserstein distance, without considering covariate shifts. This content now constitutes Section~\ref{subsec:estimator-np} of the current paper.}
}

\maketitle
\section{Introduction}\label{sec:intro}
In the current era of big data, the increasing availability of contextual information, often referred to as covariates, empowers data-driven decision-making by enabling more effective resolution of uncertainty.
This growing interest has given rise to the field of \emph{contextual (stochastic) optimization}, which has attracted considerable scholarly attention in recent years.
See \cite{sadana2024survey} for a comprehensive review of various optimization methods and their applications.
In this paper, we examine contextual optimization challenges that arise from \emph{covariate shift}, where the covariate distribution diverges between the training and test environments.

\subsection{Background}\label{sec:intro-1}
A generic formulation of a contextual optimization problem is as follows:
\begin{equation}
\label{eq:generic-cso}
\min_{z \in \Zscr} ~ \E\left[ c(z,Y) | X=x \right].
\end{equation}
In this framework,
a decision maker (DM) aims to choose a decision vector $z \in \Zscr \subseteq \R^{d_z}$
to minimize a cost function $c(z,Y)$,
which depends on a random vector $Y \in \Yscr \subseteq \R^{d_y}$ that characterizes the uncertain parameters in the problem setting.
While the precise realization of $Y$ is not revealed at the time of decision-making,
the DM has access to a covariate vector $X \in \Xscr \subseteq \R^{d_x}$, which proffers side information about $Y$ and thereby helps the DM to infer the uncertainty.
Ideally, given a newly observed covariate $X=x$, the DM is interested in minimizing the expected cost in \eqref{eq:generic-cso}
with respect to the ground-truth conditional distribution of $Y$ given $x$, denoted by $\mu_{Y|x}$.
However, in real-world applications, $\mu_{Y|x}$ is typically unknown.
Instead, there is a training dataset $\{(x_i, y_i)\}_{i = 1}^{n}$ available,
which comprises $n$ samples collected from past observations of uncertain parameters and their associated covariates.
Consequently, the DM may resort to data-driven optimization procedures to make informed decisions.
A classical approach is the \emph{estimate-then-optimize} method,
where the DM first uses the training samples $\{(x_i, y_i)\}_{i = 1}^{n}$ and the test covariate $x$
to estimate a conditional distribution $\hat\mu_{Y|x}$ as an approximation to the true $\mu_{Y|x}$,
and then minimizes the expected cost over the estimated $\hat\mu_{Y|x}$,
that is, to solve $\displaystyle{\min_{z \in \Zscr} ~ \E_{\hat\mu_{Y|x}}\left[ c(z,Y)\right]}$ (e.g., \citealp{bertsimas2020predictive}).
Moreover, to safeguard against overly optimistic decisions due to estimation errors,
the DM can select an ambiguity set $\Ascr$
and apply the distributional robust optimization (DRO) framework to minimize the corresponding worst-case expected cost
for any distribution $\mu \in \Ascr$, as follows:
\begin{equation}
\label{eq:generic-cso-dro}
\min_{z \in \Zscr} \sup_{\mu \in \Ascr} ~\E_{Y\sim\mu}[c(z,Y)].
\end{equation}
Intuitively, $\Ascr$ should include all distributions sufficiently likely to be the true one $\mu_{Y|x}$.
A popular choice is to construct $\Ascr = \{\mu: \Wscr_p(\mu, \hat\mu_{Y|x})\leq \varepsilon\}$
as a ball of distributions centered around an estimated conditional distribution $\hat\mu_{Y|x}$
with a radius $\varepsilon$ measured by some distance, such as the $p$-Wasserstein distance $\Wscr_p$
(e.g., \citealp{kannan2020residuals}).
Empirical studies have suggested that the Wasserstein-DRO approach, formulated in \eqref{eq:generic-cso-dro},
yields favorable out-of-sample performance.

While there has been significant growth in the literature on data-driven methods for contextual optimization problems,
most studies assume a stationary environment for covariates,
meaning that the covariate distribution remains the same between the training and test datasets.
In practice, however, this assumption often fails,
leading current contextual optimization methods to produce severely suboptimal decisions.
This phenomenon is known as covariate shift \citep{shimodaira2000improving}.
We illustrate this concept with the following examples.

$\bullet$ A public policymaker allocates resources to specific population groups based on income levels,
utilizing individual features $X$ (covariates such as education, occupation, etc.) to predict their income $Y$.
If the prediction model is trained using data from one city and then deployed in another,
or if the training dataset is based on a survey where younger people are more likely to participate while the test dataset consists of people from a broader range of ages,
the covariate distributions in the population might differ due to environmental or demographic factors.
As a result, the income predictions may become less effective and lead to misguided resource allocation decisions.

$\bullet$ A financial institution develops a model to manage portfolio risks by utilizing various market information $X$ (covariates such as interest rates, sectors moment, etc)
to estimate asset returns $Y$.
Suppose the model is trained using historical data from a period of economic strength.
However, when deployed, the model may face different market conditions.
For instance, an economic downturn could alter the distribution of interest rates,
potentially leading the model to suggest portfolio positions with undesired risk exposures.

In both cases, if covariate shift is not properly addressed in the estimate-then-optimize approach,
the upstream estimation step may yield inaccurate predictions given covariates significantly differ from the training samples. Consequently, the downstream optimization step is likely to produce poor decisions.
However, if the relationship between covariates and uncertain parameters remains stable
(for example, in the above two cases: higher education continues to correlate with higher income, and lower interest rates consistently drive bond prices higher),
there remains an opportunity to improve the contextual optimization model to suit a new covariate environment. 

\subsection{Motivation and Contribution}\label{sec:intro-2}
We examine contextual optimization problems in the presence of covariate shift, which is defined as follows:
\begin{definition}[Covariate Shift]\label{def:cs-shift}
Let $\mu_{X}$ denote the marginal probability measure of covariate $X$ in the training data distribution,
and let $\nu_{X}$ denote that in the test data distribution.
The following conditions hold within the training and test data distributions:
(i) $\mu_X \neq \nu_X$, i.e.,
the distribution of $X$ changes;
and (ii) the conditional distribution of $Y$ given $x$ remains same as $\mu_{Y|x}$.
\end{definition}

The assumption of a stable conditional distribution $\mu_{Y|x}$ between the training and test datasets is common in the machine learning literature \citep{sugiyama2012machine}.
Without this assumption, the patterns learned from the training data may become irrelevant for predicting uncertainty in the test phase.
In this paper, we particularly use the term \emph{covariate shift} to describe situations where the test covariates fall into regions that are
insufficiently represented by the training samples.

To mitigate the negative impacts of covariate shift on decision-making,
we develop an algorithmic framework
that integrates various estimation methods with the DRO formulation \eqref{eq:generic-cso-dro}.
Specifically, we propose a novel ambiguity set, $\setaiw$, which is
defined as the intersection of two Wasserstein balls:
\begin{align}
\label{eq:ambiguity-set}
\setaiw = \{\mu: \Wscr_p(\mu, \munp)\leq \varepsilonnp, ~\Wscr_p(\mu, \mup)\leq \varepsilonp \}.
\end{align}
Here, the reference distributions $\munp$ and $\mup$, which serve as the centers of the two Wasserstein balls,
are selected based on certain \emph{nonparametric} and \emph{parametric} estimators of $\mu_{Y|x}$, respectively.
Meanwhile, the radii of these balls, $\varepsilonnp$ and $\varepsilonp$, are chosen accordingly.
Using this ambiguity set,
decisions are derived by solving the following min-max optimization problem over $\setaiw$:
\begin{align}
\label{eq:iwdro-opt}
\min_{z \in \Zscr} \sup_{\mu \in \setaiw} ~ & \E_{Y\sim\mu}[c(z,Y)].  \tag{{\iwdro}}
\end{align}
Throughout the paper, we use the abbreviation ``IW'' to indicate the concept of Intersecting Wasserstein-balls,
e.g., see the subscript of $\setaiw$.
We refer to our proposed framework as {\iwdro} and denote the associated optimization problem using \eqref{eq:iwdro-opt} with parentheses.

In the design of $\setaiw$,
the rationale for utilizing the intersection of Wasserstein balls is twofold:
it can (i) combine the advantages of both nonparametric and parametric estimators;
and (ii) yield a more compact ambiguity set.
To see this, we first argue that there is a fundamental trade-off in statistical learning akin to the bias-variance dilemma:
a ``statistically-robust'' estimator, such as a simple parametric one,
may generalize well to shifted covariates due to its simplicity, but might not fully capture the complexity of the data when there is an abundance of training samples in the neighborhood;
on the other hand, nonparametric estimators, while potentially capturing the nuances of the local data structure, may underperform when faced with shifted covariates from less explored regions.
Therefore, by positing $\mu_{Y|x}$ within the intersection of Wasserstein balls formed by these estimators, {\iwdro} harnesses the strengths of both methods regardless of whether the new covariate has shifted.
Furthermore, in a DRO framework,
while it is important to contain $\mu_{Y|x}$,
an overly broad ambiguity set can lead to excessively conservative (and thus potentially suboptimal) decisions due to incorporating implausible worst-case scenarios.
By intersecting two Wasserstein balls,
we can selectively exclude certain improbable distributions
and maintain a high likelihood of containing $\mu_{Y|x}$.
We next consider a detailed example to further illustrate the idea.

\begin{example}\label{ex:motivate}
In a prediction task, the uncertainty $Y$ follows a normal distribution with a mean conditional on the covariate $x$, specifically, $\mu_{Y|x} \sim \mbox{N}(2x + 0.3 \sin(4\pi x), ~0.3)$.
There are a total of $n=204$ training samples, comprising two groups:
(i) 200 samples generated with $x$ values drawn from a uniform distribution $\mbox{U}(0, 1)$,
and (ii) 4 samples generated with $x$ values at $\{1.1, 1.2, 1.3, 1.4\}$.
As a result, the distribution of $X$ in the training dataset is heavily skewed towards the interval $[0, 1]$.
In contrast, the test dataset samples $X$ uniformly from $\mbox{U}(0, 1.4)$, leading to a covariate shift in the region $[1, 1.4]$,
since the training data is sparse in this region. \hfill\halmos
\end{example}

\begin{figure}[ht]
\SingleSpacedXI
\begin{center}
\caption{\hspace{-0.4cm} Comparison of the nonparametric estimator (NW kernel) with the parametric estimator (linear regression)}
\label{fig:motivation}
\subfigure[{Confidence intervals of two different estimators} ]{\includegraphics[width = 0.49\textwidth, trim = 0 10 0 10]{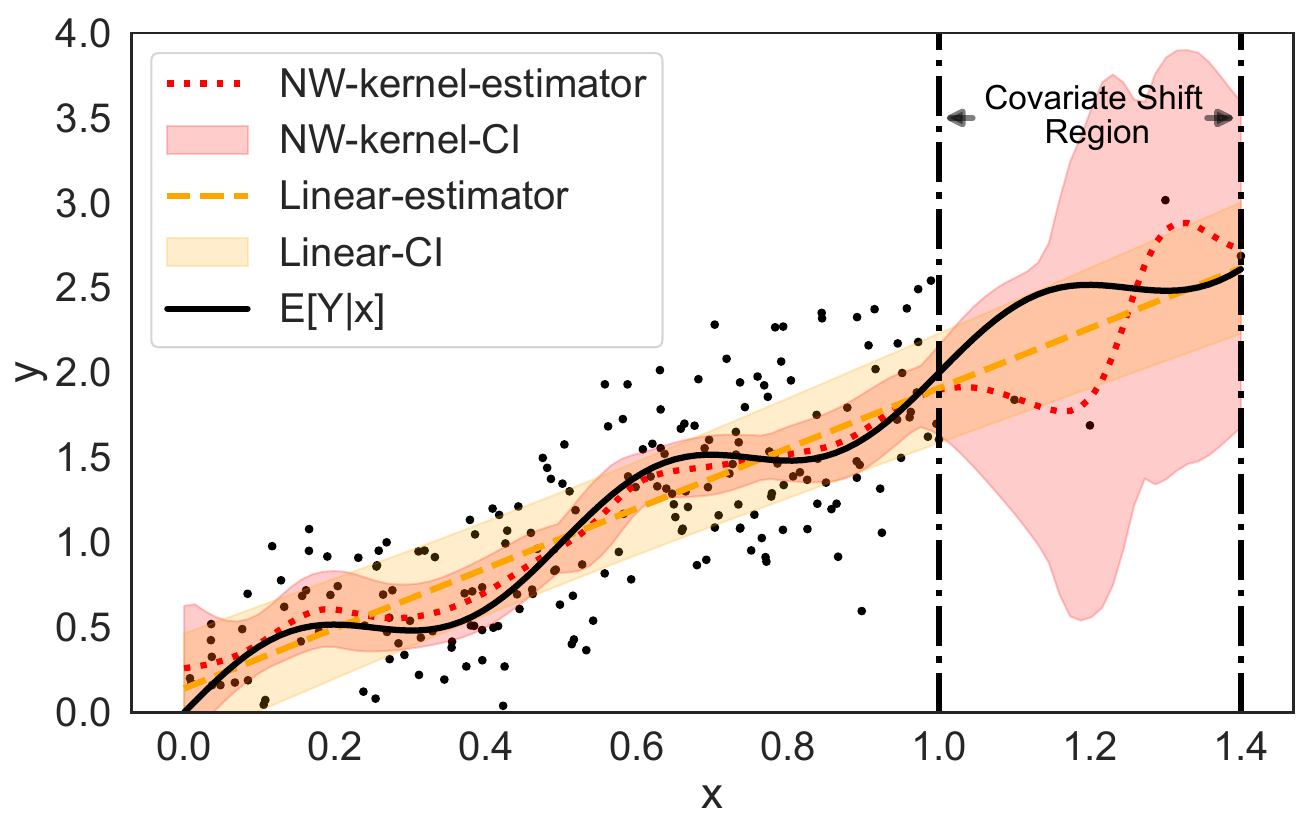}}
\subfigure[{The intersection of two confidence intervals}]{\includegraphics[width = 0.49\textwidth, trim = 0 10 0 10]{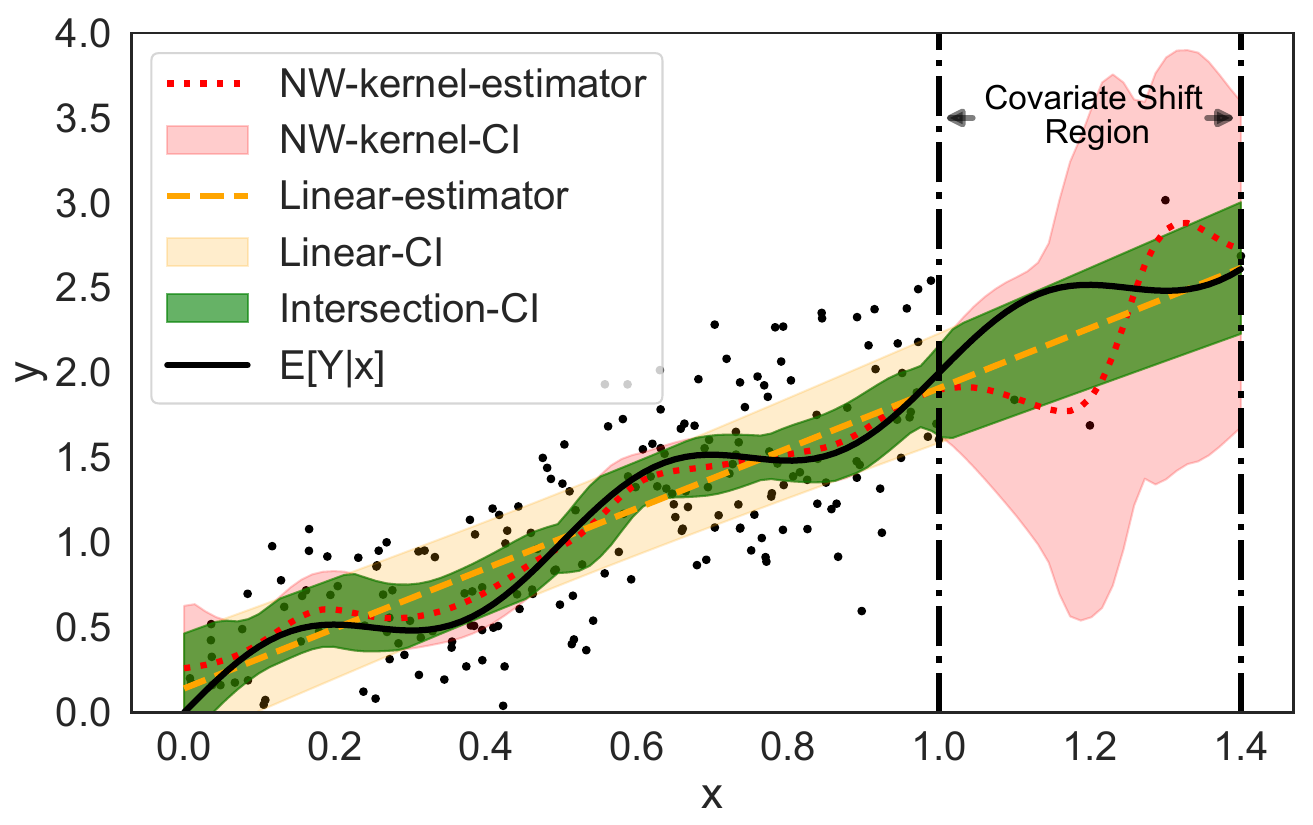}}
\end{center}
{\footnotesize {\em Note}: Each estimator is fitted using the same set of $n=204$ training samples (black points), and their respective confidence intervals (CIs) are constructed at a 90\% confidence level.}
\vspace{-0cm}
\end{figure}

We illustrate the performances of two different estimators using \Cref{ex:motivate}.
When the training dataset contains enough samples around the test covariate,
certain nonparametric estimators, such as the Nadaraya–Watson (NW) kernel estimator \citep{nadaraya1964estimating,watson1964smooth},
are well known to achieve asymptotic consistency and often outperform parametric estimators.
As shown in \Cref{fig:motivation}(a), for the region $x \in [0, 1]$ with sufficient training samples,
the NW kernel estimator produces predictions (the red dot line) that are very close to the true conditional means (the black solid line), accompanied by a tight confidence interval (the pink area).

However, when the test covariate falls outside the region well covered by the training samples, the kernel estimator becomes highly unstable due to its reliance on local data points.
In contrast, parametric estimators, such as linear regression,
may provide more reliable estimates in this situation, assuming that the regression model is specified close to the true data generating process.
As also shown in \Cref{fig:motivation}(a), when there are few training samples available for $x > 1$,
the NW kernel estimator becomes drastically ineffective and yields a very loose confidence interval.
In this case, a linear regression estimator offers more reasonable predictions (the orange dash line) with a tighter confidence interval (the yellow area).

This example illustrates that nonparametric kernel estimators and parametric regression estimators each have distinct advantages and limitations, arising from their differing sensitivity to local sample density under the influence of covariate shifts.
Notably, the intersection of their confidence intervals, depicted as the green area in \Cref{fig:motivation}(b),
creates a narrow region with a high likelihood of encompassing the true conditional mean value,
i.e., the ideal prediction, across the entire covariate domain.
When generalizing this intuition to {\iwdro}, 
it suggests that the intersection of Wasserstein balls centered on these estimators could yield a compact ambiguity set capable of effectively capturing the true conditional distribution.
Instead of using a single Wasserstein ball with radius $\varepsilon_1'$ associated with one estimator,
the DM can choose two balls formed by different estimators with radii $\varepsilon_1$ and $\varepsilon_2$, each slightly larger than $\varepsilon_1'$.
This configuration allows
their intersection to include more plausible distributions and offer greater resilience to covariate shifts,
while maintaining a smaller or comparable overall size of the ambiguity set.

To the best of our knowledge, {\iwdro} is the first model to address covariate shifts within contextual optimization problems.
Our main contributions are summarized as follows:

\textbf{Modeling contextual optimization under covariate shift:}
We aim to improve decision quality by accounting for covariate shifts at every stage of the estimate-then-optimize process, rather than focusing solely on prediction accuracy.
To achieve this, we integrate the strengths of various estimation methods within the DRO framework by introducing a novel ambiguity set, defined as the intersection of Wasserstein balls.
Our {\iwdro} model provides a unified solution for robust contextual decision-making, accommodating test covariates both near and far from the training samples.

\textbf{Computational reformulation:}
For a piecewise concave cost function, we show that the \eqref{eq:iwdro-opt} optimization problem (a min-max program over the intersection ambiguity set) can be reformulated as a tractable convex program.
This reformulation is highly nontrivial and may be of independent interest from a computational perspective.
Additionally, we propose an approximate model that significantly enhances computational efficiency for large-scale problems.

\textbf{Performance guarantees and measure concentration:}
We demonstrate that {\iwdro} can yield superior decisions compared to those derived from single Wasserstein-ball DRO models.
We establish performance guarantees by examining the coverage of the ambiguity set.
In particular, we develop a finite-sample concentration result for the NW kernel estimator, measured under the Wasserstein distance.
This constitutes the first nonasymptotic rate result,
a topic that was not previously explored even in settings without covariate shifts.

\textbf{Empirical evaluation:}
Through extensive numerical experiments on income prediction and portfolio optimization problems, using both synthetic and real-world data, we show that \textsf{IW-DRO} consistently outperforms various benchmark policies in situations involving covariate shifts.

\subsection{Literature Review}\label{sec:intro-3}

\paragraph{Optimization Models in Contextual Optimization.}
There are three streams in the literature for solving contextual optimization problems:
estimate-then-optimize (e.g., \citealp{ban2019big,bertsimas2020predictive}), integrate-estimation-and-optimization (e.g., \citealp{elmachtoub2017smart,qi2021integrated}), and end-to-end learning (e.g., \citealp{donti2017task,ho2022risk,qi2022e2e}).
For a review of various models, see \cite{sadana2024survey}.
To the best of our knowledge, none of the existing studies within these frameworks have yet to incorporate covariate shifts into their formulations.

Our study on contextual optimization aligns with the estimate-then-optimize stream.
In this paradigm, recent research has employed DRO to mitigate estimation errors and enhance decision-making,
with the Wasserstein ambiguity set being a popular choice.
(The Wasserstein-DRO framework has its roots in classical stochastic programming without covariates~\citep{esfahani2018data,blanchet2019quantifying,gao2016distributionally,long2023robust}. There is a substantial body of related studies in this field. For a comprehensive survey, see \citealp{kuhn2019wasserstein}.)
Drawing on contextual information, \cite{kannan2020residuals} suggest fitting a regression model to correlate uncertainty with covariates, followed by selecting a Wasserstein ball centered on a residual-based distribution; \cite{bertsimas2023dynamic} propose placing a Wasserstein ball around a local conditional distribution estimator, such as nearest neighbors and kernel regression.
In operations tasks with covariates data, the Wasserstein-DRO models have been applied in various areas,
including portfolio optimization \citep{nguyen2021robustifying,pun2023data}, power flow \citep{esteban2023distributionally}, newsvendor problems \citep{zhang2023optimal}, and multi-product pricing \citep{sim2024analytics}.
Besides the Wasserstein metric, other statistical distance measures have been investigated when defining ambiguity sets, such as the KL divergence in \cite{bertsimas2021bootstrap} and the causal transport distance in \cite{yang2022decision}.
It has been observed that existing DRO models maintain adequate out-of-sample performance when training and test distributions align.
However, these models struggle with covariate shifts, particularly because they lack mechanisms to adjust the ambiguity set based on the density of the training samples.
Such adjustments are crucial to ensure the inclusion of the target conditional distribution and, consequently, to provide satisfactory performance guarantees.
On the other hand, empirical studies have shown that an overly conservative ambiguity set, based on implausible distribution shifts, can result in substantial performance degradation \citep{hu2018does,liu2024need}.
To tackle these challenges, we propose {\iwdro} with the intersection ambiguity set specifically designed to account for covariate shifts,
and we adopt the Wasserstein metric due to its advantageous statistical properties.

We note that the intersection of Wasserstein ambiguity sets has recently gained attention in the research community for modeling problems with multiple data sources.
\cite{awasthi2022distributionally} and \cite{selvican} construct Wasserstein balls centered at the empirical distributions of each training data source and apply their intersection set in the context of robust logistic regression.
Concurrent to our work, \cite{rychener2024wasserstein} provide a duality result similar to our \Cref{thm:strong-duality} for a distributionally robust stochastic program.
These formulations do not incorporate contextual information and, as such, they do not address conditional distribution estimators or covariate shifts.
In addition, some studies examine the intersections of different types of ambiguity sets to reduce conservativeness, such as moment-based ambiguity sets \citep{chen2019distributionally}, moment and Wasserstein balls \citep{cheramin2022computationally}, Wasserstein and goodness of fitness balls \citep{tanoumand2023data}, leading to optimization models that differ from contextual optimization.

\paragraph{Decision-making under Distribution Shifts.}
In real-world scenarios, distribution shifts are common and pose significant challenges for decision-making \citep{quinonero2008dataset}.
Recent studies have explored statistical guarantees for optimization models that face data contamination or particular types of distribution shifts.
For example, \cite{besbes2022beyond} and \cite{gupta2022data} examine problem-dependent optimization models where training samples come from different distributions;
\cite{bennouna2022holistic} propose a DRO approach to protect against corruptions in historical data with limited precision.
These settings differ from our focus, which is on covariate shifts between the training and test data distributions.
Additionally, \cite{duchi2020distributionally} investigate a DRO model for prediction tasks under the worst-case covariate shift scenario.
However, their DRO model may be overly conservative for our purposes, as it does not adapt to new test covariate observations.
In the related field of statistical learning,
methods designed to counteract covariate shifts,
such as importance weighting \citep{sugiyama2007covariate} and distribution matching \citep{bickel2009discriminative},
typically require a set of covariate observations from the test environment.
Therefore, these methods are not suitable for our contextual optimization problems, which necessitate making decisions based on each new observation of the test covariate as it appears.

\paragraph{Statistical Estimation.}
Our {\iwdro} approach utilizes a nonparametric NW kernel estimator and a parametric linear regression estimator to learn the conditional distribution between uncertain parameters and covariates.
For general background information on these estimators, one can refer to \cite{hastie2009elements}.
In particular, we choose the NW kernel estimator due to its widespread use in practice,
including stochastic optimization \citep{hannah2010nonparametric}, dynamic programming \citep{hanasusanto2013robust}, and newsvendor problem \citep{ban2019big}.

\smallskip
\noindent
\textbf{Paper Outline.} The remainder of this paper is organized as follows.
Section \ref{sec:IWDRO-formulation} provides the computational reformulation and performance guarantees of {\iwdro}.
In Section \ref{sec:statistical-result},  we delve into the construction of the intersection ambiguity set and examine its statistical properties.
In \Cref{sec:approximation-formulation}, we provide a computational surrogate to {\iwdro} with provable performance guarantees.
\Cref{sec:numerical} presents extensive numerical results on both synthetic and real-world datasets, including portfolio optimization and income prediction problems, to demonstrate the benefit of {\iwdro}.
The detailed proofs can be found in the appendices.

\smallskip
\noindent \textbf{Notations}:
For the covariate and uncertainty,
we use uppercase letters ($X$ and $Y$) to represent the random variables
and lowercase letters ($x$ and $y$) for their samples or realizations.
The feasible sets for covariate, uncertainty, and decision
are represented by $\Xscr$, $\Yscr$, and $\Zscr$, respectively.
For $n\in \mathbb{N}$, we set $[n] \coloneqq \{1,2,...,n\}$.
We let $\|\cdot\|_p$ denote the $\ell_p$-norm,
$\delta$ denote the Dirac measure,
$\I{\cdot}$ denote the indicator function,
$\Pscr(\Sscr)$ denote the space of probability distributions within support set $\Sscr$,
and $\Xi(\mu, \nu)$ denote the set of all joint distributions with marginal distributions $\mu$ and $\nu$.
We let $\sigma_{\Sscr}$ denote  the support function of set $\Sscr$, i.e., $\sigma_\Sscr (u)=\sup_{s \in \Sscr}\{s^{\top} u\}$,
and $f^*$ denote the conjugate of a convex function $f:\R^{d}\rightarrow(-\infty, \infty)$, i.e., $f^*(u) = \sup_{s\in \R^d}\{s^{\top} u\} - f(s)$.
We use $\Theta$-notations to describe the asymptotic performance of functions.
Specifically, given two functions $f(n)$ and $g(n)$, $g(n)$ is said to be $\Theta(f(n))$ if there exist positive constants $c_1$, $c_2$, and $n_0$ such that $c_1 f(n) \leq g(n) \leq c_2 f(n)$ for all $n \geq n_0$. 
We use $\mathbf{I}$ to denote an identity matrix, and use $\mathbf{0}$ and $\mathbf{1}$ to denote vectors of zeros and ones, respectively.
The $p$-Wasserstein distance is defined as follows:
\begin{definition}[Wasserstein distance]\label{def:was-distance}
The $p$-Wasserstein distance $(p\in \mathbb{N})$ between two distributions $\mu,\nu \in \Pscr(\Yscr)$
with respect to the $l_1$-norm (i.e., $\|\cdot\|_1$) is defined as:
\begin{equation}\label{eq:pwasserstein-dist}
\Wscr_p(\mu,\nu) \coloneqq \inf_{\xi \in \Xi(\mu, \nu)} \left( \E_{(Y_1, Y_2)\sim \xi}\Big[\|Y_1 - Y_2\|_1^p \Big] \right)^{\frac{1}{p}} ,
\end{equation}
where $\xi$ is a joint distribution of $(Y_1 \in \Yscr, Y_2 \in \Yscr)$
from $\Xi(\mu, \nu)$.
For $p \le q$, $\Wscr_p(\mu,\nu) \le \Wscr_q(\mu,\nu)$.
\end{definition}
Here, we use the $l_1$-norm in definition \eqref{eq:pwasserstein-dist} for simplicity. Our computational and statistical results can be generalized to the Wasserstein distance defined with other metric functions.

\section{{\iwdro} Computation}\label{sec:IWDRO-formulation}
We now investigate the solution approach to the optimization problem \eqref{eq:iwdro-opt}:
$$\displaystyle{\min_{z \in \Zscr} \sup_{\mu \in \setaiw} \E_{Y\sim\mu}[c(z,Y)]}.$$
Recall that, as defined in \eqref{eq:ambiguity-set},
the ambiguity set $\setaiw = \{\mu: \Wscr_p(\mu, \munp)\leq \varepsilonnp, ~\Wscr_p(\mu, \mup)\leq \varepsilonp \}$ is constructed as the intersection of two Wasserstein balls centered at $\munp$ and $\mup$, which correspond to a pair of nonparametric and parametric estimators.
In this section, to focus on the computational aspects, we consider a more generic definition of the intersection set:
\begin{align}
\label{eq:ambiguity-set-generic}
\setaiw = \{\mu: \Wscr_p(\mu, \hat\mu_1)\leq \varepsilon_1, ~\Wscr_p(\mu, \hat\mu_2)\leq \varepsilon_2\},
\end{align}
which abstracts away the dependence on specific estimation methods.
Here, we assume that $\hat\mu_1$ and $\hat\mu_2$, the centers of two Wasserstein balls, follow a generic discrete distribution structure, as stated in the assumption below.
(The detailed specifications of $\munp$ and $\mup$ are provided in Section \ref{sec:statistical-result}.)

\begin{assumption}[Discrete Reference Distributions]\label{asp:dist_structure}
For $m\in\{1,2\}$, $\hat\mu_m$ in \eqref{eq:ambiguity-set-generic} is a discrete distribution on $\Yscr$
with a finite support set $\{y_{m,i}\}_{i \in [n]}$ and corresponding probability mass $\{w_{m,i}\}_{i \in [n]}$,
that is,  $\hat\mu_m = \sum_{i=1}^n w_{m,i} \delta_{y_{m,i}}$.
\end{assumption}
This assumption of discreteness is commonly adopted in the literature
since it facilitates computation.
It is also justified by
the inherently discrete characteristics of both empirical distributions (i.e., the training samples) and many estimated conditional distributions
(e.g., \citealp{bertsimas2023dynamic,kannan2020residuals}).

To solve \eqref{eq:iwdro-opt},
we begin by establishing a necessary and sufficient condition for its feasibility.
\begin{proposition}[Feasibility of \iwdro]\label{prop:intersect-ambiguity}
The ambiguity set $\setaiw$ defined in \eqref{eq:ambiguity-set-generic}
is non-empty if and only if $\Wscr_p(\hat\mu_1, \hat\mu_2)\leq \varepsilon_1 + \varepsilon_2$.
\end{proposition}
This proposition aligns with the rationale of Proposition 5.3 in \cite{taskesen2021sequential}.
Given that $\hat\mu_1$ and $\hat\mu_2$ have finite supports, as assumed in \Cref{asp:dist_structure},
let $p_{i,j}$ denote the probability of $(y_{1,i}, y_{2,j})$ for a distribution $\xi \in \Xi(\hat\mu_1, \hat\mu_2)$.
Then, by \Cref{def:was-distance},
$\Wscr_p(\hat\mu_1, \hat\mu_2)$ can be computed via the following transportation problem in linear programming
to ascertain the non-emptiness of $\setaiw$:
\begin{align}\label{eq:transportation-1}
\min_{p_{i,j} \ge 0} \Bigg\{ \sum_{i,j\in [n]}p_{i,j} \|y_{1,i} - y_{2,j}\|_1^p, ~\text{s.t.}  ~\sum_{j\in[n]}p_{i,j} = w_{1,i}, \forall i\in[n];   ~\sum_{i\in[n]}p_{i,j} = w_{2,j}, \forall j\in[n].  \Bigg\}
\end{align}

In addition, we also impose an assumption on the boundness of the cost function.
\begin{assumption}[Bounded Cost Function]\label{asp:reg-costfunc}
The expected cost is bounded, for $m\in\{1,2\}$,
$|\E_{Y\sim\hat\mu_m}[c(z,Y)]| < \infty, ~\forall z \in \Zscr$.
\end{assumption}
This assumption is standard in the Wasserstein-DRO literature, e.g., Assumption 1 in \cite{zhang2022simple},
and can be easily satisfied when $c(\cdot,\cdot)$ or $\Yscr$ is bounded.

\subsection{Reformulation of {\iwdro}}\label{sec:IWDRO-covex-reformulation}
Throughout this paper, we proceed under
Assumptions~\ref{asp:dist_structure} and \ref{asp:reg-costfunc} and assume that $\setaiw$ is non-empty.
Next, we present the main result for solving \eqref{eq:iwdro-opt}.
\begin{theorem}[Strong Duality of \iwdro]\label{thm:strong-duality}
For all $z \in \Zscr$, we have $ \JP(z) = \JD(z)$ with:
\begin{align}
\label{eq:strong-dual0-cp}
\JP(z) \coloneqq & ~ \sup_{\mu \in \setaiw} ~ \E_{Y \sim \mu}[c(z, Y)], \\
\label{eq:strong-dual-further}
\begin{split}
\JD(z) \coloneqq & \inf_{\substack{\lambda_1, \lambda_2 \geq 0, \\ u_i, v_j \in \mathbb{R}}} \bigg\{ \lambda_1 \varepsilon_1^p + \lambda_2 \varepsilon_2^p + \sum_{i \in [n]} w_{1,i}u_i + \sum_{j \in [n]} w_{2,j}v_j, \\
&\quad\quad\quad\quad \text{\textup{s.t.}} ~u_i + v_j \geq \sup_{y \in \Yscr} \Big\{c(z, y) - \lambda_1 \|y_{1,i} - y\|_1^p - \lambda_2 \|y_{2,j} - y\|_1^p \Big\}, \quad \forall i,j \in [n]. \bigg\}
\end{split}
\end{align}
\end{theorem}
\emph{Proof: }
Given a decision $z$,
the worst-case expected cost $\JP(z)$, as defined in \eqref{eq:strong-dual0-cp}, is equivalent to:
\begin{align}
\label{eq:strong-dual0-cd}
\JP(z) = \inf_{\lambda_1, \lambda_2 \geq 0} \paran{ \lambda_1 \varepsilon_1^p + \lambda_2 \varepsilon_2^p +
\hspace{-0.25cm}\sup_{\xi \in \Xi(\hat\mu_1, \hat\mu_2)}\hspace{-0.23cm}\E_{(Y_1, Y_2)\sim \xi} \Paran{\sup_{y \in \Yscr} \big\{c(z, y) - \lambda_1 \|Y_1 - y\|_1^p - \lambda_2 \|Y_2 - y\|_1^p \big\}}}.
\end{align}
This result can be derived analogously to the strong duality in the single Wasserstein-ball DRO case (e.g., \citealp{zhang2022simple}).
We provide the detailed proof in Appendix \ref{app:IWDRO-formulation}.

The optimization program \eqref{eq:strong-dual0-cd} is generally intractable due to the two nested supremum operations within the minimization problem.
Following the argument used in \eqref{eq:transportation-1}
that $\hat\mu_1$ and $\hat\mu_2$ have finite support sets,
we represent the joint probability for each pair of masses $(y_{1,i}, y_{2,j})$ with $p_{i,j}$,
and denote $\displaystyle{\ell_{i,j} = \sup_{y \in \Yscr} \big\{c(z, y) - \lambda_1 \|y_{1,i} - y\|_1^p - \lambda_2 \|y_{2,j} - y\|_1^p \big\}}$.
Then, the inner supremum problem over $\Xi(\hat\mu_1, \hat\mu_2)$ inside \eqref{eq:strong-dual0-cd}
can be restructured as a transportation problem as follows:
\begin{align}\label{eq:reformulation-inner}
\max_{p_{i,j} \ge 0} \Bigg\{ \sum_{i,j\in [n]}p_{i,j} \ell_{i,j}, \quad \text{s.t.  }  \sum_{j\in[n]}p_{i,j} = w_{1,i}, ~\forall i\in[n];   ~\sum_{i\in[n]}p_{i,j} = w_{2,j}, ~\forall j\in[n].  \Bigg\}
\end{align}
We formulate the dual program of \eqref{eq:reformulation-inner} with dual variables $u_i$ and $v_j$:
\begin{align}\label{eq:reformulation-inner-dual}
\min_{u_i, v_j} \Bigg\{ \sum_{i \in [n]} w_{1,i}u_i + \sum_{j \in [n]} w_{2,j}v_j, \quad \text{s.t.  } u_i + v_j \ge l_{i,j}, ~\forall i,j\in[n].  \Bigg\}
\end{align}
Then, we inject the dual formulation \eqref{eq:reformulation-inner-dual} into \eqref{eq:strong-dual0-cd} to finally obtain $\JP(z) = \JD(z)$,
where $\JD(z)$ is a single layer infimum program as defined in \eqref{eq:strong-dual-further}.
\hfill \halmos

In the derivation,
we focus on addressing the inner supremum problem over $\Xi(\hat\mu_1, \hat\mu_2)$ in \eqref{eq:strong-dual0-cd}
by introducing an additional layer of duality.
This approach substantiates the use of the intersection ambiguity set
$\setaiw$, making a novel contribution to the Wasserstein-DRO literature.

\Cref{thm:strong-duality} enables the reformulation of \eqref{eq:iwdro-opt} for tractable solutions.
With a piecewise concave cost function,
we can transform \eqref{eq:iwdro-opt} into a convex program by simplifying the robust constraints in \eqref{eq:strong-dual-further}.
Furthermore, when the cost function is piecewise linear,
\eqref{eq:iwdro-opt} can be further simplified to a linear program.
These results are presented in the following corollaries.

\begin{corollary}[Convex Reformulation]\label{coro:cvx}
If $\displaystyle{c(z,y) = \max_{s \in [S]} \{c_s(z, y)}\}$
and each $c_s(z, y)$ is concave and upper semi-continuous in $y \in \Yscr$,
then \eqref{eq:iwdro-opt} can be reformulated as a convex program.
\end{corollary}
\emph{Proof: }
Consider the case of $p = 1$.
It holds that:
\begin{align}
\label{eq:reformulation-cvx-obj}
\begin{split}
& \JP(z) =\inf ~\lambda_1 \varepsilon_1 +\lambda_2 \varepsilon_2 + \sum_{i \in[n]}w_{1,i} u_i + \sum_{j \in [n]} w_{2,j} v_j  \\
\text{\textup{s.t.}} ~~ & u_i + v_j \geq [-c_s(z,\cdot)]^*(\alpha_{i,s} + \beta_{j,s} - \zeta_{i,j,s}) + \sigma_{\Yscr}(\zeta_{i,j,s}) - \alpha_{i,s}^{\top}y_{1,i} - \beta_{j,s}^{\top}y_{2,j}, ~~ \forall i,j \in [n], s \in [S]; \\
& \|\alpha_{i,s}\|_{\infty} \leq \lambda_1, ~~\|\beta_{i,s}\|_{\infty} \leq \lambda_2, ~~ \alpha_{i,s}, \beta_{i,s}  \in \R^{d_y}, \quad \forall i \in [n], s \in [S]; \\
& \lambda_1, \lambda_2 \geq 0;
\quad u_i, v_i \in \R, ~~\forall i \in [n];
\quad \zeta_{i,j,s} \in \R^{d_y}, ~~\forall i, j \in [n], s \in [S].
\end{split}
\end{align}

The optimization problem \eqref{eq:reformulation-cvx-obj} is a convex program
because both the conjugate of the component function $[-c_s(z,\cdot)]^*$ and the support function $\sigma_{\Yscr}$ are convex.
Hence, it is straightforward to verify that \eqref{eq:iwdro-opt},
which is equivalent to $\displaystyle{\min_{z \in \Zscr} \JP(z)}$, can be reformulated as a convex program.

The reformulation details are provided in Appendix~\ref{app:tractable-formulation}. Analogous reformulations can be performed for other values of $p$.
\hfill \halmos

\begin{corollary}[Linear Reformulation]\label{coro:linear}
If $\displaystyle{c(z,y)= \max_{s \in [S]}\{a_{s}(z)^{\top}y + b_{s}(z)\}}$,
where $a_{s}(z)$ and $b_{s}(z)$ are linear in $z$,
and the uncertainty region $\Yscr$ is a polyhedron,
then \eqref{eq:iwdro-opt} can be reformulated as a linear program when choosing $p = 1$.
\end{corollary}
\emph{Proof: }
Assume that $\Yscr = \{y \in \R^{d_y}: Ay \leq h \}$ is a polyhedron set, where $A \in \R^{d_m \times d_y}$ and $h \in \R^{d_m}$.
If $p = 1$ and $c(z,y)$ is piecewise linear, the first constraint in program~\eqref{eq:reformulation-cvx-obj} can be further simplified to linear ones,
and thus we obtain a linear program as following.
\begin{align}
\label{eq:reformulation-linear-obj}
\begin{split}
\JP(z) =\inf ~~  &\lambda_1 \varepsilon_1 +\lambda_2 \varepsilon_2 + \sum_{i \in[n]}w_{1,i} u_i + \sum_{j \in [n]} w_{2,j} v_j  \\
\text{\textup{s.t.}} \quad & u_i + v_j \geq b_s(z) + \gamma_{i,j,s}^{\top}h - \alpha_{i,s}^{\top}y_{1,i} - \beta_{j,s}^{\top} y_{2,j}, \quad \forall i,j \in [n], s \in [S]; \\
& \alpha_{i,s} + \beta_{j,s} = A^{\top}\gamma_{i,j,s} - a_{s}(z), \quad \forall i,j \in [n], s \in [S]; \\
& \|\alpha_{i,s}\|_{\infty} \leq \lambda_1, ~~\|\beta_{i,s}\|_{\infty} \leq \lambda_2, ~~ \alpha_{i,s}, \beta_{i,s}  \in \R^{d_y}, \quad \forall i \in [n], s \in [S]; \\
& \lambda_1, \lambda_2 \geq 0;
\quad u_i, v_i \in \R, ~~\forall i \in [n];
\quad \gamma_{i,j,s} \in \R^{d_m}_{+}, ~~\forall i, j \in [n], s \in [S].
\end{split}
\end{align}
The reformulation details are provided in Appendix~\ref{app:tractable-formulation}. \hfill \halmos

Note that the computational results presented in this section
can be extended to an intersection ambiguity set involving more than two Wasserstein balls, as well as to a generalized Wasserstein distance beyond \Cref{def:was-distance}.
We provide a detailed discussion in Appendix \ref{app:IWDRO-formulation}.

\subsection{Performance Guarantee}\label{subsec:iwdro-opt-guarantees}
After solving the optimization program \eqref{eq:iwdro-opt},
our goal is to evaluate the performance of its solutions under covariate shift.
Let $\ziw$ and $\Jiw$ represent the optimal solution and objective value of \eqref{eq:iwdro-opt}, respectively.
Recall that the conditional distribution of $Y$ given $x$ remains unchanged as $\mu_{Y|x}$ across both the training and test datasets.
We denote by $z^*$ the optimal solution to $\displaystyle{\min_{z \in \Zscr}\E_{\mu_{Y|x}}\left[ c(z,Y)\right]}$.
In contextual optimization, $z^*$ is considered the theoretical ideal decision, but it is not attainable because $\mu_{Y|x}$ is unknown.
In contrast, {\iwdro} provides a practical and implementable decision $\ziw$.
Assuming that $\setaiw$ includes $\mu_{Y|x}$,
we can establish the following performance guarantees by comparing $z^*$ and $\ziw$.
\begin{proposition}[Performance of {\iwdro}]\label{prop:optimality-gap-iwdro}
For each $x\in \Xscr$, if $\mu_{Y|x} \in \setaiw$, then we have
the cost bound \eqref{eq:optimality-iwdro-bound} and the optimality gap \eqref{eq:optimality-iwdro-gap} as follows:
\begin{align}
\label{eq:optimality-iwdro-bound}
& \E_{\mu_{Y|x}}[c(z^{*}, Y)] \le \E_{\mu_{Y|x}}[c(\ziw, Y)] \le \Jiw, \\
\label{eq:optimality-iwdro-gap}
& \E_{\mu_{Y|x}}[c(\ziw, Y)] - \E_{\mu_{Y|x}}[c(z^{*}, Y)] \le 2\mathcal{L}_{z^*} \min\{\varepsilon_1,\varepsilon_2\},
\end{align}
where $\mathcal{L}_{z^*}= \|c(z^*,\cdot)\|_{\textup{Lip}}$ is the local Lipschitz constant of $c(\cdot,\cdot)$ under $z^{*}$.
\end{proposition}
\emph{Proof:} For generality, we consider a single Wasserstein-ball DRO program as follows:
\begin{align}
\label{eq:single-ball-dro}
\min_{z \in \Zscr} \sup_{\mu \in \Ascr_{\hat\mu, \varepsilon}} ~\E_{Y\sim\mu}[c(z,Y)], ~~\text {with} ~~\Ascr_{\hat\mu, \varepsilon} \coloneqq\{\mu: \Wscr_p(\mu, \hat\mu)\leq \varepsilon\}.
\end{align}
Here, the ambiguity set $\Ascr_{\hat\mu, \varepsilon}$ is constructed using its center $\hat\mu$ and radius $\varepsilon$ (as indicated by the subscript).
Let $\hat z$ and $\hat J$ denote the optimal solution and objective value to \eqref{eq:single-ball-dro}, respectively.

Given that $\mu_{Y|x} \in \Ascr_{\hat\mu, \varepsilon}$ and the definitions of optimality for $z^{*}$ and $\hat z$, we obtain the following:
\begin{align}
\label{eq:single-ball-dro-gap-1}
\E_{\mu_{Y|x}}[c(z^{*}, Y)] \le \E_{\mu_{Y|x}}[c(\hat z, Y)] \le \hat J \coloneqq \sup_{\mu \in \Ascr_{\hat\mu, \varepsilon}} \E_{Y\sim\mu}[c(\hat z, Y)]
\le \sup_{\mu\in \Ascr_{\hat\mu, \varepsilon}}\E_{Y\sim\mu}[c(z^*, Y)].
\end{align}
This directly leads to $\text{Term-1} := \E_{\mu_{Y|x}}[c(\hat z, Y)] - \sup_{\mu\in \Ascr_{\hat\mu, \varepsilon}}\E_{Y\sim\mu}[c(z^*, Y)] \le 0$.

Additionally, by the variational representation of 1-Wasserstein distance, we have:
\begin{align}
\E_{\mu}[c(z^*, Y)] - \E_{\mu_{Y|x}}[c(z^*,Y)]
\le \sup_{\|\tilde{c}\|_{\textup{Lip}} \le \mathcal{L}_{z^*}} \Big\{\E_{\mu}[\tilde{c}(Y)]- \E_{\mu_{Y|x}}[\tilde{c}(Y)] \Big\} = \mathcal{L}_{z^*} \Wscr_1(\mu, \mu_{Y|x}), \label{eq:single-ball-dro-gap-2-term2-1}
\end{align}
where the supremum is taken over all possible functions $\tilde{c}$ satisfying $\|\tilde{c}\|_{\textup{Lip}} \le \mathcal{L}_{z^*}$, including $c(z^*, \cdot)$.
Following from \eqref{eq:single-ball-dro-gap-2-term2-1},
since $\Wscr_1(\mu, \mu_{Y|x}) \le \Wscr_p(\mu, \mu_{Y|x}) \le \Wscr_p(\mu, \hat\mu) + \Wscr_p(\hat\mu, \mu_{Y|x}) \le 2\varepsilon$,
we obtain:
\begin{align}
\text{Term-2} := \sup_{\mu\in \Ascr_{\hat\mu, \varepsilon}}\E_{Y\sim\mu}[c(z^*, Y)] - \E_{\mu_{Y|x}}[c(z^*,Y)]
\le \mathcal{L}_{z^*} \Wscr_1(\mu, \mu_{Y|x})
\le 2\mathcal{L}_{z^*} \varepsilon. \nonumber
\end{align}

By combing these two terms,
we establish the following optimality gap:
\begin{align}
\label{eq:single-ball-dro-gap-2}
\E_{\mu_{Y|x}}[c(\hat z, Y)] - \E_{\mu_{Y|x}}[c(z^{*}, Y)]
= \text{Term-1} + \text{Term-2} \le 2\mathcal{L}_{z^*} \varepsilon.
\end{align}

Because $\setaiw = \{\mu: \Wscr_p(\mu, \hat\mu_1)\leq \varepsilon_1, \Wscr_p(\mu, \hat\mu_2)\leq \varepsilon_2\}$
is the intersection of two Wasserstein balls $\Ascr_{\hat\mu_1, \varepsilon_1}$ and $\Ascr_{\hat\mu_2, \varepsilon_2}$,
$\mu_{Y|x} \in \setaiw$ implies that $\mu_{Y|x} \in \Ascr_{\hat\mu_1, \varepsilon_1}$ and $\mu_{Y|x} \in \Ascr_{\hat\mu_2, \varepsilon_2}$.
By extending derivations above,
we can prove \eqref{eq:optimality-iwdro-bound} and \eqref{eq:optimality-iwdro-gap} for $\ziw$ in a similar manner.
\hfill \halmos

\Cref{prop:optimality-gap-iwdro} shows that
{\iwdro} provides an upper bound \eqref{eq:optimality-iwdro-bound} for the  contextual optimization problem.
Furthermore, in the proof of \eqref{eq:single-ball-dro-gap-2}, we observe that the performance gap between the ideal decision $z^*$ and a decision based on a DRO model
is proportionally bounded by the radius of the ambiguity set.
This finding leads to the optimality gap \eqref{eq:optimality-iwdro-gap}:
by intersecting two Wasserstein balls,
we can achieve an ambiguity set $\setaiw$ with a reduced radius,
which could, in turn, enhance the cost performance of {\iwdro} compared to the corresponding single Wasserstein-ball DRO models.

\section{Statistical Guarantees of the Ambiguity Set}\label{sec:statistical-result}
We have demonstrated that the effectiveness of {\iwdro} critically depends on the specification of the ambiguity set $\setaiw$.
As shown in \Cref{prop:optimality-gap-iwdro},
the essential condition for assuring performance guarantees is that $\setaiw$ includes the true conditional distribution $\mu_{Y|x}$.
However, an excessively large $\setaiw$ can lead to overly conservative decisions with a larger optimality gap.
In this section, we discuss the construction of $\setaiw$ and examine its statistical properties in terms of containing $\mu_{Y|x}$.

\subsection{Main Results}\label{subsec:statistical-guarantees}
In selecting $\munp$ and $\mup$, i.e., the centers of the two Wasserstein balls as defined in \eqref{eq:ambiguity-set},
we employ the two prevalent estimation approaches in the Wasserstein-DRO literature: nonparametric and parametric estimates.
These methods determine each reference distribution using an estimated conditional distribution
derived from the training samples $\{(x_i, y_i)\}_{i = 1}^n$ and the test covariate $x$.
\begin{itemize}
\item A nonparametric approach typically assigns a weight to each sample outcome $y_i$
based on the distance between $x_i$ and the test covariate $x$.
We define $\munp$ using the NW kernel estimator as
\begin{align}
\label{eq:ball-center-np}
\munp \coloneqq\sum_{i = 1}^n\frac{K((x - x_i)/h_n)}{\sum_{j = 1}^n K((x-x_j)/h_n)}\delta_{y_i},
\end{align}
where $K(\cdot)$ is a kernel function and $h_n$ is its bandwidth parameter.

\item A parametric approach utilizes samples $\{(x_i, y_i)\}_{i = 1}^n$ to fit a regression model $f_{\theta}$
(a parameterized function which predicts the outcome $y_i$ based on $x_i$),
and estimates the model parameter, denoted by $\hat{\theta}$.
Following \cite{kannan2020residuals},
we perturb the predicted value $f_{\hat{\theta}}(x)$ for the test covariate $x$
by adding the residual $y_i - f_{\hat{\theta}}(x_i)$ from the $i$-th sample,
and then construct $\mup$ as
a residual-based distribution with uniform weights $1/n$:
\begin{align}
\label{eq:ball-center-p}
\mup \coloneqq \frac{1}{n}\sum_{i = 1}^n \delta_{\{f_{\hat{\theta}}(x) + y_i - f_{\hat{\theta}}(x_i)\}}.
\end{align}
\end{itemize}
We refer to the balls centered at $\munp$ and $\mup$ as the nonparametric and parametric balls, respectively.
Note that both $\munp$ and $\mup$ yield discrete distributions that satisfy Assumptions~\ref{asp:dist_structure} and \ref{asp:reg-costfunc},
and thus the convex reformulation to \eqref{eq:iwdro-opt} derived in Section \ref{sec:IWDRO-formulation} applies.

To describe covariate shift,
we denote the probability densities of the covariate $X$ in the training and test data distributions by $\mu_X(x)$ and $\nu_X(x)$, respectively.
We then use the Radon-Nikodym derivative $\rho_{x} \coloneqq {\nu_X(x)}/{\mu_X(x)}$ to quantify the degree of covariate shift.
For example, a large value of $\rho_{x}$ for a test covariate $x$ indicates that there are few nearby training samples.
For analytical purposes, we consider cases where the covariate distributions are continuous
and impose some standard technical assumptions on the conditional distribution.

\begin{assumption}[Continuity of Covariate Distribution]\label{asp:regular_marginal_dist}
$\mu_X(x)$ is Lipschitz continuous.
\end{assumption}

\begin{assumption}[Regularity of the Conditional Distribution]\label{asp:regular_dist}
\hspace{-0.2cm}The conditional distribution
follows $\mu_{Y|x} \sim y=f(x) + \eta$, where:
(i) the function $f$ satisfies the $(L, \beta)$-H\"older condition under the $\ell_2$-norm for constants $L > 0$ and $\beta \in (0, 1]$,
i.e., $\|f(x_1) - f(x_2)\|_1\leq L\|x_1 - x_2\|_2^{\beta}, ~\forall x_1, x_2 \in \Xscr$;
and (ii) the random noise vector $\eta \in \R^{d_y}$ is independent of $x$,
has zero means $\E[\eta] = \mathbf{0}$ and $\E[\eta \eta^{\top}]\preceq \sigma^2 \mathbf{I}$,
and is light-tailed distributed with $\E[\exp(\|\eta\|_2^{p + 1})] < \infty$.
\end{assumption}

Assumption~\ref{asp:regular_dist} specifies the data generating process $y = f(x) + \eta$ that links the uncertainty outcome $y$ to the covariate $x$.
In particular,
the H\"older condition on $f$ is commonly adopted for the analysis of nonparametric estimates \citep{gyorfi2006distribution,tsybakov2008introduction},
while the light-tailed property of $\eta$ is a typical prerequisite for bounding the Wasserstein distance \citep{fournier2015rate}.
Under these conditions, along with \Cref{asp:kernel} (introduced later for kernel functions),
we can establish high-probability guarantees for the inclusion of $\mu_{Y|x}$ in $\setaiw$.

\begin{theorem}[Coverage Guarantees of $\setaiw$]\label{thm:choice-size}
Under Assumptions \ref{asp:regular_marginal_dist}, \ref{asp:regular_dist}, and \ref{asp:kernel},
for any $x \in \Xscr$ and a risk level $\alpha \in (0, 1)$,
we have $\P(\mu_{Y|x} \in \setaiw) \geq 1- \alpha,$
if we choose
\begin{align}
\label{eq:ball-radius}
\varepsilonnp = \Theta \bigg( \bigg(\frac{\rho_x}{n}\bigg)^{\rnp}\bigg)
~~\text{and}~~
\varepsilonp = \Theta \bigg( \bigg(\frac{1}{n}\bigg)^{\rp} \bigg) + \eapx.
\end{align}
Here, 
$\rnp \in (0,1)$ and $\rp \in (0,1)$ are some constants,
and $\eapx$ denotes the model misspecification error associated with the parametric estimator $\mup$.
\end{theorem}

To highlight the relationship between the necessary radius and the degree of covariate shift,
we focus on the dependence on $\rho_{x}$ in expression \eqref{eq:ball-radius},
where we use $\rnp$ and $\rp$ to abstract away certain constant terms
and use $\Theta$-notation to suppresses the polynomial dependence of $\log(1/\alpha)$, $L$, $p$, and other parameters.
\Cref{thm:choice-size} shows that, to include $\mu_{Y|x}$,
the radius $\varepsilonnp$ for the nonparametric ball must increase with $\rho_{x}$,
whereas the radius $\varepsilonp$ for the parametric ball is influenced not by $\rho_{x}$ but by a model mispecification error $\eapx$.
This finding explains the benefits of combining different estimators in $\setaiw$
by comparing $\varepsilonnp$ and $\varepsilonp$.
Specifically, when $\rho_{x}$ is small and there are adequate nearby samples,
the nonparametric ball can have a smaller radius $\varepsilonnp$.
Conversely, under a significant covariate shift indicated by a large $\rho_{x}$,
the parametric ball requires a smaller radius $\varepsilonp$.
By intersecting the two balls, we achieve a more compact ambiguity set $\setaiw$ that still effectively contains $\mu_{Y|x}$ for all $x \in \Xscr$.
In practice, although the exact values of $\rho_{x}$ and $\eapx$ are typically unobservable,
we can estimate their scales or upper bounds by analyzing the test covariates and training samples,
guiding us in selecting appropriate values for $\varepsilonnp$ and $\varepsilonp$.
In the remainder of this section, we derive \Cref{thm:choice-size} by analyzing the measure concentration of $\munp$ and $\mup$, respectively.

\subsection{Measure Concentration of the Nonparametric Estimate}\label{subsec:estimator-np}
For the nonparametric ball center $\munp$, as defined in \eqref{eq:ball-center-np},
we employ the NW kernel estimator to assign a weight $\frac{K((x - x_i)/h_n)}{\sum_{j = 1}^n K((x-x_j)/h_n)}$ to each sample outcome $y_i$.
Here, $K:\R^{d_x} \rightarrow [0,+\infty)$ is a kernel function, and $h_n\in \R^+$ denotes its bandwidth parameter which depends on the sample size $n$.
Typically, $K(\tau)$ diminishes as $\|\tau\|_2$ increases,
implying that $K((x-x_i)/h_n)$ yields larger values for $x_i$ closer to $x$.
As a result, $\mup$ assigns greater weights to the corresponding $y_i$.
For analysis purposes, we impose some standard technical conditions on the kernel function
to ensure that the value of $K(\tau)$ is bounded and decreases to zero faster than any polynomial rate when $\|\tau\|_2 \to \infty$.

\begin{assumption}[Regularity \hspace{-0.05cm}of \hspace{-0.05cm}the \hspace{-0.05cm}Kernel \hspace{-0.05cm}Function]\label{asp:kernel}
\hspace{-0.2cm}There exist constants $r > 0$, \hspace{-0.05cm}$0< b_r \hspace{-0.05cm}\le b_R$,
and $C_K > 0$
such that, for all $\tau \in \R^{d_x}$, the kernel function $K(\tau)$ satisfies:
(i) $b_r \I{\| \tau \|_2\le r} \le K(\tau) \le b_R$;
and (ii) $ K(\tau) < C_K \exp(-\|\tau\|_2)$.
\end{assumption}
\Cref{asp:kernel} is satisfied for most kernels with bounded supports and Gaussian kernels with lighter tails than polynomial decreasing rates. For example, the common forms of $K(\tau)$ include $\I{\|\tau\|_2\le 1}$, $(1-\|\tau\|_2^2)_+$, and $\exp(-\|\tau\|_2^2)$,
which all satisfy \Cref{asp:kernel}.

Existing theoretical studies on the NW kernel estimator (e.g., \citealp{gyorfi2006distribution, tsybakov2008introduction})
have primarily focused on its statistical properties for estimating the conditional mean $\E[Y|x]$.
In contrast, our work is dedicated to adapting it to estimate the conditional distribution $\mu_{Y|x}$, making technical contributions to the analysis of finite-sample measure concentration.
To the best of our knowledge, this is the first result that establishes a nonasymptotic rate for the NW kernel estimator using the Wasserstein distance.
(In Appendix~\ref{app:kNN}, we extend our measure concentration analysis to the nearest neighbors, another classical nonparametric estimator.)



\begin{theorem}[Measure \hspace{-0.05cm}Concentration of $\munp$]\label{thm:concentration-np}
\hspace{-0.2cm}Under Assumptions \ref{asp:regular_marginal_dist}, \ref{asp:regular_dist}, and \ref{asp:kernel},
for any $ x\in \Xscr$ and $t\in (0,1]$,
if we choose $h_n = C_h n^{-\psi}$  with $\psi \in (0, d_x^{-1})$, and $n \ge C_{\mu} \mu_X(x)^{-\frac{1}{\psi d_x}}$,
then we have:
\begin{equation}\label{eq:convergence_rate_np}
\hspace{-0.15cm}\P\big(\Wscr_p(\mu_{Y|x},\munp)> t + C_0 L h_n^{\beta} \big) \le
\begin{cases} C_1 \exp(-C_2n \mu_X(x) h_n^{d_x} t^2), & \text{if } p > d_y / 2, \\
C_1 \exp(-C_2n\mu_X(x) h_n^{d_x} (t/\log (2 + 1/t))^2), & \text{if } p = d_y / 2, \\
C_1 \exp(-C_2n \mu_X(x) h_n^{d_x} t^{d_y/p}), & \text{if } p \in [1, d_y / 2),\\
\end{cases}
\end{equation}
where $C_0$, $C_1$, $C_2$, $C_h$ and $C_{\mu}$ are some constants independent of $x$, $n$, and $t$.
\end{theorem}
\emph{Proof:}
Here, we present a brief proof sketch.
The detailed proof is provided in Appendix~\ref{app:proof-concentration}.
To upper bound $\Wscr_p(\mu_{Y|x},\munp)$, we construct a discrete distribution:
\begin{align}
\mur  = \sum_{i = 1}^n\frac{K((x - x_i)/h_n)}{\sum_{j = 1}^n K((x-x_j)/h_n)}\delta_{\{f(x)+\eta_i\}}, \nonumber
\end{align}
which shares the same weights as $\munp$ and has supports $\{f(x) + \eta_i\}_{i \in [n]}$.
Using $\mur$ as an intermediate to apply the triangle inequality $\Wscr_p(\mu_{Y|x},\munp) \leq \Wscr_p(\mu_{Y|x}, \mur) + \Wscr_p(\mur, \munp)$,
we derive the probabilistic bound
$\P \big(\Wscr_p(\mu_{Y|x},\munp)>t + C_0 Lh_n^{\beta}\big) \le \P\big(\Wscr_p(\mu_{Y|x},\mur)>t\big)+ \P \big(\Wscr_p(\mur,\munp)> C_0 Lh_n^{\beta}\big).$

The introduction of $\mur$ simplifies the analysis.
When analyzing $\Wscr_p(\mu_{Y|x},\mur)$, the support points of $\mur$, i.e., $\{f(x) + \eta_i\}_{i \in [n]}$, can be viewed as i.i.d. samples from $\mu_{Y|x}$.
This allows us to apply an extended result from \cite{fournier2015rate} to derive an upper bound.

On the other hand, when analyzing $\Wscr_p(\mur, \munp)$,
we note that $\mur$ and $\munp$ share the same probability mass,
and the difference between their supports is given by $f(x) + \eta_i - y_i = f(x) - f(x_i)$.
Using the definition of the Wasserstein distance and the H\"older continuity on $f(x)$ from \Cref{asp:regular_dist},
we can bound this term accordingly.

In the analysis, the condition $n \ge C_{\mu} \mu_X(x)^{-\frac{1}{\psi d_x}}$ ensures a sufficient number of nearby samples for $\munp$ to achieve the desired concentration rate \eqref{eq:convergence_rate_np}.
\hfill \halmos

\Cref{thm:concentration-np} establishes a link between the efficacy of $\munp$ and the training sample density $\mu_X(x)$.
As indicated by the concentration rate \eqref{eq:convergence_rate_np},
a low value $\mu_X(x)$ corresponds to a high probability of observing a large distance $\Wscr_p(\mu_{Y|x},\munp)$.
This implies that, in covariate-shifted regions where training samples are sparse,
the nonparametric estimator $\munp$ is more likely to deviate from the true conditional distribution $\mu_{Y|x}$.
Based on this result, we can determine the necessary radius $\varepsilonnp$ for the nonparametric ball to include $\mu_{Y|x}$, as stated in \Cref{thm:choice-size}.

\noindent\emph{Derivation of $\varepsilonnp$ for \Cref{thm:choice-size}:\label{proof:radius_np} }
For the first case of \eqref{eq:convergence_rate_np} where $p > d_y / 2$, it holds that:
\begin{equation}\label{eq:radius_np_step1}
\P\big(\Wscr_p(\mu_{Y|x},\munp)> t + C_0 L h_n^{\beta} \big)
\le C_1 \exp\Para{-C_2 n \mu_X(x) h_n^{d_x} t^{2}}.
\end{equation}
Given a risk level $\alpha \in (0, 1)$,
we set the right side of \eqref{eq:radius_np_step1} equal to $\frac{\alpha}{2}$ and solve for $t=\Para{\frac{\log(2C_1/\alpha)}{C_2 n \mu_X(x) h_n^{d_x}}}^{\frac{1}{2}}$.
Consequently, we establish that $\P\big(\Wscr_p(\mu_{Y|x}, \munp) > \varepsilonnp \big) \le \frac{\alpha}{2}$,
with the radius chosen as:
\begin{equation}\label{eq:radius_np_step2}
\varepsilonnp = \Para{\frac{\log(2C_1/\alpha)}{C_2 n \mu_X(x) h_n^{d_x}}}^{\frac{1}{2}} + C_0 L h_n^{\beta}.
\end{equation}
We minimize $\varepsilonnp$, as given in \eqref{eq:radius_np_step2}, with respect to $n\mu_X(x)$
by choosing $h_n = \Theta\big((n\mu_X(x))^{-\frac{1}{2\beta + d_x}}\big)$,
and we also substitute $\rho_x = \frac{\nu_X(x)}{\mu_X(x)} = \Theta\Big(\frac{1}{\mu_X(x)}\Big)$.
Then, we obtain the necessary radius:
\begin{equation}\label{eq:radius_np_step_final}
\varepsilonnp = \Theta\Big(\max\{\log(1/\alpha), L\} \left(\frac{\rho_x}{n}\right)^{\rnp}\Big)
\text{ with } \rnp = \frac{\beta}{2\beta + d_x}.
\end{equation}
This leads to $\varepsilonnp = \Theta \big( \big(\frac{\rho_x}{n}\big)^{\rnp}\big)$, as presented in \eqref{eq:ball-radius}.
(Here, for $h_n$ and $\varepsilonnp$,
we focus on their dependence on $n\mu_X(x)$ (or equivalently $\rho_x/n$) and use $\Theta$-notation to suppress other constant parameters.
Aligning with the setting $h_n = C_h n^{-\psi}$ in \Cref{thm:concentration-np}, we set $\psi=\frac{1}{2\beta + d_x}$ and omit $C_h$ for simplicity.)

The analogous procedures apply for the two remaining cases of \eqref{eq:convergence_rate_np}.
We can prove that \eqref{eq:radius_np_step_final} holds for $p = d_y / 2$ and $ p \in [1, d_y / 2)$,
by choosing $\rnp = \frac{\beta}{2\beta + d_x}$
and $\rnp = \frac{\beta}{\beta d_y/p + d_x}$, respectively.
\hfill \halmos

The concentration rate of $\munp$ and the resulting size of $\varepsilonnp$, as presented in \eqref{eq:convergence_rate_np} and \eqref{eq:radius_np_step_final}, respectively, are proven to be tight.
This can be demonstrated with the following example.

\begin{example}\label{ex:tightness}
Consider a data generating process $y = f(x) + \eta$ that satisfies \Cref{asp:regular_dist},
with the covariate domain $\Xscr = [0, 1]^{d_x}$.
We let $\hat f(x) = \E_{\munp}[Y] = \sum_{i=1}^{n}\frac{K((x - x_i)/h_n)}{\sum_{j=1}^{n} K((x - x_j)/h_n)} y_i$ denote the mean of $\munp$,
and also note that $\E_{\mu_{Y|x}}[Y] = f(x)$.
Since the Wasserstein distance between two distributions is lower bounded by the difference of their means, we have:
\begin{equation}\label{eq:example2-lower-bound}
\Wscr_p(\mu_{Y|x}, \munp) \geq \Wscr_1(\mu_{Y|x}, \munp) \geq \|f(x) - \hat{f}(x)\|_1.
\end{equation}

Let $\Fscr$ denote the set of all possible functions $f$ that satisfy \Cref{asp:regular_dist}, specifically the $(L, \beta)$-H\"older condition under the $\ell_2$-norm. 
Given a test covariate $x$, let $\tilde{f}(x, \{(x_i, y_i)\}_{i = 1}^{n})$ represent an estimate of the outcome $y$ based on the training samples.
Then, by Theorem 3.2 in \cite{gyorfi2006distribution}, we have:
\begin{equation}\label{eq:example2-lower-bound2}
\inf_{\tilde{f}} \sup_{f \in \Fscr} \E\big[\|f(x) - \tilde{f}(x, \{(x_i, y_i)\}_{i = 1}^{n})\|_2\big] = \Theta(n^{-\frac{\beta}{2\beta + d_x}}),
\end{equation}
where the expectation is taken over the randomness of training samples,
and the infimum is taken over all possible estimates $\tilde{f}$, i.e., over all measurable functions of the training data, including $\hat{f}(x)$.

Finally, by substituting $\rho_x = \Theta(\frac{1}{\mu_X(x)}) = \Theta(1)$ and comparing~\eqref{eq:example2-lower-bound} with~\eqref{eq:example2-lower-bound2}, we can show that for any $\munp$, there exists a data generating process $y = f(x) + \eta$, where $f(x) \in \Fscr$, such that $\E[\Wscr_p(\mu_{Y|x}, \munp)] = \Omega\big(\Para{\frac{\rho_x}{n}}^{\frac{\beta}{2\beta + d_x}}\big)$,
which matches the order of $\varepsilonnp$, as pretested in \eqref{eq:radius_np_step_final}.
\hfill \halmos
\end{example}

Additionally, we note that \Cref{thm:concentration-np} can remain valid under more general conditions,
including scenarios with covariate-dependent noise.
This is further elaborated in the following remark.
\begin{remark}[Relaxation \hspace{-0.1cm}of \hspace{-0.1cm}\Cref{asp:regular_dist}]
\hspace{-0.1cm}Consider a data generating process $y = f(x) + \eta_x$,
where the noise $\eta_x$ dependents on $x$ (e.g., the variance is a function of $x$).
By imposing an additional condition:
$\Wscr_p(\eta_{x_1}, \eta_{x_2}) \leq L_1\|x_1 - x_2\|_2$ for some constant $L_1 > 0$ and for all $x_1, x_2 \in \Xscr$,
we can show that the measure concentration results in \Cref{thm:concentration-np} still hold.
This setting is less restrictive than the i.i.d noise condition posited in \Cref{asp:regular_dist},
distinguishing our result from \cite{kannan2020residuals}, which requires independent noise.
\hfill \halmos
\end{remark}

\subsection{Measure Concentration of the Parametric Estimator}\label{subsec:estimator-p}
For the parametric ball center $\mup$, as defined in \eqref{eq:ball-center-p},
we employ a regression model to estimate the conditional distribution.
Unlike nonparametric estimates which rely on local samples,
parametric regression performs a global estimation.
As a result, $\mup$ may exhibit greater resilience to the scarcity of nearby samples under covariate shifts.
Nevertheless, it remains susceptible to errors from model misspecification.
To elaborate, recall that
\Cref{asp:regular_dist} specifies $\mu_{Y|x}$ through the data generating process $y = f(x) + \eta$.
In the regression estimate, since the true model $f$ is typically unknown, we select an approximate model $f_{\theta}$ and fit the parameter $\theta$.
Let $\theta^*$ denote the limiting parameter.
Consequently, $f_{\theta^*}$ represents the best possible outcome of such a parametric estimate.
The model misspecification error can then be quantified by a constant $\eapx$,
which uniformly bounds the difference between $f_{\theta}(x)$ and  $f_{\theta^*}(x)$ across the covariate domain.
The application of parametric regression methods in contextual optimization has been thoroughly studied in \cite{kannan2020residuals},
which provides a counterpart result to \Cref{thm:concentration-np} as follows.
\begin{theorem}[Measure \hspace{0.1cm}Concentration \hspace{0.1cm}of \hspace{0.1cm}$\mup$]\label{thm:concentration-p}
\hspace{0.1cm}Under Assumption \ref{asp:regular_dist},
for any $ x\in \Xscr$ and $t\in (0,1]$,
there exists a parametric regression estimator $\mup$ such that
\begin{equation}\label{eq:convergence_rate_p}
\P\big(\Wscr_p(\mu_{Y|x}, \mup) \geq t + \eapx\big) \leq
\begin{cases} C_4 \exp(-C_5n t^2), & \text{if } p > d_y / 2, \\
C_4 \exp(-C_5 (t/\log (2 + 1/t))^2), & \text{if } p = d_y / 2, \\
C_4 \exp(-C_5 n t^{d_y/p}), & \text{if } p \in [1, d_y / 2),\\
\end{cases}
\end{equation}
where $C_4$ and $C_5$ are some constants independent of $x$, $n$, and $t$.
\end{theorem}

\Cref{thm:concentration-p} shows that the convergence rate in \eqref{eq:convergence_rate_p} is independent of the density of nearby samples.
This highlights the advantage of $\mup$ in yielding more reliable estimates under covariate shift, provided that the regression model is properly specified.
The regression models and the requisite conditions ensuring that $\mup$ and $\eapx$ satisfy \Cref{thm:concentration-p}
are discussed in detail in \cite{kannan2020residuals} and \cite{iyengar2023hedging}.
In the following, we illustrate these concepts through the application of an OLS regression model.

\begin{example}\label{prop:p-estimator-example}
Consider a data generating process $y = f(x) + \eta$.
In an OLS regression, we assume a linear model $f_{\theta} = \theta^{\top} x$
and define the limiting parameter as $\theta^* \in \argmin_{\theta}\E[\|Y - \theta^{\top} X\|_2^2]$,
where the expectation is taken over the joint distribution of $(X,Y)$.
Thus, $f_{\theta^*} = \theta^{*\top}x$ represents the best linear approximation to $f(x)$,
and the unmodeled nonlinear term is defined as $g(x) = f(x) - \theta^{*\top} x$.
Within this setting, the model misspecification error can be quantified as $\eapx = 2 \sup_{x \in \Xscr} \|g(x)\|_1$.
In {\iwdro} , the regression parameter $\hat\theta \in \argmin_{\theta}\frac{1}{n}\sum_{i = 1}^n\|y_i - \theta^{\top} x_i\|_2^2$ is fitted using training samples,
and then the distribution estimator $\mup =\frac{1}{n}\sum_{i = 1}^n\delta_{\{\hat\theta^{\top} x + y_i - \hat\theta^{\top} x_i\}}$ is constructed based on $\hat\theta$.
Finally, $\mup$ and $\eapx$ together ensure the convergence rate given in \eqref{eq:convergence_rate_p},
under some additional conditions outlined in Assumption~\ref{asp:add-asp} in Appendix~\ref{app:p-estimator-proof}.
\end{example}


Based on \Cref{thm:concentration-p}, we can determine the necessary radius $\varepsilonp$ for the parametric ball to include $\mu_{Y|x}$, as stated in \Cref{thm:choice-size}

\noindent
\emph{Derivation of $\varepsilonp$ for \Cref{thm:choice-size}:\label{proof:radius_p} }
We follow the derivation of $\varepsilonnp$.
For the first case of \eqref{eq:convergence_rate_p} where $p > d_y / 2$,
we have $\P\big(\Wscr_p(\mu_{Y|x}, \mup) > \varepsilonp \big) \le \frac{\alpha}{2}$ with
\begin{equation}\label{eq:radius_p_step1}
\varepsilonp =  \Para{\frac{\log(2C_4/\alpha)}{C_5 n}}^{\frac{1}{2}} + \eapx.
\end{equation}
Then, we obtain the necessary radius $\varepsilonp = \Theta \big( \big(\frac{1}{n}\big)^{\rp}\big) + \eapx$, as presented in \eqref{eq:ball-radius}, with $\rp = \frac{1}{2}$.
Similarily, we can show that \eqref{eq:ball-radius} holds for $p = d_y / 2$ and $ p \in [1, d_y / 2)$,
by choosing $\rnp = \frac{1}{2}$ and $\rnp = \frac{p}{d_y}$, respectively.
\hfill \halmos

\smallskip
Finally, given that $\P\big(\Wscr_p(\mu_{Y|x}, \mup) > \varepsilonnp \big) \le \frac{\alpha}{2}$  and $\P\big(\Wscr_p(\mu_{Y|x}, \mup) > \varepsilonp \big) \le \frac{\alpha}{2}$ hold for $\varepsilonnp$ and $\varepsilonp$ chosen in \eqref{eq:ball-radius},
we prove \Cref{thm:choice-size} that
the intersection ambiguity set $\setaiw$ contains the true conditional distribution $\mu_{Y|x}$ with the following guarantee:
$$\P(\mu_{Y|x} \in \setaiw) \ge 1- \P\big(\Wscr_p(\mu_{Y|x}, \munp) > \varepsilonp \big) - \P\big(\Wscr_p(\mu_{Y|x}, \mup) > \varepsilonp \big) = 1-  \alpha.$$

\section{Approximation to {\iwdro}}\label{sec:approximation-formulation}

\subsection{An Approximation Model}\label{sec:approximation-formulation-model}
We now examine the computational complexity of the solution approach to \eqref{eq:iwdro-opt} presented in Section \ref{sec:IWDRO-formulation}.
Recall that when the cost function is piecewise concave, \eqref{eq:iwdro-opt} can be transformed into a tractable convex program.
Nevertheless, the reformulated program requires $\Theta(n^2)$ decision variables and constraints
(see the convex program \eqref{eq:reformulation-cvx-obj} in \Cref{coro:cvx}
and the linear program \eqref{eq:reformulation-linear-obj} in \Cref{coro:linear}).
In contrast, a conventional single Wasserstein-ball DRO program involves only $\Theta(n)$ decision variables and constraints.
To improve computational efficiency for large training datasets (i.e., when $n$ is large),
we introduce a single Wasserstein-ball DRO model, denoted as {\iwdroapx}, as an approximation to {\iwdro}.

Following the observation that the intersection set combines the advantages of different estimators,
it is natural to consider constructing a distribution ball centered around a mixture estimator.
Specifically, we define an ambiguity set
$\setame \coloneqq \{\mu: \Wscr_p(\mu, \mume) \leq \varepsilonme\}$,
where the center $\mume$ and the radius $\varepsilonme$ are determined by interpolating between the nonparametric and parametric approaches:
\begin{align}
\label{eq:mixture-estimator}
\mume  & \coloneqq \kappa_x \munp + (1-\kappa_x) \mup, \\
\label{eq:mixture-radius}
\varepsilonme & \coloneqq \big(\kappa_x \varepsilonnp^p + (1-\kappa_x) \varepsilonp^p\big)^{\frac{1}{p}}, \\
\label{eq:mixture-weight}
\kappa_x & \coloneqq \max\bigg\{1 - \tau\bigg(\sum_{i=1}^{n} \I{\|x - x_i\|_2 \leq r h_n}\bigg)^{\rme}, ~0\bigg\}.
\end{align}
Using the bandwidth $h_n$ specified in \Cref{thm:concentration-np} and a parameter $\tau >0$,
we define the interpolation weight $\kappa_x$, as given in \eqref{eq:mixture-weight}, with
$\rme = -\frac{p^2}{d_y}$ for $p < \frac{d_y}{2}$ and $\rme = -\frac{p}{2}$ for  $p \geq \frac{d_y}{2}$.
Clearly, $\kappa_x$ is influenced by covariate shifts.
It biases the mixture estimator $\mume$ towards $\munp$ when there is a sufficient number of nearby samples, and conversely towards $\mup$ when such samples are scarce.
Building on this interpolation,
we can establish guarantees for the inclusion of $\mu_{Y|x}$ in $\setame$ 
\begin{theorem}[Coverage Guarantees of $\setame$]\label{thm:coverage-iwdro-apx}
Under Assumptions \ref{asp:regular_dist} and \ref{asp:kernel},
for any $x \in \Xscr$ and a risk level $\alpha \in (0, 1)$,
we have $\P(\mu_{Y|x} \in \setame) \geq 1- \alpha$.
\end{theorem}

Given that both $\munp$ and $\mup$ are discrete distributions over $n$ points,
the mixture estimator $\mume$, as defined in \eqref{eq:mixture-estimator},
effectively becomes a discrete distribution across $2n$ points:
\begin{align}
\label{eq:mixture-estimator-weight}
\mume = \kappa_x\sum_{i = 1}^n\frac{K((x - x_i)/h_n)}{\sum_{j = 1}^n K((x-x_j)/h_n)}\delta_{y_i}
+ (1-\kappa_x)\sum_{i = 1}^n \frac{1}{n} \delta_{\{f_{\hat{\theta}}(x) + y_i - f_{\hat{\theta}}(x_i)\}}.
\end{align}
Hence, the approximation problem \eqref{eq:iwdro-opt-approximate} below is a single Wasserstein-ball DRO program
\begin{align}
\label{eq:iwdro-opt-approximate}
\min_{z \in \Zscr} \sup_{\mu \in \setame} \E_{Y\sim\mu}[c(z,Y)] \tag{{\iwdroapx}},
\end{align}
which involves $\Theta(n)$ decision variables and constraints.
Letting $\zme$ denote the optimal solution to \eqref{eq:iwdro-opt-approximate},
we can establish the following performance guarantees.
\begin{proposition}[Performance of {\iwdroapx}]\label{prop:optimality-gap-iwdro-approximate}
For any $ x\in \Xscr$ and a risk level $\alpha \in (0, 1)$,
if $\mu_{Y|x} \in \setame$, then we have:
\begin{align}
\label{eq:optimality-iwdro-approximate-bound}
& \E_{\mu_{Y|x}}[c(z^{*}, Y)] \le \sup_{\mu \in \setaiw} \E_{Y\sim\mu}[c(\ziw, Y)] \le \sup_{\mu \in \setame} \E_{Y\sim\mu}[c(\zme, Y)], \\
\label{eq:optimality-iwdro-approximate-gap}
& \E_{\mu_{Y|x}}[c(\zme, Y)] - \E_{\mu_{Y|x}}[c(z^{*}, Y)] \le 2 \mathcal{L}_{z^*} \varepsilonme \leq \Cme \mathcal{L}_{z^*} \min\{\varepsilonnp, \varepsilonp\},
\end{align}
where $\Cme > 2$ is a constant independent of $n$.
\end{proposition}

\Cref{prop:optimality-gap-iwdro-approximate} shows that {\iwdroapx} also provides an upper bound \eqref{eq:optimality-iwdro-approximate-bound} for the contextual optimization problem.
In terms of the performance gap \eqref{eq:optimality-iwdro-approximate-gap}, {\iwdroapx} is less effective than {\iwdro} by a certain constant margin.

\subsection{Connection to the Ensemble Estimator}
Our approximation model {\iwdroapx} is related to ensemble learning,
a widely used strategy in statistical and machine learning to enhance prediction quality.
Ensemble learning typically involves forming a convex combination of multiple estimators.
For example, in our context, a conventional ensemble estimator (e.g., \citealp{breiman2001random}) might assign a fixed weight of $0.5$ to each of $\munp$ and $\mup$.
While both the ensemble method and our intersection-based {\iwdro} utilize two estimators, their effectiveness in addressing covariate shifts may differ.
To illustrate this distinction, we revisit the prediction task from \Cref{ex:motivate} and compare their confidence intervals.
\begin{figure}[ht]
\SingleSpacedXI
\caption{Performances of ensemble estimators under the same setup as \Cref{fig:motivation}.}
\begin{center}
\label{fig:motivation2}
\subfigure[Ensemble estimator with fixed equal weights]{\includegraphics[width = 0.49\textwidth, trim = 0 10 0 10]{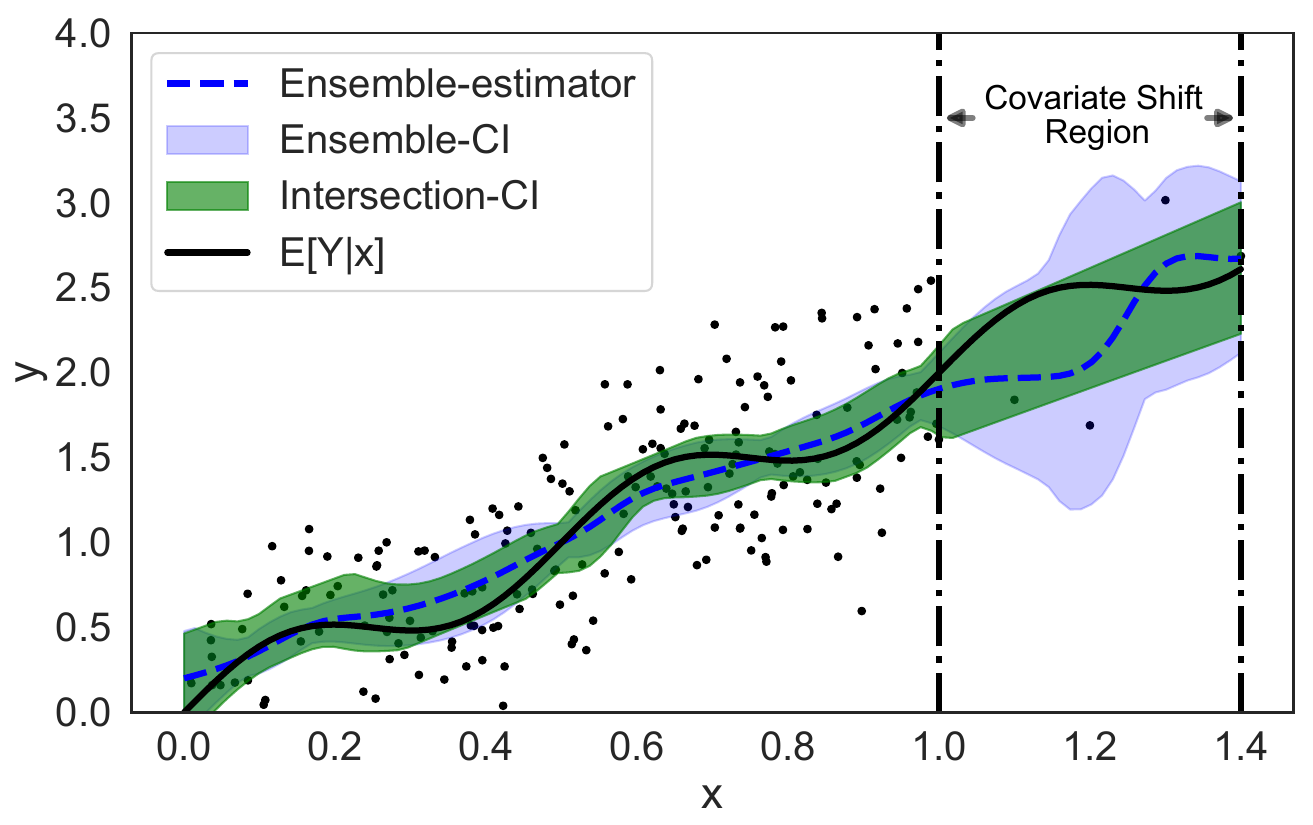}}
\subfigure[Ensemble estimator with adaptive weights]{\includegraphics[width = 0.49\textwidth, trim = 0 10 0 10]{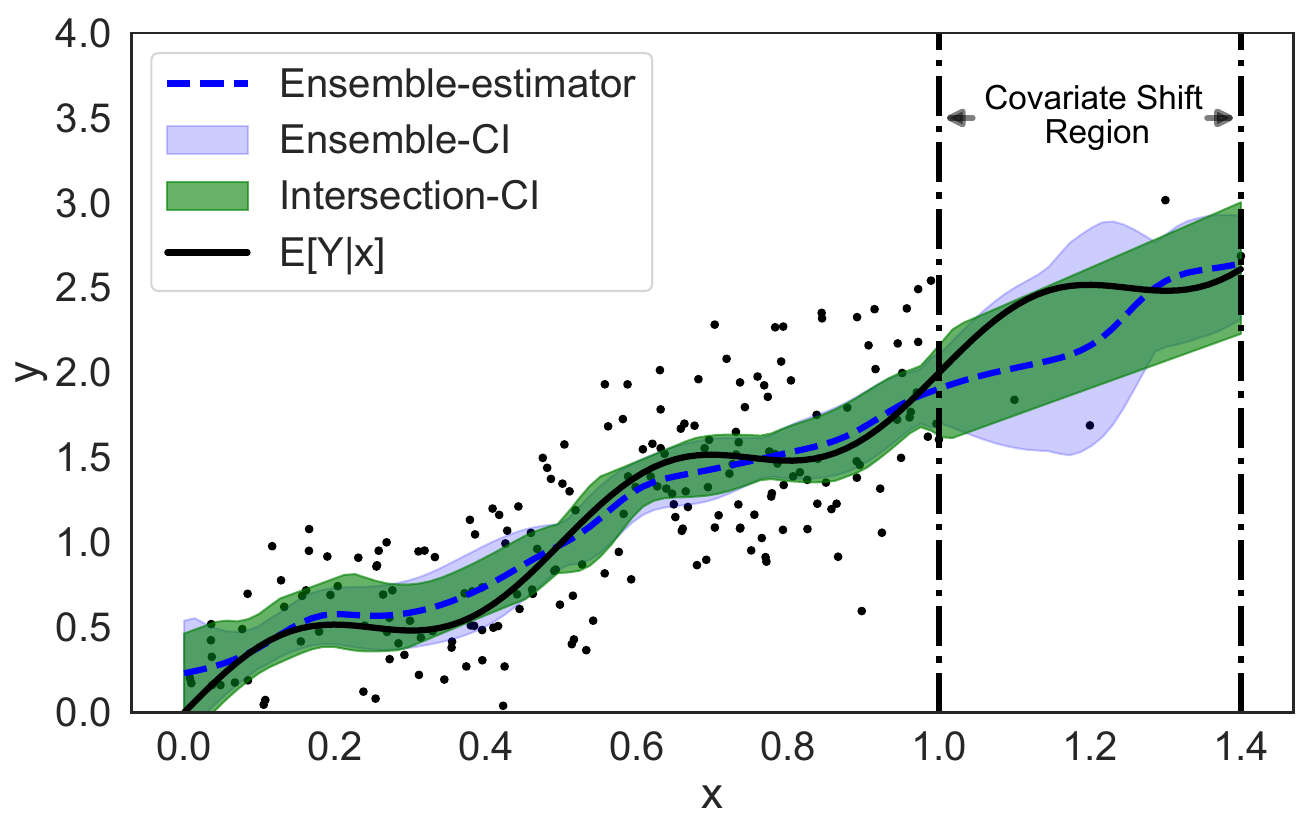}}
\end{center}
\end{figure}

As shown in \Cref{fig:motivation2}(a), in the covariate-shift region of $x>1$,
the fixed-weight ensemble estimator produces predictions (the blue dash line) that still noticeably deviate from the true conditional means (the black solid line).
This is accompanied by a wider confidence interval (the violet area) compared to the intersection-based interval (the green area).
This inefficacy can be attributed to the simple averaging of $\munp$ and $\mup$ in the fixed-weight ensemble estimator,
which consistently retains the limitations of the nonparametric estimator under covariate shifts.

The performance of the ensemble estimator can be improved through adaptive construction.
For example, $\mume$, the ball center used in {\iwdroapx}, exemplifies such an approach.
It assigns a larger weight to the parametric estimator (the linear regression model)
in cases of severe covariate shifts ($x > 1$) and vice versa.
As shown in \Cref{fig:motivation2}(b),
$\mume$ produces more accurate predictions than the fixed-weight ensemble estimator.
Nevertheless, it still leads to a wider confidence interval compared to the intersection-based approach.
Given that our objective is to make high-quality decisions in the downstream DRO model,
a narrower confidence interval is indicative of a more compact ambiguity set, leading to less conservative outcomes.
Hence, these empirical observations suggest that
the intersection-based approach may be better suited for contextual optimization problems than ensemble learning.
This concurs with Propositions \ref{prop:optimality-gap-iwdro} and \ref{prop:optimality-gap-iwdro-approximate}, which demonstrate that {\iwdro} provides a stronger performance guarantee compared to {\iwdroapx}.

\section{Numerical Studies}\label{sec:numerical}
We evaluate our proposed main model {\iwdro} and the approximation model {\iwdroapx} across two applications: an income prediction problem and a portfolio allocation problem.

\subsection{Income Prediction}\label{subsec:income-predict}
We first test {\iwdro} in an income prediction task, following the illustrative example discussed in the introduction section.
This prediction problem fits straightforwardly within the contextual optimization framework \eqref{eq:generic-cso},
where the realized uncertainty $y$ denotes the actual annual income per person, and
the decision $z$ represents the predicted income based on the individual's features (covariates) $x$.
To measure prediction errors,
we use a cost function $c(z,y) = |z - y|$,
and thus the optimization problem \eqref{eq:iwdro-opt}
can be reformulated as a linear program when we choose the 1-Wasserstein distance ($p = 1$) for computation, according to \Cref{coro:linear}.

\paragraph{Datasets.}
We use the ACS-Income dataset for the state of California \citep{ding2021retiring} for our experiments.
Specifically, we categorize individuals in the dataset as either young or old using an age threshold of 25 years.
For each experiment instance, we construct a training dataset by randomly selecting $n=50$ samples,
ensuring that it comprises $(1-m) \times 100\%$ of old individuals and $m \times 100\%$ of young individuals.
Similarly, we create a test dataset by randomly drawing another distinct $50$ samples, with a reversed proportion:
$m \times 100\%$ old individuals and $(1-m) \times 100\%$ young individuals.
We examine three scenarios of covariate shifts by choosing $m \in \{0.8, 0.9, 0.95\}$,
and generate $100$ instances for each scenario.
It is clear that a larger $m$ value leads to a greater divergence in the covariate distribution between the training and test datasets,
indicating more severe covariate shifts.

\paragraph{Benchmark Policies.}
To assess the performance of {\iwdro},
we compare it with conventional DRO models that use a single Wasserstein-ball ambiguity set $\Ascr = \{\mu: \Wscr_1(\mu, \hat\mu)\leq \varepsilon\}$.
In selecting the ball center $\hat\mu$,
we consider an \textsf{NP-DRO} model which uses the nonparametric estimator $\munp$ defined in \eqref{eq:ball-center-np}
and a \textsf{P-DRO} model which uses the parametric estimator $\mup$ defined in \eqref{eq:ball-center-p}.
For each model, we conduct a 4-fold cross-validation to choose the radius $\varepsilon$ that minimizes the average objective value within the training dataset.
The calibration details are provided in Appendix \ref{app:numerical-prediction}.
We then evaluate each model's decisions with the test dataset and compute the average objective, denoted by \emph{OBJ}, over the 50 test samples.

\paragraph{Policy Comparison.}
In \Cref{table:numerical-prediction-result},
we report the average OBJ for each of the three policies across the 100 instances within each covariate-shift scenario,
with the relative performance of each benchmark policy compared to {\iwdro} shown in parentheses.
We note that between \textsf{NP-DRO} and \textsf{P-DRO}, the former performs better under small covariate shifts, while the latter outperforms when the shift becomes large. In contrast, {\iwdro} consistently achieves a lower average OBJ than both \textsf{NP-DRO} and \textsf{P-DRO} across all scenarios, demonstrating its robustness to covariate shifts.
These findings confirm the effectiveness of our intersection-based approach,
aligning with the discussions on nonparametric kernel estimators and parametric regression estimators in Example 1.

\begin{table}[htb]
\SingleSpacedXI
\small
\renewcommand{\arraystretch}{1.1}
\renewcommand{\tabcolsep}{3mm}
\caption{Policies performances comparisons for the income prediction task}
\label{table:numerical-prediction-result}
\begin{center}
\begin{tabular}{l ccc cc}
\toprule
\multicolumn{1}{c}{Covariate Shift} &\textsf{NP-DRO} &\textsf{P-DRO}  &{\iwdro} \\
\midrule
large  shift ($m=0.95$) & 59.19 \scriptsize(119.3\%) & 58.58 \scriptsize(118.1\%)   & \textbf{49.61} \\
medium shift ($m=0.9$)  & 56.05 \scriptsize(123.5\%) & 56.91 \scriptsize(125.4\%)   & \textbf{45.39} \\
small  shift ($m=0.8$)  & 49.69 \scriptsize (116.6\%) & 55.26 \scriptsize (129.7\%) & \textbf{42.60} \\
\bottomrule
\end{tabular}
\end{center}
\smallskip
\footnotesize
{\em Note}:
Each row reports the average OBJ for three policies across the 100 instances within the same scenario.
Boldfaced value highlights the best policy with the smallest average OBJ.
The relative performance of each policy compared to {\iwdro} are shown in parentheses.
\end{table}

\subsection{Portfolio Allocation}\label{sec:numerical-portfolio}
We consider a portfolio allocation problem
where $Y \in \mathbb{R}^{d_y}$ represents the uncertain returns of $d_y$ assets
and $z \in \mathcal{Z}= \{z \in \mathbb{R}^{d_z} | \sum_{i=1}^{d_z} z_i = 1, z_i \ge 0\}$ (with $d_y=d_z$) denotes the allocation weights, subject to the no-short-sale constraint.
Within our {\iwdro} framework,
we aim to construct an intersection ambiguity set $\setaiw$
and minimize the worst-case expected portfolio loss while incorporating the conditional value-at-risk (CVaR) at a confidence level of $\phi =0.05$,
\begin{align}\label{eq:mean-cvar-problem1}
\min_{z \in \Zscr} \sup_{\mu \in \setaiw} ~ \Big\{\E_{Y\sim\mu}[-z^{\top}Y] + \text{CVaR}_{\phi}(-z^{\top}Y)\Big\}.
\end{align}
Using the linear definition of CVaR from \cite{rockafellar2002conditional},
the mean-CVaR objective can be reformulated as a piecewise linear function by introducing an auxiliary variable $r$, which represents the $\phi$-quantile value-at-risk.
The {\iwdro} formulation \eqref{eq:mean-cvar-problem1} is then equivalent to:
\begin{align}\label{eq:mean-cvar-problem2}
\min_{z\in \Zscr, r} \sup_{\mu \in \setaiw} ~\E_{Y\sim\mu} \bigg[\max\Big\{-\big(1 + 1/\phi\big) z^{\top}Y + \big(1 - 1/\phi\big) r, ~ -z^{\top}Y + r\Big\}\bigg].
\end{align}
This formulation can be simplified to a linear program when the 1-Wasserstein distance ($p = 1$) is used for computation.
The reformulation details are provided in Appendix \ref{app:numerical-portfolio}.

\subsubsection{A Synthetic-Data Study}\label{sec:numerical-portfolio-synthetic}
In the first experiment, we use synthetic data to simulate covariate shift and model misspecification error at various magnitudes.

\paragraph{Data Generation.} We consider a setting involving five assets ($d_y = 5$)
where the returns $Y$ are influenced by three factors, represented by a covariate vector $X \in \R^{d_x}$ ($d_x = 3$).
The relationship between the uncertain returns (in percentage terms) and covariates is defined as follows:
\begin{align}\label{eq:mean-cvar-return-dynamic}
y = m B_0 + B_1 x + m B_2 x\odot x + \eta,
\end{align}
where $\odot$ denotes elemental-wise multiplication,
$\eta \in \R^{d_y}$ is an idiosyncratic noise vector following a normal distribution $\mbox{N}(\mathbf{0}, 0.5 \mathbf{I})$,
$B_0 = (0.1, 0.1, -0.1)^{\top}$, and $B_1, B_2\in\R^{d_y\times d_x}$ are matrices with values
\begin{align}
B_1 = \left[\begin{array}{rrrrr}
0.1  ~& 0.1 ~& 0.1 ~& -0.1 ~& -0.1 \\
0 ~& 0 ~& 0  ~& 0 ~& 0 \\
0  ~& 0 ~& 0 ~& 0 ~& 0 \\
\end{array}\right]^\top,
\qquad
B_2 = \left[\begin{array}{rrrrr}
-0.1  ~& -0.1 ~&  0.1 ~& 0.1 ~& -0.1 \\
-0.1  ~&  0.1 ~& -0.1 ~& 0.1 ~& -0.1 \\
-0.1  ~& -0.1 ~&  0.1 ~& 0.1 ~&  0.1 \\
\end{array}\right]^\top. \nonumber
\end{align}
The data generating process \eqref{eq:mean-cvar-return-dynamic} includes a nonlinear term $mB_2 x\odot x$,
which results in model misspecification when an OLS estimator is applied (see \Cref{prop:p-estimator-example}).
We select $m \in \{0.1, 0.4\}$ for two cases that correspond to small and large model misspecification errors, respectively.

Recall that $\mu_{X}$ and $\nu_{X}$ denote the marginal distributions of covariate $X$ in the training and test datasets, respectively.
We examine three degrees of covariate shift:
(1) no shift: both $\mu_{X}$ and $\nu_X$ follow the same normal distribution $\mbox{N}(\bm u_1, \Sigma)$;
(2) severe shift: $\mu_{X} \sim \mbox{N}(\bm u_1, \Sigma)$ and $\nu_{X} \sim \mbox{N}(\bm u_2, \Sigma)$;
(3) mild shift: $\mu_{X} \sim \mbox{N}(\bm u_1, \Sigma)$ and $\nu_{X} \sim \mbox{N}(\frac{\bm u_1 + \bm u_2}{2}, \Sigma)$.
Here, $\bm u_1 = (-0.5, -0.5, -0.5)^{\top}$, $\bm u_2 = (1,1,1)^{\top}$, and $\Sigma = 2 \mathbf{I}$.
This setup, together with \eqref{eq:mean-cvar-return-dynamic}, yields the annual return $Y$
with a mean in the range of -40\% - 40\% and a volatility in the range of 10\% - 20\% for each asset in the training dataset.

Considering the three degrees of covariate shift and the two cases of model misspecification error,
we have a total of six scenarios.
For each scenario, we generate 100 instances.
In each instance, we simulate 100 covariate data points from distribution $\mu_{X}$ and produce the corresponding asset returns using \eqref{eq:mean-cvar-return-dynamic}.
These samples are then split to form a training dataset of 50 samples and a validation dataset with the remaining 50 samples.
We also simulate 20 samples based on covariate drawn from $\nu_{X}$ to create a test dataset.

\paragraph{Benchmark Policies.}
Besides \textsf{NP-DRO} and \textsf{P-DRO},
we test the approximation model {\iwdroapx}. 
Given an instance,
we train the estimator for a policy (e.g., $\munp$ or $\mup$) using the training dataset,
and then calibrate the ambiguity set radius through a 4-fold cross-validation.
Specifically, for {\iwdro}, we set $\varepsilonnp = k_1(1 + k_2) \Wscr_1(\munp, \mup)$ and $\varepsilonp = (1 - k_1)(1 + k_2) \Wscr_1(\munp, \mup)$,
where the hyperparameters are chosen from
$k_1 \in \{0.025, 0.05,\ldots, 0.95, 0.975\}$ and $k_2 \in \{0.01, 0.02, 0.05, 0.1\}$.
(This specification ensures that $\varepsilonnp + \varepsilonp > \Wscr_1(\munp, \mup)$ and thus results in a nonempty $\setaiw$ according to \Cref{prop:intersect-ambiguity}.)
For \textsf{NP-DRO}, 
we set $\varepsilon = k/\sum_{i = 1}^n K((x - x_i)/h_n)$ with a hyperparameter $k \in \{1,2,5,10\}$ and a kernel bandwidth $h_n= 10 n^{-\frac{1}{d_x + 2}}$.
For \textsf{P-DRO}, we set $\varepsilon \in \{0.05, 0.1,0.2,0.5\}$.
For each policy, we select the best radius from the search grids that minimizes the empirical mean-CVaR objective across the validation dataset.
After calibration, we evaluate the policy's decision with the test dataset,
where we compute the average mean-CVaR objective, denoted by OBJ, over the 20 test samples.

\paragraph{Policy Comparison.}
In Table \ref{table:numerical-portfolio-synthetic-result},
we report the average OBJ for each of the four policies across the 100 instances within each scenario,
with the relative performance of each policy compared to {\iwdro} shown in parentheses.

\begin{table}[htb]
\SingleSpacedXI
\small
\renewcommand{\arraystretch}{1.1}
\renewcommand{\tabcolsep}{3mm}
\caption{Policies performances comparisons for the portfolio allocation study with synthetic-data}
\label{table:numerical-portfolio-synthetic-result}
\begin{center}
\begin{tabular}{cc cccc}
\toprule
Covariate Shift & Model Error &\textsf{NP-DRO} &\textsf{P-DRO} &{\iwdroapx} &{\iwdro} \\
\midrule
severe shift & small & 0.5320 \scriptsize(122.4\%) & 0.4438 \scriptsize(102.1\%) & 0.4357 \scriptsize(100.3\%) & \textbf{0.4346} \\
& large & 0.5528 \scriptsize(115.6\%) & 0.5739 \scriptsize(120.1\%) & 0.4793 \scriptsize(100.2\%) & \textbf{0.4782} \\
\midrule
mild shift & small & 0.5161 \scriptsize(116.7\%) & 0.4535 \scriptsize(102.8\%) & 0.4496 \scriptsize(101.7\%) & \textbf{0.4423} \\
& large & 0.5318 \scriptsize(117.5\%) & 0.5053 \scriptsize(111.6\%) & 0.4562 \scriptsize(100.8\%) & \textbf{0.4527} \\
\midrule
no shift & small & 0.4895 \scriptsize (118.6\%) & 0.4209 \scriptsize (102.0\%) & \textbf{0.4087} \scriptsize (99.1\%) & 0.4126 \\
& large & 0.4174 \scriptsize (111.6\%) & 0.4254 \scriptsize(113.7\%) & 0.3767 \scriptsize(100.7\%) & \textbf{0.3741} \\
\bottomrule
\end{tabular}
\end{center}
\smallskip
\footnotesize
{\em Note}:
Each row reports the average OBJ for four policies across the 100 instances within the same scenario.
Boldfaced value highlights the best policy with the smallest average OBJ.
The relative performance of each policy compared to {\iwdro} are shown in parentheses.
\end{table}

\vspace{-0.2cm}
We note that \textsf{NP-DRO} tends to underperform as the degree of covariate shift increases.
While \textsf{P-DRO} demonstrates certain resilience to covariate shift when the model misspecification error is small,
its performance significantly deteriorates under large such errors.
These observations illustrate the limitations of DRO models that rely on a single Wasserstein ball centered on $\munp$ or $\mup$,
because the nonparametric and parametric estimators are highly sensitive to the density of local samples and the accuracy of the regression model, respectively.
In contrast, we observe that {\iwdro} consistently outperforms both \textsf{NP-DRO} and \textsf{P-DRO} by achieving a lower average OBJ across all scenarios.
This finding substantiates our rationale for {\iwdro}'s ability to tackle covariate shifts by employing an intersection ambiguity set, which combines the strengths of various estimators while also mitigating excessive conservatism.
Additionally, we find that {\iwdroapx} serves as a robust proxy for {\iwdro}, demonstrating very close performances.

\subsubsection{A Real-Data Study}\label{sec:numerical-portfolio-real}
In the second experiment, we utilize real-world data from Kenneth French's data library
website (\url{http://mba.tuck.dartmouth.edu/pages/faculty/ken.french/data_library.html}).
We consider four datasets:
(1) 6 Portfolios Formed on Size and Book-to-Market (6-FF);
(2) 10 Industry Portfolios (10-Ind);
(3) 25 Portfolios Formed on Book-to-Market and Investment (25-FF);
and (4) 30 Industry Portfolios (30-Ind).
For each dataset,
we treat the asset returns as $Y$ ($d_y = 6, 10, 25, 30$ respectively) and take the concurrently observed factors from the Fama-French three-factor model as the covariate $X$ ($d_x=3$).
All data are on a monthly basis.
We employ the rolling-sample approach following \cite{demiguel2007optimal} on data spanning $T = 666$ months from July 1963 to December 2018, with an estimation window of $n=60$. For each test month $t \in \{n+1,...,T\}$, we use the data from the preceding $n$ months, i.e., months $t-n$ to $t-1$, as training samples to solve the portfolio allocation problem of $t$.

\begin{figure}[ht]
\SingleSpacedXI
\begin{center}
\caption{The evolution of effective number of samples}
\label{fig:portfolio-real-covariate-shift}
\includegraphics[width=0.75\textwidth, trim = 0 0 0 0]{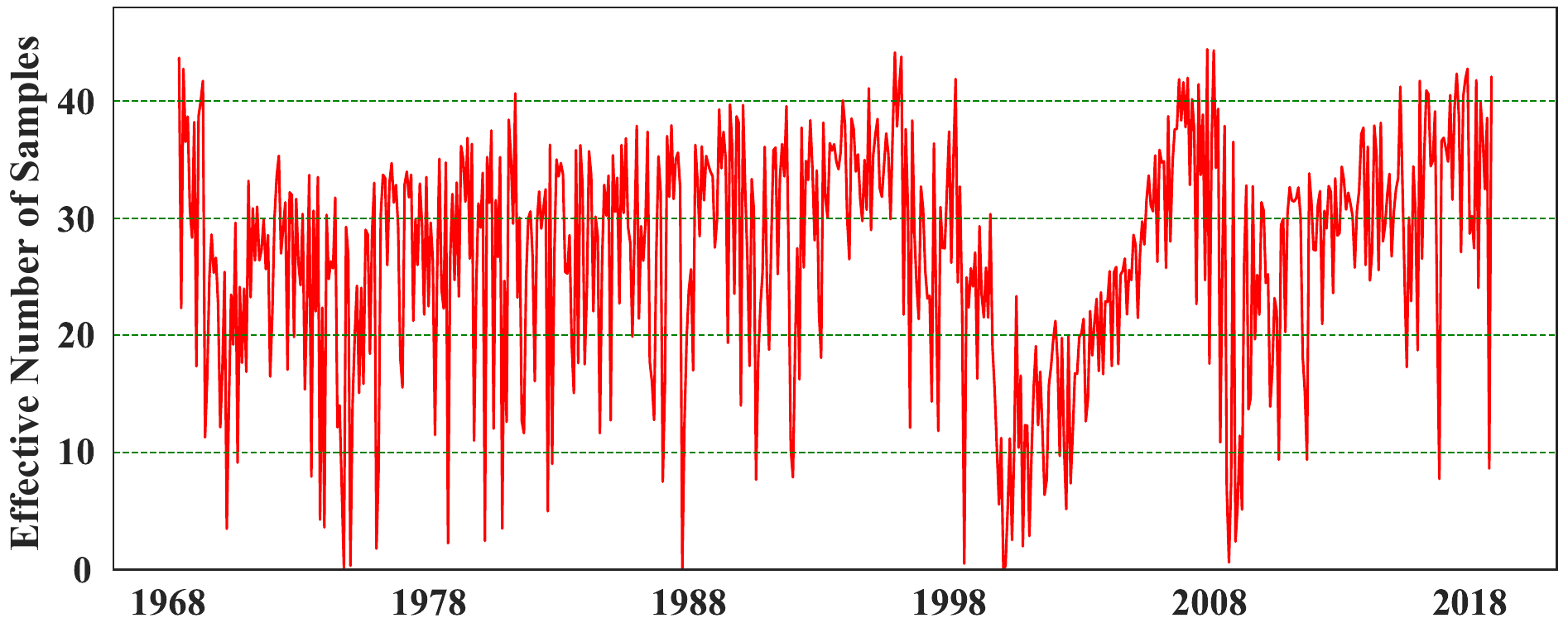}
\end{center}
{\footnotesize {\em Note}: We display the evolution of the quantity $N_t = \sum_{i=t-n}^{t-1}K((x_t - x_i)/h_n)$,
where we choose $K(\tau) = \exp(-\|\tau\|^2)$ and $h_n = 0.1 n^{-\frac{1}{d_x + d_y}}$ with $n=60$, $d_x=3$, and $d_y=10$ for 10 assets in the 10-Ind dataset.}
\vspace{0mm}
\end{figure}

To illustrate the occurrence of covariate shift within real-world datasets,
we examine the effective number of samples defined as $N_t \coloneqq \sum_{i=t-n}^{t-1}K((x_t - x_i)/h_n)$,
i.e., the denominator of the NW kernel estimator.
This measure reflects the density of training samples near the test covariate $x_t$,
thus with low values suggesting potential covariate shifts.
(Recall that the kernel function $K$ assigns a larger value to sample $x_i$ that are closer to $x_t$.)
As shown in \Cref{fig:portfolio-real-covariate-shift}, the effective number of samples changes over time, which indicates that the extent of covariate shift varies as well.
Notably, there are very few effective samples in the years 1998 and 2008,
which align with periods of global financial crises, suggesting severe covariate shifts during those times.

\paragraph{Other Benchmarks.}
Besides the \textsf{NP-DRO} and \textsf{P-DRO} evaluated in the synthetic-data study in \Cref{sec:numerical-portfolio-synthetic},
we benchmark several robust contextual optimization approaches from the literature on CVaR portfolio optimization.
We include:
(1) the residual DRO model (\textsf{R-DRO}) from \cite{kannan2020residuals} which select a Wasserstein-ball center $\hat\mu = \frac{1}{n}\sum_{i = 1}^n \delta_{\{\hat f(x) + y_i - \hat f(x_i)\}}$ with an NW kennel regression estimator $\hat{f}(x) = \sum_{i = 1}^n \frac{K((x-x_i)/h_n)}{\sum_{j = 1}^n K((x - x_j)/h_n)} y_i$;
(2) the joint-covariate-uncertainty DRO model (\textsf{OTCMV}) from Proposition 3.4 in \cite{nguyen2021robustifying};
and (3) the robust contextual optimization model (\textsf{RNW}) from Proposition 4 in \cite{srivastava2021data}.
We also include the equally-weighted model (\textsf{EW}) 
due to its robust empirical performance as documented in \cite{demiguel2007optimal}.
For each DRO model, we select the radius of the ambiguity set through a 4-fold cross-validation
by choosing the best configuration that maximizes the empirical Sharpe Ratio over the training samples.
The calibration details are provided in Appendix \ref{app:numerical-portfolio}.

We evaluate portfolio policies through three metrics:
the mean-CVaR objective (OBJ), Sharpe ratio (SR), and the certainty-equivalent return (CER).
Specifically, we run each policy using the rolling-sample approach to obtain a return series $\Rscr: = \{r_t\}_{t=n+1}^{T}$,
and then we compute the above three metrics with the formula from \cite{demiguel2007optimal} as follows:
\begin{align}\label{eq:portfolio-real-metric}
\text{OBJ} = - \widehat{\text{mean}}(\Rscr) + \widehat{\text{CVaR}}(\Rscr),
\quad \text{SR} = \widehat{\text{mean}}(\Rscr)/\widehat{\text{std}}(\Rscr),
\quad \text{CER} = \widehat{\text{mean}}(\Rscr) - (\widehat{\text{std}}(\Rscr))^2,
\end{align}
where $\widehat{\text{mean}}(\Rscr)$, $\widehat{\text{std}}(\Rscr)$, and $\widehat{\text{CVaR}}(\Rscr)$ denote the empirical computations of the mean, standard deviation, and CVaR on the return samples in $\Rscr$.

\paragraph{Policy Comparison.}
In Table \ref{table:numerical-portfolio-real-result},
we report the overall performances of various policies. 
We find that {\iwdro} achieves the best performance across all three metrics
(i.e., the lowest OBJ, and the highest SR and CER) in datasets 6-FF, 10-Ind, and 30-Ind,
and it yields the highest CER and the second highest SR in dataset 25-FF.
We also note that {\iwdroapx} perform closely to {\iwdro} and generally outperforms other benchmark policies across the four datasets.
These observations suggest that {\iwdro} and {\iwdroapx} can achieve better empirical performance compared with others benchmarks in the long run.

\begin{table}[htb]
\SingleSpacedXI
\small
\renewcommand{\arraystretch}{1.1}
\renewcommand{\tabcolsep}{2.mm}
\caption{Policies performances comparisons for the portfolio allocation study with real-data}
\label{table:numerical-portfolio-real-result}
\begin{center}
\begin{tabular}{ll cccc cccc}
\toprule
Dataset  & Metric &\textsf{EW} &\textsf{R-DRO} &\textsf{OTCMV} &\textsf{RNW} &\textsf{NP-DRO} &\textsf{P-DRO} &\textsf{IW-DRO-Apx} &\iwdro\ \\
\midrule
6-FF & OBJ &0.0987 &0.0912 &0.0853 &0.0868 &0.0892 &0.0903 & 0.0878 & \bf{0.0822} \\
& SR &0.2045 &0.2279 &0.2321 &0.2331 &0.2269 &0.2285 &0.2284 &\textbf{0.2399} \\
& CER &0.0076 &0.0083 &0.0082 &0.0086 &0.0083 &0.0085 &0.0083 &\textbf{0.0089} \\
\midrule
10-Ind & OBJ &0.0824 &0.0728 &0.0765 &0.0735 &0.0720 &0.0759 & {0.0718} &\bf{0.0702} \\
& SR &0.2203 &0.2597 &0.2365 &0.2290 &0.2545 &0.2527 &0.2577 &\textbf{0.2647} \\
& CER &0.0075 &0.0084 &0.0076 &0.0077 &0.0082 &0.0085 &0.0084 &\textbf{0.0090} \\
\midrule
25-FF & OBJ &0.1017 &0.0900 &0.0859 &\textbf{0.0857} &0.0906 &0.0909 &0.0883 &0.0915 \\
& SR &0.2029 &0.2317 & \bf{0.2447} &0.2362 &0.2318 &0.2407 &0.2418 & 0.2430 \\
& CER &0.0078 &0.0085 &0.0087 &0.0083 &0.0084 &0.0090 &0.0089 &\textbf{0.0092} \\
\midrule
30-Ind & OBJ &0.0948 &0.0748 &0.0826 & 0.0767 &0.0763 &0.0761 &\bf{0.0740} &0.0766 \\
& SR &0.1981 &0.2297 &0.2428 &0.2310 &0.2476 &0.2498 & \bf{0.2681} & {0.2673} \\
& CER &0.0072 &0.0074 &0.0084 &0.0074 &0.0082 &0.0090 &0.0090 &\textbf{0.0098} \\
\bottomrule
\end{tabular}
\end{center}
\smallskip
\footnotesize
{\em Note:} The metrics of the mean-CVaR objective (OBJ), Sharpe ratio (SR), and the  certainty-equivalent return (CER) are defined in \eqref{eq:portfolio-real-metric}.
Boldfaced value highlights the policy with the best metric, i.e., the lowest OBJ, the highest SR or CER.
\end{table}

Next, we investigate the performance improvement for {\iwdro}.
We sort the effective number of samples $N_t = \sum_{i=t-n}^{t-1}K((x_t - x_i)/h_n)$ over test months $t \in \{n+1,...,T\}$,
and divide these months into two subsets $\Tscr_1$ and $\Tscr_2$,
which contain months associated with the upper and lower halves of $N_t$ values, respectively.
Hence, $\Tscr_1$ corresponds to periods with potentially large covariate shifts, while $\Tscr_2$ corresponds to periods with potentially small covariate shifts.
In \Cref{table:numerical-portfolio-real-htest2},
we report results of one-tailed tests that assess the Sharpe ratio improvement of {\iwdro} compared to other policies in $\Tscr_1$ and $\Tscr_2$.
We find that {\iwdro} consistently achieves a higher SR than \textsf{NP-DRO} and \textsf{P-DRO} across all tests.
More notably, the results clearly show that the magnitude of SR improvement is greater and the tests of statistical significance are more frequent in $\Tscr_1$ than in $\Tscr_2$.
These findings indicate the efficacy of {\iwdro} in addressing covariate shifts, especially when these shifts are large.

\begin{table}[htb]
\SingleSpacedXI
\small
\renewcommand{\arraystretch}{1.1}
\renewcommand{\tabcolsep}{1.8mm}
\caption{\small Sharpe ratio improvement of {\iwdro} in periods with large and small covariate shifts}
\label{table:numerical-portfolio-real-htest2}
\begin{center}
\begin{tabular}{c | llll | llll}
\toprule
SR improvement & \multicolumn{4}{c|}{$\Tscr_1$: large covariate shifts} & \multicolumn{4}{c}{$\Tscr_2$: small covariate shifts} \\
{\iwdro} compared to& \multicolumn{1}{c}{6-FF} & \multicolumn{1}{c}{10-Ind} & \multicolumn{1}{c}{25-FF} & \multicolumn{1}{c|}{30-Ind}
& \multicolumn{1}{c}{6-FF} & \multicolumn{1}{c}{10-Ind} & \multicolumn{1}{c}{25-FF} & \multicolumn{1}{c}{30-Ind}\\
\midrule
\textsf{NP-DRO} & $0.0084$      & $0.0555^{***}$ & $0.0238^{*}$ & $0.0739^{***}$
& $0.0192^{**}$ & $0.0088$     & $0.0041$ & $0.0107$ \\
\textsf{P-DRO}  & $0.0163^{**}$ & $0.0384^{***}$ & $0.0059$     & $0.0568^{***}$
& $0.0109^{*}$  & $0.0178^{*}$ & $0.0017$ & $0.0117^{**}$ \\
\bottomrule
\end{tabular}
\end{center}
\smallskip
{\emph{Notes.} Superscripts $^{***}$, $^{**}$ and $^{*}$ indicate
that {\iwdro} has a statistically significant higher Sharpe ratio at the 1\%, 5\%, and 10\% levels in one-tailed tests, respectively.
}
\end{table}

\section{Conclusion and Future Research} \label{sec:conclusion}
We study contextual optimization under covariate shift,
a common real-world phenomenon that poses significant performance challenges to existing decision models.
To address this, we develop a new Wasserstein-DRO framework with a novel intersection ambiguity set that combines the benefits of various estimation methods.
Our model improves the quality of contextual decision-making and provides insights from both computational and statistical perspectives.
There are several promising directions for future research. One is to explore other tractable ambiguity sets that better capture different types of distribution shifts. Another is to investigate the trade-offs between performance utility and computational efficiency when intersecting different types of ambiguity sets.







\bibliographystyle{informs2014}
\bibliography{refs_iwdro}

\begin{thebibliography}{58}
\providecommand{\natexlab}[1]{#1}
\providecommand{\url}[1]{\texttt{#1}}
\providecommand{\urlprefix}{URL }

\bibitem[{Awasthi et~al.(2022)Awasthi, Jung, \protect\BIBand{}
  Morgenstern}]{awasthi2022distributionally}
Awasthi P, Jung C, Morgenstern J (2022) Distributionally robust data join.
  \emph{arXiv preprint arXiv:2202.05797} .

\bibitem[{Ban \protect\BIBand{} Rudin(2019)}]{ban2019big}
Ban GY, Rudin C (2019) The big data newsvendor: Practical insights from machine
  learning. \emph{Operations Research} 67(1):90--108.

\bibitem[{Bennouna \protect\BIBand{} Van~Parys(2022)}]{bennouna2022holistic}
Bennouna A, Van~Parys B (2022) Holistic robust data-driven decisions.
  \emph{arXiv preprint arXiv:2207.09560} .

\bibitem[{Bertsimas \protect\BIBand{} Kallus(2020)}]{bertsimas2020predictive}
Bertsimas D, Kallus N (2020) From predictive to prescriptive analytics.
  \emph{Management Science} 66(3):1025--1044.

\bibitem[{Bertsimas et~al.(2023)Bertsimas, McCord, \protect\BIBand{}
  Sturt}]{bertsimas2023dynamic}
Bertsimas D, McCord C, Sturt B (2023) Dynamic optimization with side
  information. \emph{European Journal of Operational Research} 304(2):634--651.

\bibitem[{Bertsimas \protect\BIBand{} Van~Parys(2021)}]{bertsimas2021bootstrap}
Bertsimas D, Van~Parys B (2021) Bootstrap robust prescriptive analytics.
  \emph{Mathematical Programming} 1--40.

\bibitem[{Besbes et~al.(2022)Besbes, Ma, \protect\BIBand{}
  Mouchtaki}]{besbes2022beyond}
Besbes O, Ma W, Mouchtaki O (2022) Beyond iid: data-driven decision-making in
  heterogeneous environments. \emph{Advances in Neural Information Processing
  Systems} 35:23979--23991.

\bibitem[{Bickel et~al.(2009)Bickel, Br{\"u}ckner, \protect\BIBand{}
  Scheffer}]{bickel2009discriminative}
Bickel S, Br{\"u}ckner M, Scheffer T (2009) Discriminative learning under
  covariate shift. \emph{Journal of Machine Learning Research} 10(9).

\bibitem[{Blanchet \protect\BIBand{} Murthy(2019)}]{blanchet2019quantifying}
Blanchet J, Murthy K (2019) Quantifying distributional model risk via optimal
  transport. \emph{Mathematics of Operations Research} 44(2):565--600.

\bibitem[{Breiman(2001)}]{breiman2001random}
Breiman L (2001) Random forests. \emph{Machine learning} 45:5--32.

\bibitem[{Chen et~al.(2019)Chen, Sim, \protect\BIBand{}
  Xu}]{chen2019distributionally}
Chen Z, Sim M, Xu H (2019) Distributionally robust optimization with infinitely
  constrained ambiguity sets. \emph{Operations Research} 67(5):1328--1344.

\bibitem[{Cheramin et~al.(2022)Cheramin, Cheng, Jiang, \protect\BIBand{}
  Pan}]{cheramin2022computationally}
Cheramin M, Cheng J, Jiang R, Pan K (2022) Computationally efficient
  approximations for distributionally robust optimization under moment and
  wasserstein ambiguity. \emph{INFORMS Journal on Computing} 34(3):1768--1794.

\bibitem[{DeMiguel et~al.(2007)DeMiguel, Garlappi, \protect\BIBand{}
  Uppal}]{demiguel2007optimal}
DeMiguel V, Garlappi L, Uppal R (2007) Optimal versus naive diversification:
  How inefficient is the 1/n portfolio strategy? \emph{The Review of Financial
  Studies} 22(5):1915--1953.

\bibitem[{Ding et~al.(2021)Ding, Hardt, Miller, \protect\BIBand{}
  Schmidt}]{ding2021retiring}
Ding F, Hardt M, Miller J, Schmidt L (2021) Retiring adult: New datasets for
  fair machine learning. \emph{Advances in neural information processing
  systems} 34:6478--6490.

\bibitem[{Donti et~al.(2017)Donti, Amos, \protect\BIBand{}
  Kolter}]{donti2017task}
Donti P, Amos B, Kolter JZ (2017) Task-based end-to-end model learning in
  stochastic optimization. \emph{Advances in neural information processing
  systems} 30.

\bibitem[{Duchi et~al.(2023)Duchi, Hashimoto, \protect\BIBand{}
  Namkoong}]{duchi2020distributionally}
Duchi J, Hashimoto T, Namkoong H (2023) Distributionally robust losses for
  latent covariate mixtures. \emph{Operations Research} 71(2):649--664.

\bibitem[{Elmachtoub \protect\BIBand{} Grigas(2022)}]{elmachtoub2017smart}
Elmachtoub AN, Grigas P (2022) Smart “predict, then optimize”.
  \emph{Management Science} 68(1):9--26.

\bibitem[{Esteban-P{\'e}rez \protect\BIBand{}
  Morales(2023)}]{esteban2023distributionally}
Esteban-P{\'e}rez A, Morales JM (2023) Distributionally robust optimal power
  flow with contextual information. \emph{European Journal of Operational
  Research} 306(3):1047--1058.

\bibitem[{Fournier \protect\BIBand{} Guillin(2015)}]{fournier2015rate}
Fournier N, Guillin A (2015) On the rate of convergence in wasserstein distance
  of the empirical measure. \emph{Probability Theory and Related Fields}
  162(3):707--738.

\bibitem[{Gao \protect\BIBand{} Kleywegt(2023)}]{gao2016distributionally}
Gao R, Kleywegt A (2023) Distributionally robust stochastic optimization with
  wasserstein distance. \emph{Mathematics of Operations Research}
  48(2):603--655.

\bibitem[{Gupta \protect\BIBand{} Kallus(2022)}]{gupta2022data}
Gupta V, Kallus N (2022) Data pooling in stochastic optimization.
  \emph{Management Science} 68(3):1595--1615.

\bibitem[{Gy{\"o}rfi et~al.(2006)Gy{\"o}rfi, Kohler, Krzyzak, \protect\BIBand{}
  Walk}]{gyorfi2006distribution}
Gy{\"o}rfi L, Kohler M, Krzyzak A, Walk H (2006) \emph{A distribution-free
  theory of nonparametric regression} (Springer Science \& Business Media).

\bibitem[{Hanasusanto \protect\BIBand{} Kuhn(2013)}]{hanasusanto2013robust}
Hanasusanto GA, Kuhn D (2013) Robust data-driven dynamic programming.
  \emph{Advances in Neural Information Processing Systems} 26:827--835.

\bibitem[{Hannah et~al.(2010)Hannah, Powell, \protect\BIBand{}
  Blei}]{hannah2010nonparametric}
Hannah L, Powell W, Blei D (2010) Nonparametric density estimation for
  stochastic optimization with an observable state variable. \emph{Advances in
  Neural Information Processing Systems} 23:820--828.

\bibitem[{Hastie et~al.(2009)Hastie, Tibshirani, Friedman, \protect\BIBand{}
  Friedman}]{hastie2009elements}
Hastie T, Tibshirani R, Friedman JH, Friedman JH (2009) \emph{The elements of
  statistical learning: data mining, inference, and prediction}, volume~2
  (Springer).

\bibitem[{Ho-Nguyen \protect\BIBand{}
  K{\i}l{\i}n{\c{c}}-Karzan(2022)}]{ho2022risk}
Ho-Nguyen N, K{\i}l{\i}n{\c{c}}-Karzan F (2022) Risk guarantees for end-to-end
  prediction and optimization processes. \emph{Management Science}
  68(12):8680--8698.

\bibitem[{Hsu et~al.(2012)Hsu, Kakade, \protect\BIBand{} Zhang}]{hsu2012random}
Hsu D, Kakade SM, Zhang T (2012) Random design analysis of ridge regression.
  \emph{Conference on learning theory}, 9--1 (JMLR Workshop and Conference
  Proceedings).

\bibitem[{Hu et~al.(2018)Hu, Niu, Sato, \protect\BIBand{}
  Sugiyama}]{hu2018does}
Hu W, Niu G, Sato I, Sugiyama M (2018) Does distributionally robust supervised
  learning give robust classifiers? \emph{International Conference on Machine
  Learning}, 2029--2037 (PMLR).

\bibitem[{Iyengar et~al.(2023)Iyengar, Lam, \protect\BIBand{}
  Wang}]{iyengar2023hedging}
Iyengar G, Lam H, Wang T (2023) Hedging against complexity: Distributionally
  robust optimization with parametric approximation. \emph{International
  Conference on Artificial Intelligence and Statistics}, 9976--10011 (PMLR).

\bibitem[{Kannan et~al.(2023)Kannan, Bayraksan, \protect\BIBand{}
  Luedtke}]{kannan2020residuals}
Kannan R, Bayraksan G, Luedtke JR (2023) Residuals-based distributionally
  robust optimization with covariate information. \emph{Mathematical
  Programming} 1--57.

\bibitem[{Kantorovich \protect\BIBand{}
  Rubinshtein(1958)}]{kantorovich1958space}
Kantorovich LV, Rubinshtein S (1958) On a space of totally additive functions.
  \emph{Vestnik of the St. Petersburg University: Mathematics} 13(7):52--59.

\bibitem[{Kuhn et~al.(2019)Kuhn, Esfahani, Nguyen, \protect\BIBand{}
  Shafieezadeh-Abadeh}]{kuhn2019wasserstein}
Kuhn D, Esfahani PM, Nguyen VA, Shafieezadeh-Abadeh S (2019) Wasserstein
  distributionally robust optimization: Theory and applications in machine
  learning. \emph{Operations Research \& Management Science in the Age of
  Analytics}, 130--166 (INFORMS).

\bibitem[{Liu et~al.(2024)Liu, Wang, Cui, \protect\BIBand{}
  Namkoong}]{liu2024need}
Liu J, Wang T, Cui P, Namkoong H (2024) On the need for a language describing
  distribution shifts: Illustrations on tabular datasets. \emph{Advances in
  Neural Information Processing Systems} 36.

\bibitem[{Long et~al.(2023)Long, Sim, \protect\BIBand{} Zhou}]{long2023robust}
Long DZ, Sim M, Zhou M (2023) Robust satisficing. \emph{Operations Research}
  71(1):61--82.

\bibitem[{Mohajerin~Esfahani \protect\BIBand{} Kuhn(2018)}]{esfahani2018data}
Mohajerin~Esfahani P, Kuhn D (2018) Data-driven distributionally robust
  optimization using the wasserstein metric: Performance guarantees and
  tractable reformulations. \emph{Mathematical Programming} 171(1):115--166.

\bibitem[{Nadaraya(1964)}]{nadaraya1964estimating}
Nadaraya EA (1964) On estimating regression. \emph{Theory of Probability \& Its
  Applications} 9(1):141--142.

\bibitem[{Nguyen et~al.(2024)Nguyen, Zhang, Blanchet, Delage, \protect\BIBand{}
  Ye}]{nguyen2021robustifying}
Nguyen VA, Zhang F, Blanchet J, Delage E, Ye Y (2024) Robustifying conditional
  portfolio decisions via optimal transport. \emph{Operations Research} .

\bibitem[{Pun et~al.(2023)Pun, Wang, \protect\BIBand{} Yan}]{pun2023data}
Pun CS, Wang T, Yan Z (2023) Data-driven distributionally robust cvar portfolio
  optimization under a regime-switching ambiguity set. \emph{Manufacturing \&
  Service Operations Management} 25(5):1779--1795.

\bibitem[{Qi et~al.(2021)Qi, Grigas, \protect\BIBand{} Shen}]{qi2021integrated}
Qi M, Grigas P, Shen ZJM (2021) Integrated conditional estimation-optimization.
  \emph{arXiv preprint arXiv:2110.12351} .

\bibitem[{Qi et~al.(2023)Qi, Shi, Qi, Ma, Yuan, Wu, \protect\BIBand{}
  Shen}]{qi2022e2e}
Qi M, Shi Y, Qi Y, Ma C, Yuan R, Wu D, Shen ZJM (2023) A practical end-to-end
  inventory management model with deep learning. \emph{Management Science}
  69(2):759--773.

\bibitem[{Quinonero-Candela et~al.(2008)Quinonero-Candela, Sugiyama,
  Schwaighofer, \protect\BIBand{} Lawrence}]{quinonero2008dataset}
Quinonero-Candela J, Sugiyama M, Schwaighofer A, Lawrence ND (2008)
  \emph{Dataset shift in machine learning} (Mit Press).

\bibitem[{Rockafellar \protect\BIBand{}
  Uryasev(2002)}]{rockafellar2002conditional}
Rockafellar RT, Uryasev S (2002) Conditional value-at-risk for general loss
  distributions. \emph{Journal of Banking \& Finance} 26(7):1443--1471.

\bibitem[{Rychener et~al.(2024)Rychener, Esteban-P{\'e}rez, Morales,
  \protect\BIBand{} Kuhn}]{rychener2024wasserstein}
Rychener Y, Esteban-P{\'e}rez A, Morales JM, Kuhn D (2024) Wasserstein
  distributionally robust optimization with heterogeneous data sources.
  \emph{arXiv preprint arXiv:2407.13582} .

\bibitem[{Sadana et~al.(2025)Sadana, Chenreddy, Delage, Forel, Frejinger,
  \protect\BIBand{} Vidal}]{sadana2024survey}
Sadana U, Chenreddy A, Delage E, Forel A, Frejinger E, Vidal T (2025) A survey
  of contextual optimization methods for decision-making under uncertainty.
  \emph{European Journal of Operational Research} 320(2):271--289.

\bibitem[{Selvi et~al.(2024)Selvi, Kreacic, Ghassemi, Potluru, Balch,
  \protect\BIBand{} Veloso}]{selvican}
Selvi A, Kreacic E, Ghassemi M, Potluru V, Balch T, Veloso M (2024)
  Distributionally and adversarially robust logistic regression via
  intersecting wasserstein balls. \emph{arXiv preprint arXiv:2407.13625} .

\bibitem[{Shimodaira(2000)}]{shimodaira2000improving}
Shimodaira H (2000) Improving predictive inference under covariate shift by
  weighting the log-likelihood function. \emph{Journal of statistical planning
  and inference} 90(2):227--244.

\bibitem[{Sim et~al.(2024)Sim, Tang, Zhou, \protect\BIBand{}
  Zhu}]{sim2024analytics}
Sim M, Tang Q, Zhou M, Zhu T (2024) The analytics of robust satisficing:
  predict, optimize, satisfice, then fortify. \emph{Operations Research} .

\bibitem[{Srivastava et~al.(2021)Srivastava, Wang, Hanasusanto,
  \protect\BIBand{} Ho}]{srivastava2021data}
Srivastava PR, Wang Y, Hanasusanto GA, Ho CP (2021) On data-driven prescriptive
  analytics with side information: A regularized nadaraya-watson approach.
  \emph{arXiv preprint arXiv:2110.04855} .

\bibitem[{Sugiyama \protect\BIBand{} Kawanabe(2012)}]{sugiyama2012machine}
Sugiyama M, Kawanabe M (2012) \emph{Machine learning in non-stationary
  environments: Introduction to covariate shift adaptation} (MIT press).

\bibitem[{Sugiyama et~al.(2007)Sugiyama, Krauledat, \protect\BIBand{}
  M{\"u}ller}]{sugiyama2007covariate}
Sugiyama M, Krauledat M, M{\"u}ller KR (2007) Covariate shift adaptation by
  importance weighted cross validation. \emph{Journal of Machine Learning
  Research} 8(5).

\bibitem[{Tanoumand et~al.(2023)Tanoumand, Bodur, \protect\BIBand{}
  Naoum-Sawaya}]{tanoumand2023data}
Tanoumand N, Bodur M, Naoum-Sawaya J (2023) Data-driven distributionally robust
  optimization: Intersecting ambiguity sets, performance analysis and
  tractability. \emph{Working Paper} .

\bibitem[{Taskesen et~al.(2021)Taskesen, Yue, Blanchet, Kuhn, \protect\BIBand{}
  Nguyen}]{taskesen2021sequential}
Taskesen B, Yue MC, Blanchet J, Kuhn D, Nguyen VA (2021) Sequential domain
  adaptation by synthesizing distributionally robust experts.
  \emph{International Conference on Machine Learning}, 10162--10172 (PMLR).

\bibitem[{Tsybakov(2008)}]{tsybakov2008introduction}
Tsybakov AB (2008) \emph{Introduction to nonparametric estimation} (Springer
  Science \& Business Media).

\bibitem[{Wang et~al.(2021)Wang, Chen, \protect\BIBand{}
  Wang}]{wang2021distributionally}
Wang T, Chen N, Wang C (2021) Distributionally robust prescriptive analytics
  with wasserstein distance. \emph{arXiv preprint arXiv:2106.05724} .

\bibitem[{Watson(1964)}]{watson1964smooth}
Watson GS (1964) Smooth regression analysis. \emph{Sankhy{\=a}: The Indian
  Journal of Statistics, Series A} 359--372.

\bibitem[{Yang et~al.(2022)Yang, Zhang, Chen, Gao, \protect\BIBand{}
  Hu}]{yang2022decision}
Yang J, Zhang L, Chen N, Gao R, Hu M (2022) Decision-making with side
  information: A causal transport robust approach. \emph{Optimization Online} .

\bibitem[{Zhang et~al.(2024{\natexlab{a}})Zhang, Yang, \protect\BIBand{}
  Gao}]{zhang2023optimal}
Zhang L, Yang J, Gao R (2024{\natexlab{a}}) Optimal robust policy for
  feature-based newsvendor. \emph{Management Science} 70(4):2315--2329.

\bibitem[{Zhang et~al.(2024{\natexlab{b}})Zhang, Yang, \protect\BIBand{}
  Gao}]{zhang2022simple}
Zhang L, Yang J, Gao R (2024{\natexlab{b}}) A short and general duality proof
  for wasserstein distributionally robust optimization. \emph{Operations
  Research} .

\end{thebibliography}


\ECSwitch
\OnePointTwoSpacedXI


\vspace{-0.5cm}
\begin{center}
{\noindent{\large\textbf{E-Companion to Contextual Optimization under Covariate Shift:\\
A Robust Approach by Intersecting Wasserstein Balls}}}
\end{center}

\section{Proofs for Section~\ref{sec:IWDRO-formulation}}\label{app:IWDRO-formulation}

\subsection{Proof of Proposition~\ref{prop:intersect-ambiguity}.}

\begin{definition}[Geodesic Metric Space]
A metric space $(\Sscr, d)$ is a geodesic metric space if, for any two points $s_1, s_2 \in \Sscr$,
there exists a continuous curve $\gamma: [0, 1] \to \Sscr$
connecting $s_1$ to $s_2$ such that, for any $t_1, t_2 \in [0, 1]$,
the distance metric $d(\cdot,\cdot)$ satisfies $d(\gamma(t_1), \gamma(t_2)) = |t_1 - t_2| d(s_1, s_2)$.
\end{definition}

When $\Sscr$ is a subset of Euclidean space and $d$ is the $p$-Wasserstein distance, $(\Sscr, \Wscr_p)$ forms a geodesic metric space.
Consequently, we can obtain \Cref{prop:intersect-ambiguity} by applying Proposition 5.3 in \cite{taskesen2021sequential}.  $\hfill \square$

\subsection{Proof of Theorem~\ref{thm:strong-duality} and Results for the Extended {\iwdro}}\label{app:strongduality}
In this proof, we show that \Cref{thm:strong-duality} holds for
an intersection ambiguity set involving more than two Wasserstein balls, as well as for a generalized Wasserstein distance beyond \Cref{def:was-distance}.
\begin{itemize}
\item For any metric $d$,
we define a distance  $\Wscr_d$ between two distributions $\mu$ and $\nu$ as:
\begin{equation}\label{eq:dwasserstein-dist}
\Wscr_d(\mu,\nu) \coloneqq \inf_{\xi \in \Xi(\mu, \nu)} \E_{(Y_1, Y_2)\sim \xi}[ d(Y_1, Y_2) ],
\end{equation}
Here, we recall that $\Xi(\mu, \nu)$ denotes the set of all joint distributions with marginals $\mu$ and $\nu$.
Compared to the $p$-Wasserstein distance $\Wscr_p$ defined in \eqref{eq:pwasserstein-dist},
when $d(Y_1, Y_2) = \|Y_1 - Y_2\|_1^p$, the following two sets (i.e.,  Wasserstein balls) are equivalent:
$$\{\mu: \Wscr_p(\mu, \hat\mu)\leq \varepsilon\} = \{\mu: \Wscr_d(\mu, \hat\mu)\leq \varepsilon^p\}.$$
\item We define a generalized intersection ambiguity set formed by $M$ Wasserstein balls:
\begin{equation}\label{eq:aug-a}
\setag = \{\mu: \Wscr_d(\mu, \hat\mu_m) \leq \varepsilon_m, ~ \forall m \in [M]\}.
\end{equation}
Clearly, our proposed $\setaiw$, as defined in \eqref{eq:ambiguity-set-generic}, is a special case of $\setag$ with $M = 2$ balls.
\end{itemize}

We extend the {\iwdro} framework by using $\setag$ and $\Wscr_d$ in place of $\setaiw$ and $\Wscr_p$, respectively.
Similarly, we adhere to Assumptions~\ref{asp:dist_structure} and \ref{asp:reg-costfunc},
which state that each reference distribution $\hat\mu_m$, i.e., the center of the ball in \eqref{eq:aug-a}, is a discrete distribution of the form $\hat\mu_m = \sum_{i=1}^n w_{m,i} \delta_{y_{m,i}}$,
where the expected cost is bounded as $|\E_{Y\sim\hat\mu_m}[c(z,Y)]| < \infty$ for all $z \in \Zscr$.
We also assume that $\setag$ is non-empty.
Under these conditions, we can establish the following result.
\begin{proposition}\label{thm:strong-duality2}
For all $z \in \Zscr$, let $\JP(z) \coloneqq \sup_{\mu \in \setag}\E_{Y \sim \mu}[c(z, Y)]$.
Then, we have
\begin{align}
\hspace{-0.2cm}\JP(z) = & \inf_{\lambda_m \geq 0}\Bigg\{\sum_{m\in[M]}\lambda_m \varepsilon_m +
\sup_{\xi \in \Xi(\hat\mu_1,...,\hat\mu_M)}
\E_{(\tilde{Y}_1,..., \tilde{Y}_M)\sim \xi}
\bigg[\sup_{y \in \Yscr}\bigg\{c(z, y) -
\sum_{m\in[M]} \lambda_m d(\tilde Y_m, y) \bigg\}\bigg] \Bigg\}. \label{eq:strong-dual-extend-1}
\end{align}
\end{proposition}
Here, $(\tilde{Y}_1,..., \tilde{Y}_M)$ denotes random vectors that follow
a joint distribution $\xi$ from $\Xi(\hat\mu_1,...,\hat\mu_M)$,
which is the set of all joint distributions with marginals $\hat\mu_m$ for $m =1,...,M$,
such that $\tilde{Y}_m \sim \hat\mu_m$.

We defer the proof of \Cref{thm:strong-duality2} later.
Using \Cref{thm:strong-duality2}, it immediately follows that
$\JP(z)$ for \Cref{thm:strong-duality}, as defined in \eqref{eq:strong-dual0-cp}, is equivalent to \eqref{eq:strong-dual0-cd}, since $\setaiw$ is a special case of $\setag$.
As shown in the main text, by formulating the dual program for the inner supremum problem using \eqref{eq:reformulation-inner} and \eqref{eq:reformulation-inner-dual}, we have established \Cref{thm:strong-duality}.
$\hfill \square$

\smallskip
For the extended {\iwdro} problem with $\setag$ formed by multiple Wasserstein balls,
by applying the same duality techniques, we can reformulate \eqref{eq:strong-dual-extend-1} to derive the corresponding optimization program, analogous to \eqref{eq:strong-dual-further}.
We denote $\Iscr_{M,n}:=[n]^M$ as the set of $M$-dimensional multi-indices,
i.e., the set of permutations with repetition over $[n]^M$.
Consequently, each element $\bm{i} \in \Iscr_{M, n}$ uniquely corresponds to a combination of atoms from the $M$ discrete reference distributions $\{\hat\mu_m\}_{m \in [M]}$.
Let $\bm{i}_m$ represent the $m$-th component of $\bm{i}$.
By introducing a set of dual variables $\{u_{m,i}\}_{m \in [M], i\in[n]}$, we can establish the strong duality as follows:
\begin{theorem}[Strong Duality of the Extended {\iwdro}]\label{thm:strong-duality3}
We have $\JP(z) = \JD(z)$ with
\begin{align}
\JD(z) ~\coloneqq~
\inf \quad & \sum_{m\in[M]} \lambda_m \varepsilon_m + \sum_{m\in[M]} \sum_{i\in[n]} w_{m,i} u_{m,i} \label{eq:strong-dual-extend-2} \\
\text{\textup{s.t.}} \quad & \sum_{m\in[M]} u_{m,\bm{i}_m} \geq \sup_{y \in \Yscr} \Big\{c(z, y) - \sum_{m\in[M]} \lambda_m d(y_{m,\bm{i}_m}, y) \Big\}, \quad \forall \bm{i} \in \Iscr_{M, n}; \nonumber \\
& \lambda_m \geq 0, \quad u_{m,\bm{i}_m} \in \R, \quad \forall m \in [M], ~\bm{i} \in \Iscr_{M,n}. \nonumber
\end{align}
\end{theorem}

\smallskip
\noindent
\textbf{Proof of \Cref{thm:strong-duality2}}

Our analysis begins with the application of the Legendre transformation.
For a function $g(\bm{x})$, where $\bm{x} = (x_1,..., x_M)$, we denote its Legendre transformation as
$g^*(\bm{\lambda}) \coloneqq \sup_{\bm{x}} \{\bm{\lambda}^{\top}\bm{x} - g(\bm{x})\} $ with $\bm{\lambda}=(\lambda_1,..., \lambda_M)$.
For $\JP(z)$,
we rewrite the function in the form $J_{z}(\bm{\varepsilon}) := \sup_{\mu \in \setag}\E_{Y \sim \mu}[c(z, Y)]$,
which emphasizes the dependence of the supremum program
on the radii $\bm{\varepsilon} = (\varepsilon_1,..., \varepsilon_M)$ of the ambiguity set $\setag$.
We then apply the Legendre transformation to $-J_{z}(\bm{\varepsilon})$, yielding
\begin{align}
(-J_{z})^*(-\bm{\lambda})
& = \sup_{\bm{\varepsilon} \geq \bm{0}} \bigg\{\sup_{\mu \in \setag}\E_{Y \sim \mu}[c(z, Y)] - \bm{\lambda}^{\top}\bm{\varepsilon} \bigg\} \nonumber \\
& = \sup_{\bm{\varepsilon} \geq \bm{0}} \sup_{\mu \in \Pscr(\Yscr)} \bigg\{\E_{Y\sim\mu}[c(z,Y)] - \sum_{m\in[M]}\lambda_m \varepsilon_m,  \quad \text{s.t. } \Wscr_d(\mu, \hat\mu_m) \leq \varepsilon_m, ~ \forall m \in [M] \bigg\} \nonumber \\
& = \sup_{\mu \in \Pscr(\Yscr)} \sup_{\bm{\varepsilon} \geq \bm{0}} \bigg\{\E_{Y\sim\mu}[c(z,Y)] - \sum_{m\in[M]}\lambda_m \varepsilon_m, \quad \text{s.t. } \Wscr_d(\mu, \hat\mu_m) \leq \varepsilon_m, ~ \forall m \in [M] \bigg\} \nonumber \\
& = \sup_{\mu \in \Pscr(\Yscr)} \bigg\{\E_{Y\sim\mu}[c(z,Y)] - \sum_{m\in[M]}\lambda_m \Wscr_d(\hat\mu_m, \mu) \bigg\}, \label{eq:extended-dual-step1}
\end{align}
which follows from the maximization over $\bm \varepsilon \geq 0$.
By substituting definition \eqref{eq:dwasserstein-dist} for $\Wscr_d$, we have
\begin{align}
\eqref{eq:extended-dual-step1}
& = \sup_{\mu \in \Pscr(\Yscr)} \bigg\{\E_{Y\sim\mu}[c(z,Y)] - \sum_{m\in[M]}\lambda_m \inf_{\zeta_m \in \Xi(\hat\mu_m, \mu)} \E_{(\tilde{Y}_m, Y) \sim \zeta_m} [d(\tilde{Y}_m, Y)] \bigg\} \nonumber \\
& = \sup_{\mu \in \Pscr(\Yscr), ~\zeta_m \in \Xi(\hat\mu_m, \mu)} \E_{Y\sim\mu, ~(\tilde{Y}_m, Y) \sim \zeta_m} \bigg[c(z,Y) - \sum_{m\in[M]}\lambda_m d(\tilde{Y}_m, Y) \bigg]. \label{eq:extended-dual-step2}
\end{align}

In \eqref{eq:extended-dual-step2}, each $\zeta_m \in \Xi(\hat\mu_m, \mu)$ represents a joint distribution with marginals $\hat\mu_m$ and $\mu$.
Hence, the supremum program is effectively taken over a collective joint distribution $(\hat\mu_1,...,\hat\mu_M, \mu)$.
Given $\xi \in \Xi(\hat\mu_1,...,\hat\mu_M)$, let $\Xi(\xi, \Pscr(\Yscr))$ denote the set of all joint distributions on $\Yscr^{M} \times \Yscr$ such that the first $M$ marginals are $\hat\mu_1,...,\hat\mu_M$, respectively.
We note that the following two distribution spaces are equivalent
\begin{equation}
\big\{(\hat\mu_1,...,\hat\mu_M, \mu) ~ \big| ~\forall \mu \in \Pscr(\Yscr), ~\zeta_m \in \Xi(\hat\mu_m, \mu), \forall m \in [M] \big\} = \big\{ \Xi(\xi, \Pscr(\Yscr)) ~ \big| ~ \forall \xi \in \Xi(\hat\mu_1,...,\hat\mu_M) \big\}. \nonumber
\end{equation}
We then reformulate \eqref{eq:extended-dual-step2} with this distribution equivalence as follows
\begin{align}
(-J_{z})^*(-\bm{\lambda})
& = \sup_{\xi \in \Xi(\hat\mu_1,...,\hat\mu_M)} ~\sup_{\zeta \in \Xi(\xi, \Pscr(\Yscr))} \E_{(\tilde{Y}_1,...,\tilde{Y}_M, Y) \sim \zeta} \bigg[c(z,Y) - \sum_{m\in[M]}\lambda_m d(\tilde{Y}_m, Y) \bigg].
\label{eq:extended-dual-step3}
\end{align}

Next, we introduce the interchangeability principle, which allows us to further simplify \eqref{eq:extended-dual-step3}.
\begin{definition}[Interchangeability Principle]
A measurable function $h: \Yscr_1 \times \Yscr_2 \to \R$
and a distribution $\xi \in \Pscr(\Yscr_1)$ satisfy
\begin{equation}\label{eq:interchangeability-principle}
\sup_{\zeta \in \Xi(\xi, \Pscr(\Yscr_2))} \E_{(\tilde{Y}, Y) \sim \zeta} [h(\tilde{Y}, Y)] = \E_{\tilde{Y} \sim \xi} \bigg[ \sup_{y \in \Yscr_2} h(\tilde{Y}, y)\bigg].
\end{equation}
\end{definition}
Given any cost function $c(z,Y)$ satisfying \Cref{asp:reg-costfunc}, for a decision $z$,
the interchangeability principle \eqref{eq:interchangeability-principle} holds for $h_z(\tilde{\bm{Y}}, Y) := c(z,Y) - \sum_{m\in[M]}\lambda_m d(\tilde{Y}_m, Y)$, where $\tilde{\bm{Y}} = (\tilde{Y}_1,..., \tilde{Y}_M)$,
together with any discrete distribution $\xi \in \Xi(\hat\mu_1,...,\hat\mu_M)$.
This result follows from Example 1 in \cite{zhang2022simple}.
We apply the interchangeability principle to \eqref{eq:extended-dual-step3} and obtain
\begin{align}
(-J_{z})^*(-\bm{\lambda})
& = \sup_{\xi \in \Xi(\hat\mu_1,...,\hat\mu_M)} \E_{(\tilde{Y}_1,...,\tilde{Y}_M) \sim \xi} \bigg[ \sup_{y \in \Yscr} \bigg\{c(z,Y) - \sum_{m\in[M]}\lambda_m d(\tilde{Y}_m, Y)\bigg\} \bigg].
\label{eq:extended-dual-step4}
\end{align}

Finally, we apply Legendre transform again to restore $J_{z}(\bm{\varepsilon})$ as follows
\begin{align}
J_{z}(\bm{\varepsilon}) & = -(-J_{z})^{**}(\bm{\varepsilon})
= -\sup_{\bm{\lambda} \geq \bm{0}} \{ -\bm{\lambda}^{\top}\bm \varepsilon - (-J_{z})^*(-\bm{\lambda})\}
= \inf_{\bm{\lambda} \geq \bm{0}} \{ \bm{\lambda}^{\top}\bm \varepsilon + (-J_{z})^*(-\bm{\lambda})\}.
\label{eq:extended-dual-step5}
\end{align}
By substituting \eqref{eq:extended-dual-step4} into \eqref{eq:extended-dual-step5}, we have proved \eqref{eq:strong-dual-extend-1} that
\begin{align}
J_{z}(\bm{\varepsilon}) = & \inf_{\lambda_m \geq 0}\Bigg\{\sum_{m=1}^{M}\lambda_m \varepsilon_m +
\sup_{\xi \in \Xi(\hat\mu_1,...,\hat\mu_M)}
\E_{(\tilde{Y}_1,..., \tilde{Y}_M)\sim \xi}
\bigg[\sup_{y \in \Yscr}\bigg\{c(z, y) -
\sum_{m=1}^{M} \lambda_m d(\tilde Y_m, y) \bigg\}\bigg] \Bigg\}. \quad \hfill \square \nonumber
\end{align}

\subsection{Derivation of Tractable Formulations}\label{app:tractable-formulation}

\noindent
\textbf{Proof of \Cref{coro:cvx}}

Given $c(z, y) = \max_{s \in [S]}\{c_s(z, y)\}$ and $p=1$, the constraints in program \eqref{eq:strong-dual-further} are equivalent to
\begin{align}
\label{eq:reformulation-cvx-constraint-1}
u_i + v_j \geq \sup_{y \in \Yscr} \Big\{c_s(z, y) - \lambda_1 \|y_{1,i} - y\|_1 - \lambda_2 \|y_{2,j} - y\|_1 \Big\}, \quad \forall i,j \in [n], s \in [S].
\end{align}

For a tuple $(i,j,s)$, the right-hand side of \eqref{eq:reformulation-cvx-constraint-1} can be reformulated using $\alpha_{i,s}, \beta_{j,s}  \in \R^{d_y}$
\begin{align}
\text{RHS of } \eqref{eq:reformulation-cvx-constraint-1}
& = \sup_{y \in \Yscr} \min_{\substack{\|\alpha_{i,s}\|_{\infty} \leq \lambda_1\\ \|\beta_{j,s}\|_{\infty} \leq \lambda_2}} \paran{c_{s}(z, y) - \alpha_{i,s}^{\top}(y - y_{1,i}) - \beta_{j,s}^{\top}(y - y_{2,j})} \nonumber  \\
& = \min_{\substack{\|\alpha_{i,s}\|_{\infty} \leq \lambda_1\\ \|\beta_{j,s}\|_{\infty} \leq \lambda_2}} \sup_{y \in \Yscr} \paran{c_{s}(z, y) - \alpha_{i,s}^{\top}(y - y_{1,i}) - \beta_{j,s}^{\top}(y - y_{2,j})} \nonumber \\
& = \min_{\substack{\|\alpha_{i,s}\|_{\infty} \leq \lambda_1\\ \|\beta_{j,s}\|_{\infty} \leq \lambda_2}} \bigg\{\sup_{y \in \Yscr} \Big\{c_{s}(z, y) - (\alpha_{i,s} + \beta_{j,s})^{\top}y \Big\} + \alpha_{i,s}^{\top}y_{1,i} + \beta_{j,s}^{\top} y_{2,j}\bigg\} \nonumber \\
& = \min_{\substack{\|\alpha_{i,s}\|_{\infty} \leq \lambda_1\\ \|\beta_{j,s}\|_{\infty} \leq \lambda_2}} \Big\{ [-c_{s}(z, \cdot) + \sigma_{\Yscr}]^*(\alpha_{i,s} + \beta_{j,s}) - \alpha_{i,s}^{\top}y_{1,i} - \beta_{j,s}^{\top}y_{2,j}\Big\}. \label{eq:reformulation-cvx-constraint-2}
\end{align}
Here, the first equality follows from the H\"older inequality for the $1$-norm.
The second equality is due to Sion's minimax theorem, as both $\Yscr$ and the set $\{\|\alpha_{i,s}\|_{\infty} \leq \lambda_1, \|\beta_{j,s}\|_{\infty} \leq \lambda_2\}$ are compact.
The fourth equality follows directly from the property of the convex conjugate.
We also have
\begin{align}
[-c_{s}(z,\cdot) + \sigma_{\Yscr}]^*(\alpha_{i,s} + \beta_{j,s}) = \inf_{{\zeta_{i,j,s}} \in \R^{d_y}} \Big\{[-c_s(z,\cdot)]^*(\alpha_{i,s} + \beta_{j,s} - \zeta_{i,j,s}) + \sigma_{\Yscr}({\zeta_{i,j,s}}) \Big\}.
\label{eq:reformulation-cvx-constraint-3}
\end{align}

Finally, by substituting \eqref{eq:reformulation-cvx-constraint-3} into \eqref{eq:reformulation-cvx-constraint-2}
and replacing constraints in program \eqref{eq:strong-dual-further},
we establish the following constraints for the convex program \eqref{eq:reformulation-cvx-obj}.
$\hfill \square$
\begin{numcases}{}
u_i + v_j \geq [-c_s(z,\cdot)]^*(\alpha_{i,s} + \beta_{j,s} - \zeta_{i,j,s}) + \sigma_{\Yscr}(\zeta_{i,j,s}) - \alpha_{i,s}^{\top}y_{1,i} - \beta_{j,s}^{\top}y_{2,j}, ~~ \forall i,j \in [n], s \in [S]; \label{eq:reformulation-cvx-constraint-4} \\
\|\alpha_{i,s}\|_{\infty} \leq \lambda_1, ~~\|\beta_{i,s}\|_{\infty} \leq \lambda_2, ~~ \alpha_{i,s}, \beta_{i,s}  \in \R^{d_y}, \quad \forall i \in [n], s \in [S]; \nonumber \\
\lambda_1, \lambda_2 \geq 0;
\quad u_i, v_i \in \R, ~~\forall i \in [n];
\quad \zeta_{i,j,s} \in \R^{d_y}, ~~\forall i, j \in [n], s \in [S]. \nonumber
\end{numcases}

\smallskip
\noindent
\textbf{Proof of \Cref{coro:linear}}

Given that $c(z,y)= \max_{s \in [S]}\{a_{s}(z)^{\top}y + b_{s}(z)\}$ is piecewise linear
($\forall s \in [S]$, both $a_{s}(z):\Zscr\rightarrow\R^{d_y}$ and $b_{s}(z):\Zscr\rightarrow\R$ are linear in $z$)
and $\Yscr = \{y \in \R^{d_y}: Ay \leq h \}$ is a polyhedron,
we can further simplify constraint \eqref{eq:reformulation-cvx-constraint-4} into linear ones.
Specifically, by applying the definition of conjugate function and strong duality, we have
\begin{align}
[-c_{s}(z,\cdot)]^*(\alpha_{i,s} + \beta_{j,s} - \zeta_{i,j,s})
& = [-a_s(z)^{\top}(\cdot) - b_s(z)]^*(\alpha_{i,s} + \beta_{j,s} - \zeta_{i,j,s}) \nonumber \\
& =\begin{cases}
b_s(z), & \text{if } \alpha_{i,s} + \beta_{j,s} - \zeta_{i,j,s}= -a_s(z), \\
\infty, & \text{otherwise.}
\end{cases} \label{eq:reformulation-linear-constraint-1}
\end{align}
\begin{align}
\sigma_{\mathcal{Y}}(\zeta_{i,j,s}) = \sup_{y} \Big\{\zeta_{i,j,s}^{\top} y, \text{~~s.t.~} Ay \leq h \Big\} = \inf_{\gamma_{i,j,s} \in \R^{d_m}_{+}} \Big\{\gamma_{i,j,s}^{\top}h, \text{~~s.t.~} A^{\top} \gamma_{i,j,s} = \zeta_{i,j,s}\Big\}. \label{eq:reformulation-linear-constraint-2}
\end{align}

The second case of \eqref{eq:reformulation-linear-constraint-1} results in an infeasible program.
For the first case,
by substituting \eqref{eq:reformulation-linear-constraint-1} and \eqref{eq:reformulation-linear-constraint-2} into \eqref{eq:reformulation-cvx-constraint-4}
and eliminating $\zeta_{i,j,s}$,
we establish the first two constraints for the linear program \eqref{eq:reformulation-linear-obj} as follows:
$\hfill \square$
\begin{numcases}{}
u_i + v_j \geq b_s(z) + \gamma_{i,j,s}^{\top}h - \alpha_{i,s}^{\top}y_{1,i} - \beta_{j,s}^{\top} y_{2,j}, \quad \forall i,j \in [n], s \in [S]; \nonumber \\
\alpha_{i,s} + \beta_{j,s} = A^{\top}\gamma_{i,j,s} - a_{s}(z), \quad \forall i,j \in [n], s \in [S]. \nonumber
\end{numcases}

\smallskip
\noindent
\textbf{Results for the Extended {\iwdro}}

Similar to \Cref{coro:cvx},
given a piecewise concave cost function,
we can derive a convex formulation for the extended {\iwdro} problem with $\setag$, which is formed by multiple Wasserstein balls.
\begin{corollary}[Convex Reformulation of the Extended {\iwdro}]\label{coro:cvx-extend}
Let $d^*$ denote the convex conjugate of $d(y_1, y_2)$ with respect to the first component $y_1$.
If $c(z,y) = \max_{s \in [S]} \{c_s(z, y)\}$
and each $c_s(z, y)$ is concave and upper semi-continuous in $y \in \Yscr$,
the optimization program $\JD(z)$ given in \eqref{eq:strong-dual-extend-2} (or equivalently $\JP(z)$) can be reformulated as the following convex program:
\begin{align}
\label{eq:reformulation-general-cvx-obj}
\inf ~ & \sum_{m\in[M]} \lambda_m \varepsilon_m + \sum_{m\in[M]} \sum_{i\in[n]} w_{m,i} u_{m,i}, \\
\text{\textup{s.t.}} ~
& \sum_{m\in[M]} u_{m,\bm{i}_m} \geq
[-c_s(z,\cdot)]^* \bigg(\sum_{m \in[M]} \alpha_{m, \bm{i}_m, s} - \zeta_{\bm{i},s}\bigg) + \sigma_{\Yscr}(\zeta_{\bm{i},s})
- \sum_{m \in[M]} \lambda_m d^* \bigg(\frac{\alpha_{m,\bm{i}_m,s}}{\lambda_m}, y_{m, \bm{i}_m}\bigg),
~\forall \bm{i} \in \Iscr_{M, n},s \in [S]; \nonumber \\
& \lambda_m \geq 0, \quad u_{m,\bm{i}_m} \in \R, \quad \alpha_{m,\bm{i}_m,s} \in \R^{d_y}, \quad \zeta_{\bm{i},s} \in \R^{d_y}, \quad \forall m \in [M], \bm{i} \in \Iscr_{M, n},s \in [S]. \nonumber
\end{align}
\end{corollary}

When the $p$-Wasserstein distance $\Wscr_p$ is used,
the distance metric is given by $d(y_1, y_2) = \|y_1 - y_2\|_1^p$,
which results in $d^*(y_1, y_2) = y_1^{\top}y_2 + \frac{p - 1}{p^{\frac{p}{p - 1}}}\|y_1\|_1^{\frac{p}{p - 1}}$.
Hence, the first constraint in program \eqref{eq:reformulation-general-cvx-obj} can be expressed as
\begin{align}
\sum_{m\in[M]} u_{m,\bm{i}_m} \geq
[-c_s(z,\cdot)]^* \bigg(\sum_{m \in[M]} \alpha_{m, \bm{i}_m, s} - \zeta_{\bm{i},s}\bigg) + \sigma_{\Yscr}(\zeta_{\bm{i},s})
- \sum_{m \in[M]} \bigg(\alpha_{m, \bm{i}_m, s}^{\top} y_{m,\bm{i}_m} + \frac{p - 1}{p^{\frac{p}{p - 1}}} \lambda_m \Big\|\frac{\alpha_{m,\bm{i}_m,s}}{\lambda_m} \Big\|_1^{\frac{p}{p - 1}}\bigg). \nonumber
\end{align}
Moreover, similar to \Cref{coro:linear}, when the cost function is piecewise linear,
we can further simplify program \eqref{eq:reformulation-general-cvx-obj}.

\begin{corollary}[Simplified Reformulation of the Extended {\iwdro}]\label{coro:linear-extend}
If $\Wscr_p$ is used in $\setag$,
$\displaystyle{c(z,y)= \max_{s \in [S]}\{a_{s}(z)^{\top}y + b_{s}(z)\}}$ where $a_{s}(z)$ and $b_{s}(z)$ are linear in $z$,
\hspace{-0.15cm} and $\Yscr = \{y \in \R^{d_y}: Ay \leq h \}$
for some $A \in \R^{d_m \times d_y}$ and $h \in \R^{d_m}$,
then program \eqref{eq:reformulation-general-cvx-obj} simplifies to
\begin{align}
\label{eq:reformulation-general-linear-obj}
\inf ~ & \sum_{m\in[M]} \lambda_m \varepsilon_m + \sum_{m\in[M]} \sum_{i\in[n]} w_{m,i} u_{m,i}, \\
\text{\textup{s.t.}} ~
& \sum_{m\in[M]} u_{m,\bm{i}_m} \geq
b_s(z) + \gamma_{i,j,s}^{\top}h
- \sum_{m \in[M]} \bigg(\alpha_{m, \bm{i}_m, s}^{\top} y_{m,\bm{i}_m} + \frac{p - 1}{p^{\frac{p}{p - 1}}} \lambda_m \Big\|\frac{\alpha_{m,\bm{i}_m,s}}{\lambda_m} \Big\|_1^{\frac{p}{p - 1}}\bigg), ~\forall \bm{i} \in \Iscr_{M, n},s \in [S]; \nonumber \\
& \sum_{m \in[M]} \alpha_{m, \bm{i}_m, s} = A^{\top}\gamma_{i,j,s} - a_{s}(z), \quad \forall \bm{i} \in \Iscr_{M, n},s \in [S]; \nonumber \\
& \lambda_m \geq 0, \quad u_{m,\bm{i}_m} \in \R, \quad \alpha_{m,\bm{i}_m,s} \in \R^{d_y}, \quad \gamma_{i,j,s} \in \R^{d_m}_{+}, \quad \forall m \in [M], \bm{i} \in \Iscr_{M, n},s \in [S]. \nonumber
\end{align}
\end{corollary}

\smallskip
\textbf{Remark}:
For the extended {\iwdro}, we have studied the DRO problem with an intersection ambiguity set formed by multiple Wasserstein balls.
Our results (i.e., \Cref{thm:strong-duality3} and Corollaries \ref{coro:cvx-extend} and \ref{coro:linear-extend})
yield the same problem formulation as the concurrent work \cite{rychener2024wasserstein}.

\section{Proofs in Section~\ref{sec:statistical-result}}\label{app:proof-concentration}

We first introduce some technical lemmas in \ref{app:proof-concentration-lemmas},
followed by the proofs of Theorem~\ref{thm:concentration-np} in \ref{app:proof-ws-concentrate}.
Additionally, we extend the concentration result to the $k$-nearest neighbor ($k$NN) method in \ref{app:kNN}.

\subsection{Technical Lemmas}\label{app:proof-concentration-lemmas}

\begin{lemma}[Multiplicative Chernoff Bound]\label{lemma:multi-cb}
Let $U_1,...,U_n$ be i.i.d. Bernoulli random variables.
Denote $U = \sum_{i = 1}^n U_i$ and $\mu = \E[U]$.
Then, for all $0 < \delta \le 1$, we have:
\begin{align}
\P(U \leq (1-\delta)\mu) \leq \exp(-\delta^2 \mu /2), \quad \P( U \ge (1+\delta)\mu) \leq \exp(-\delta^2 \mu /3). \nonumber
\end{align}
\end{lemma}

\smallskip
Recall that, in \Cref{thm:concentration-np}, we choose the bandwidth parameter $h_n = C_h n^{-\psi}$  with $\psi \in (0, d_x^{-1})$ and require $n \ge C_{\mu} \mu_X(x)^{-\frac{1}{\psi d_x}}$, which ensures $h_n^{d_x} \le C_h^{d_x} C_{\mu}^{-\psi d_x} \mu_X(x)$.
Since we focus on the asymptotic results relates to $n$, we write $h_n = \Theta(n^{-\psi})$ and $h_n^{d_x}=O(\mu_X(x))$
for simplicity in the proof.

Next, we introduce a key quantity in the proof of Theorem~\ref{thm:concentration-np}.
Given training samples $\{(x_i, y_i)\}_{i = 1}^{n}$,
we define the number of samples near the test covariate $x$ as:
\begin{align}
\label{eq:number-near-sample}
N_r := \sum_{i = 1}^n \I{\|x - x_i\|_2 \leq r h_n},
\end{align}
where $r>0$ is the constant assumed in \Cref{asp:kernel} and $h_n$ is the bandwidth parameter.
Using this definition, we can establish the following tail estimation results.
\begin{lemma}[Statistical Properties of Local Estimators]\label{lemma:lower-bd-Nnr}
There exists some constant $\tilde{C}_0$ independent of $n$ and $x$, such that
\begin{align}
\label{eq:lemma-number-near-sample-2}
& \P(N_r\le \tilde{C}_0 n \mu_X(x) h_n^{d_x}/4) \le \exp(-\tilde{C}_0 n\mu_X(x) h_n^{d_x} / 8), \\
\label{eq:lemma-number-near-sample-3}
& \P(N_r \ge 4\tilde{C}_0  n \mu_X(x) h_n^{d_x}) \le \exp(-\tilde{C}_0 n \mu_X(x) h_n^{d_x} / 8).
\end{align}
\end{lemma}
\emph{Proof: }
For any measurable set $\Bscr \subseteq \R^d$, its volume is defined as $\text{Vol}(\Bscr) \coloneqq \int_{\R^d} \I{x\in\Bscr} dx$,
where the integration is taken with respect to the Lebesgue measure on $\R^d$.
Given a test covariate $x \in \Xscr$, we consider the set $\Bscr = \{\tilde{x}: \|x-\tilde{x}\|_2\le rh_n\}$.
If $\Bscr \subseteq \Xscr$, the volume of this $d_x$-dimensional ball is given by
$\text{Vol}(\Bscr) = \frac{\pi^{\frac{d_x}{2}}}{\Gamma(\frac{d_x}{2} + 1)} r^{d_x} h_n^{d_x}$.
We denote $\tilde{C}_0 = \frac{\pi^{\frac{d_x}{2}}}{\Gamma(\frac{d_x}{2} + 1)} r^{d_x}$, which is independent of $n$ and $x$.
The probability that a sample covariate $x_i \in \Bscr$ is bounded by
\begin{align*}
\P(x_i \in \Bscr) = \int_{\tilde{x} \in \Bscr}\mu_X(\tilde{x}) d\tilde{x}
\in \bigg[\inf_{\tilde{x} \in \Bscr}\mu_X(\tilde{x})\text{Vol}(\Bscr), ~\sup_{\tilde{x} \in \Bscr}\mu_X(\tilde{x})\text{Vol}(\Bscr) \bigg].
\end{align*}

Given $h_n^{d_x}=O(\mu_X(x))$ and the continuity of $\mu_X(x)$ in \Cref{asp:regular_marginal_dist},
we can select a small $h_n$ such that,
for any $\tilde{x} \in \Bscr$ satisfying $\|x-\tilde{x}\|_2 \leq r h_n$, we have $\mu_X(\tilde{x}) = \mu_X(x) + O(h_n^{d_x}) \in [\mu_X(x)/2, 2\mu_X(x)]$.
This leads to $\P(x_i \in \Bscr) \in [\tilde{C}_0 \mu_X(x) h_n^{d_x}/2, ~2\tilde{C}_0 \mu_X(x) h_n^{d_x}]$.
Because $N_r$, as defined in \eqref{eq:number-near-sample}, is the sum of $n$ i.i.d. Bernoulli random variables with success probability $\P(x_i \in \Bscr)$,
we obtain:
\begin{align}\label{eq:lemma-number-near-sample-1-range}
\E[N_r] = n \P(x_i \in \Bscr) \in [\tilde{C}_0 n\mu_X(x) h_n^{d_x}/2, ~2\tilde{C}_0 n\mu_X(x) h_n^{d_x}].
\end{align}

Finally, by applying Lemma~\ref{lemma:multi-cb} and using $\E[N_r] \ge \tilde{C}_0 n\mu_X(x) h_n^{d_x}/2$, we arrive at \eqref{eq:lemma-number-near-sample-2}
$$\P(N_r\le \tilde{C}_0 n \mu_X(x) h_n^{d_x}/4) = \P(N_r \le (1-1/2)\E[N_r]) \le \exp\left(-\tilde{C}_0 n\mu_X(x) h_n^{d_x} / 8\right).$$
Similarly, by using $\E[N_r] \le 2\tilde{C}_0 n\mu_X(x) h_n^{d_x}$, we obtain \eqref{eq:lemma-number-near-sample-3} $\hfill \square$
$$\P(N_r\ge 4 \tilde{C}_0 n \mu_X(x) h_n^{d_x}) = \P(N_r \ge (1+1)\E[N_r]) \le \exp\left(-\tilde{C}_0 n\mu_X(x) h_n^{d_x} / 3\right) \le \exp\left(-\tilde{C}_0 n\mu_X(x) h_n^{d_x} / 8\right).$$

Next, we present Lemmas~\ref{lemma:was-dual} - \ref{lemma:was-convexity}, which pertain to the statistical properties of the Wasserstein distance.
\begin{lemma}[Duality of the Wasserstein Distance]\label{lemma:was-dual}
Let $f: \Yscr \rightarrow \R$ and $g: \Yscr \rightarrow \R$ denote two functions.
The $p$-Wasserstein distance between two distributions $\mu$ and $\nu$, as defined in \eqref{eq:pwasserstein-dist}, can be computed by:
\begin{equation}\label{eq:generaldualwasserstein}
\Wscr_p^p(\mu, \nu) = \sup_{f,g}\left\{\int_{\Yscr} f(y_1)d\mu(y_1) + \int_{\Yscr} g(y_2)d\nu(y_2), \text{~~s.t.~~}  f(y_1) + g(y_2)\leq \|y_1 - y_2\|_1^p \right\}.
\end{equation}
\end{lemma}
\emph{Proof: }
This result is from Theorem 5.10 in \cite{kantorovich1958space}.
$\hfill \square$

Recall that Assumption~\ref{asp:regular_dist} specifies the data generating process $y = f(x) + \eta$.
Given training samples $\{(x_i, y_i)\}_{i  \in [n]}$,
let $\eta_i = y_i-f(x_i)$ denote the realization of noise for sample $i$.
For a test covariate $x$, we define a discrete reference distribution with supports $\{f(x) + \eta_i\}_{i \in [n]}$ as:
\begin{equation}\label{eq:mur}
\mur = \sum_{i = 1}^n w_i \delta_{\{f(x) + \eta_i\}}.
\end{equation}
We use $\mur$ as an intermediate measure in analyzing $\Wscr_p(\mu_{Y|x},\munp)$.
When choosing equal weights $w_i = 1/n$ for $\mur$ in \eqref{eq:mur}, we denote the resulting uniform distribution as: $\muru = \sum_{i = 1}^n \frac{1}{n} \delta_{\{f(x) + \eta_i\}}.$

\begin{lemma}\label{lemma:kernel-naive-relationship}
Under \Cref{asp:regular_dist}, it holds that:
\begin{equation}\label{eq:kernel-naive-relationship}
\Wscr_p(\mu_{Y|x},\mur) \leq \Big(n \max_{i \in [n]}\{w_i\}\Big)^{\frac{1}{p}}\Wscr_p(\mu_{Y|x},\muru).
\end{equation}
\end{lemma}
\emph{Proof: }
By \Cref{lemma:was-dual},
let $f^*$ and $g^*$ denote the optimal solutions to the dual formulation \eqref{eq:generaldualwasserstein} for $\Wscr_p^p(\mu_{Y|x},\mur)$.
It follows that $f^*(y_1) + g^*(y_2) \leq \|y_1 - y_2\|_1^p$ for all $y_1, y_2 \in \Yscr$, and:
\begin{align}
\Wscr_p^p(\mu_{Y|x},\mur) &= \int_{\Yscr}f^*(y_1)\mu_{Y|x}(d y_1) + \int_{\Yscr}g^*(y_2)\mur(dy_2) \nonumber \\
& = \int_{\Yscr}f^*(y_1)\mu_{Y|x}(d y_1) + \sum_{i=1}^{n} w_i g^*(f(x) + \eta_i) \nonumber \\
& = \sum_{i=1}^{n} w_i \left(\int_{\Yscr}f^*(y_1)\mu_{Y|x}(d y_1) + g^*(f(x) + \eta_i) \right) \nonumber \\
& \leq n\max_{i \in [n]}\{w_i\} \sum_{i = 1}^n \frac{1}{n} \left(\int_{\Yscr}f^*(y_1)\mu_{Y|x}(d y_1) + g^*(f(x) + \eta_i) \right). \label{eq:kernel-naive-relationship-step1}
\end{align}

It is easy to see that $f^*$ and $g^*$ are feasible to the dual formulation \eqref{eq:generaldualwasserstein} for $\Wscr_p^p(\mu_{Y|x},\muru)$,
and thus we have
\begin{align}
\Wscr_p^p(\mu_{Y|x},\muru) & \ge \int_{\Yscr}f^*(y_1)\mu_{Y|x}(d y_1) + \sum_{i=1}^{n} \frac{1}{n} g^*(f(x) + \eta_i) \nonumber \\
& = \sum_{i = 1}^n \frac{1}{n} \left(\int_{\Yscr}f^*(y_1)\mu_{Y|x}(d y_1) + g^*(f(x) + \eta_i) \right). \label{eq:kernel-naive-relationship-step2}
\end{align}
Substituting \eqref{eq:kernel-naive-relationship-step2} into \eqref{eq:kernel-naive-relationship-step1}, we obtain $\Wscr_p(\mu_{Y|x},\mur) \leq \Big(n \max_{i \in [n]}\{w_i\}\Big)^{\frac{1}{p}}\Wscr_p(\mu_{Y|x},\muru).$
$\hfill \square$

\begin{lemma}\label{lemma:was-convexity}
Let $\mu_0 = \sum_{s \in [S]} \lambda_s \mu_s$ denote a convex combination of the measures $\{\mu_s\}_{s \in [S]}$,
where $0 \le \lambda_s \le 1$ and $\sum_{s \in [S]} \lambda_s = 1$.
Then, for any measure $\mu$, we have:
\begin{equation}\label{eq:was-convexity}
\Wscr_p^p(\mu, \mu_0) \leq \sum_{s \in [S]} \lambda_s \Wscr_p^p(\mu, \mu_s).
\end{equation}
\end{lemma}
\emph{Proof: }
Let $\xi_s \in \Xi(\mu, \mu_s)$ be the joint measure that attains the infimum
in the definition of $\Wscr_p(\mu, \mu_s)$,
i.e., $\Wscr_p^p(\mu, \mu_s) = \int \|y_1 - y_2\|_1^p \xi_s(dy_1,dy_2).$
It is straightforward to verify that the convex combination $\xi_0 = \sum_{s \in [S]}\lambda_s\xi_s$ satisfies $\xi_0 \in \Xi(\mu, \mu_0)$.
Then, we have:
\begin{align*}
\Wscr_p^p(\mu, \mu_0) & \leq \int \|y_1 - y_2\|_1^p \xi_0(dy_1, dy_2)
= \sum_{s \in [S]}\lambda_s \int \|y_1 - y_2\|_1^p \xi_s(dy_1, dy_2) = \sum_{s \in [S]} \lambda_s \Wscr_p^p(\mu, \mu_s). \quad  \square
\end{align*}

\subsection{Proof of Theorem~\ref{thm:concentration-np}.}\label{app:proof-ws-concentrate}

To upper bound $\Wscr_p(\mu_{Y|x},\munp)$, we use the reference distribution $\mur$ defined in \eqref{eq:mur},
with weights assigned according to the NW kernel estimator $\munp$, as given in \eqref{eq:ball-center-np}.
Specifically, we define
\begin{align}
\label{eq:mur-kenel}
\mur  = \sum_{i = 1}^n\frac{K((x - x_i)/h_n)}{\sum_{j = 1}^n K((x-x_j)/h_n)}\delta_{\{f(x)+\eta_i\}},
\end{align}
and use it as an intermediate
in the triangle inequality $\Wscr_p(\mu_{Y|x},\munp) \leq \Wscr_p(\mu_{Y|x}, \mur) + \Wscr_p(\mur, \munp)$.
This enables us to prove \Cref{thm:concentration-np} by applying the following probabilistic bound:
$$\P \big(\Wscr_p(\mu_{Y|x},\munp)>t + C_0 Lh_n^{\beta}\big) \le \P\big(\Wscr_p(\mu_{Y|x},\mur)>t\big)+ \P \big(\Wscr_p(\mur,\munp)> C_0 Lh_n^{\beta}\big).$$

The introduction of $\mur$ facilitates the analysis.
When analyzing $\Wscr_p(\mu_{Y|x},\mur)$,
the support of $\mur$, \\
i.e, $\{f(x) + \eta_i\}_{i \in [n]}$, can be interpreted  as i.i.d. samples from $\mu_{Y|x}$.
When analyzing $\Wscr_p(\mur, \munp)$, $\mur$ shares the same probability mass as $\munp$.
Next, we analyze these two terms separately.

\smallskip
\noindent
$\bullet~$ \textbf{Part 1.} For $\Wscr_p(\mu_{Y|x}, \mur)$, we consider the following three cases of the kernel function.

\underline{Case (1)}: A naive kernel function $K(\tau) = \I{\|\tau\|_2\le r}$.
By definition \eqref{eq:mur-kenel}, $\mur$ takes the form
\begin{align}
\label{eq:mur-kenel-naive}
\mur = \muru := \sum_{i=1}^{n} \frac{ \I{\|  x-  x_i\|_2\le rh_n}}{\sum_{j=1}^n \I{\|  x-  x_j\|_2\le rh_n}} \delta_{\{f(x)+\eta_i\}}.
\end{align}

Recall that we defined $N_r = \sum_{i = 1}^n \I{\|x - x_i\|_2 \leq r h_n}$ in \eqref{eq:number-near-sample},
which represents the number of samples near the test covariate $x$.
Under \Cref{asp:regular_dist},
$N_r$ is independent of the noise terms $\{\eta_i\}_{i\in[n]}$,
as the realization of $x$ from the covariate distribution is independent of the noise distribution.
Conditioned on the event $N_r=k$,
$\muru$, as given in \eqref{eq:mur-kenel-naive}, is a uniform distribution over
$k$ points from $\{f(x) + \eta_i\}_{i \in [n]}$,
where the corresponding covariates satisfy $\|x- x_i\|_2\le rh_n$.
By Theorem 2 in \cite{fournier2015rate}, we have the following result.

\begin{lemma}[Conditional Concentration]\label{lemma:empirical-was-concentration}
Given $\muru$ as defined in \eqref{eq:mur-kenel-naive},
for any $ x\in \Xscr$ and $t\in (0,1]$, we have
\begin{equation}\label{eq:condition-naivekernel}
\P(\Wscr_p(\mu_{Y|x}, \muru)>t | N_r = k) \le
\begin{cases} \tilde{C}_1 \exp(-\tilde{C}_2 k t^2), & \text{if } p > d_y / 2, \\
\tilde{C}_1 \exp(-\tilde{C}_2k (t/\log (2 + 1/t))^2), & \text{if } p = d_y / 2, \\
\tilde{C}_1 \exp(-\tilde{C}_2 k t^{d_y/p}), & \text{if } p \in [1, d_y / 2),\\
\end{cases}
\end{equation}
where $\tilde{C}_1$ and $\tilde{C}_2$ are some constants independent of $n$, $x$, and $t$.
\end{lemma}

Given that \eqref{eq:condition-naivekernel} is decreasing in $k$, we have
\begin{align}
\label{eq:condition-naivekernel-1}
\begin{split}
\P(\Wscr_p(\mu_{Y|x}, \muru)>t , N_r > k) & = \sum_{k'>k} \P(\Wscr_p(\mu_{Y|x}, \muru)>t | N_r = k')\P(N_r = k') \\
& \le \P(\Wscr_p(\mu_{Y|x}, \muru)>t | N_r = k).
\end{split}
\end{align}

Recall that $\mur = \muru$ in this case.
When $d_y = 1$, we combine \eqref{eq:condition-naivekernel-1} with the first case of \eqref{eq:condition-naivekernel},
and then apply \eqref{eq:lemma-number-near-sample-2} from \Cref{lemma:lower-bd-Nnr},
yielding the following result:
\begin{align}
\P(\Wscr_p(\mu_{Y|x}, \mur)>t) & \le \P(\Wscr_p(\mu_{Y|x},\mur)>t, N_r> \tilde{C}_0 n \mu_X(x) h_n^{d_x}/4) + \P(N_r \le \tilde{C}_0 n \mu_X(x) h_n^{d_x}/4)  \nonumber \\
& \le \tilde{C}_1 \exp(-\tilde{C}_0 \tilde{C}_2 /4 \cdot n \mu_X(x) h_n^{d_x} t^2) + \exp(-\tilde{C}_0 / 8 \cdot n\mu_X(x) h_n^{d_x}) \nonumber \\
& \le C_1^* \exp(-C_2^* n \mu_X(x) h_n^{d_x} t^2), \label{eq:mur-kenel-naive-case1}
\end{align}
where we set $C_1^* = \tilde{C}_1 + 1$ and $C_2^* = \min\{\tilde{C}_0 \tilde{C}_2/4, ~\tilde{C}_0/8 \}$.

\smallskip
\underline{Case (2)}: A box kernel function $K(\tau)$ that satisfies the assumption below.
\begin{assumption}\label{asp:kernel2}
There exist some constants $0 < r \le R$ and $0< b_r \le b_R$,  such that $K(\tau)$ satisfies
\begin{align}\label{eq:box-kernel}
b_r \I{\| \tau\|_2\le r}\le K(\tau)\le b_R \I{\| \tau\|_2\le R}.
\end{align}
\end{assumption}
Under \Cref{asp:kernel2}, box kernel functions generalize naive kernel functions and constitute a subset of kernel functions that fulfill Assumption~\ref{asp:kernel}.

Using the constant $R$ in \Cref{asp:kernel2}, we define $N_R := \sum_{i = 1}^n \I{\|x - x_i\|_2 \leq R h_n}$.
Since \eqref{eq:box-kernel} implies that $K((x - x_i)/h_n)=0$ for $\|x -x_i\|_2 > R h_n$,
conditioned on the event $N_R = k$,
there are $k$ nonzero weights in $\mur$, as defined in \eqref{eq:mur-kenel}.
Moreover, for these nonzero weights, we have
$$\frac{K((x - x_i)/h_n)}{\sum_{j = 1}^n K((x-x_j)/h_n)} \leq \frac{b_R}{k b_r}, \quad \forall i\in [n], \text{~~s.t.~} \|x - x_i\|_2 \leq R h_n.$$

Following \eqref{eq:mur-kenel-naive},
we redefine $\muru = \sum_{i = 1}^n \frac{\I{\|x - x_i\|_2\leq R h_n}}{\sum_{j = 1}^n \I{\|x - x_j\|_2\leq R h_n}}\delta_{\{f(x) + \eta_i\}}$
using the radius $R$.
Conditioned on $N_R = k$, $\muru$ is a uniform distribution over
$k$ points from $\{f(x) + \eta_i\}_{i \in [n]}$,
where the associated covariates satisfy $\|x- x_i\|_2\le Rh_n$.
Applying Lemma~\ref{lemma:kernel-naive-relationship}, we obtain:
\begin{equation}\label{eq:nonparam-upper}
\Wscr_p(\mu_{Y|x}, \mur)\I{N_R = k} \leq (b_R/b_r)^{\frac{1}{p}} \Wscr_p(\mu_{Y|x}, \muru)\I{N_R = k}, \quad \forall k \in [n].
\end{equation}

It is straightforward to note that: (i) \Cref{lemma:empirical-was-concentration} applies to the redefined $\muru$ with $N_R$,
and (ii) \Cref{lemma:lower-bd-Nnr} applies to $N_R$ with $\tilde{C}_0 = \frac{\pi^{\frac{d_x}{2}}}{\Gamma(\frac{d_x}{2} + 1)} R^{d_x}$.
Following a similar approach to the derivation of \eqref{eq:mur-kenel-naive-case1}, for $d_y = 1$, we have:
\begin{align}
\P(\Wscr_p(\mu_{Y|x}, \mur)>t) & \le \P(\Wscr_p(\mu_{Y|x},\mur)>t, N_R> \tilde{C}_0 n \mu_X(x) h_n^{d_x}/4) + \P(N_R \le \tilde{C}_0 n \mu_X(x) h_n^{d_x}/4)  \nonumber \\
& \le \P(\Wscr_p(\mu_{Y|x},\muru)> (b_r/b_R)^{\frac{1}{p}}t, N_R> \tilde{C}_0 n \mu_X(x) h_n^{d_x}/4) + \P(N_R \le \tilde{C}_0 n \mu_X(x) h_n^{d_x}/4)  \nonumber \\
& \le \tilde{C}_1 \exp(-\tilde{C}_0 \tilde{C}_2 (b_r/b_R)^{\frac{2}{p}}/4 \cdot n \mu_X(x) h_n^{d_x} t^2) + \exp(-\tilde{C}_0/8 \cdot n\mu_X(x) h_n^{d_x}) \nonumber \\
& \le C_1^* \exp(-C_2^* n \mu_X(x) h_n^{d_x} t^2),  \label{eq:mur-kenel-box-case1}
\end{align}
where we set $C_1^* = \tilde{C}_1 + 1$ and $C_2^* = \min\{\tilde{C}_0 \tilde{C}_2 (b_r/b_R)^{\frac{2}{p}} /4, ~\tilde{C}_0/8 \}$.

\smallskip
\underline{Case (3)}: A general kernel function $K(\tau)$ that satisfies \Cref{asp:kernel}.
We sort the training samples $\{(x_i, y_i)\}_{i = 1}^{n}$ by the distance $\|x - x_i\|_2$ in non-decreasing order.
As a result, the first $N_r$ sample covariates satisfy $\|x - x_i\|_2 \leq r h_n$ for $i \in [N_r]$.
We reformulate $\mur$, as defined in \eqref{eq:mur-kenel}, using a convex combination of $n - N_r + 1$ component measures as follows:
\begin{align}
\label{eq:mur-kenel-general}
& \mur  = \sum_{i = 1}^n\frac{K((x - x_i)/h_n)}{\sum_{j = 1}^n K((x-x_j)/h_n)}\delta_{\{f(x)+\eta_i\}} = \sum_{s=1}^{n - N_r + 1} \lambda_s \mur^{(s)},
\end{align}
where $\displaystyle{\mur^{(1)} = \sum_{i \in [N_r]}\frac{K((x - x_i)/h_n)}{\sum_{j \in [N_r]} K((x - x_j)/h_n)} \delta_{\{f(x)+\eta_i\}}}$
and $\displaystyle{\mur^{(s)} = \frac{1}{N_r + s - 1}\sum_{i = 1}^{N_r + s - 1}\delta_{\{f(x)+\eta_i\}}}, ~\forall s \geq 2$.
Given that the value of $\frac{K((x - x_i)/h_n)}{\sum_{j \in [n]} K((x - x_j)/h_n)}$ decreases with the index $i$,
we can always identify a suitable set of coefficients $0 \le \lambda_s \le 1$ that satisfy both \eqref{eq:mur-kenel-general} and $\sum_{s \in [n - N_r + 1]} \lambda_s = 1$.
Consequently, by Lemma~\ref{lemma:was-convexity}, we have:
\begin{align}\label{eq:general-kernel-bd}
\Wscr_p^p(\mu_{Y|x}, \mur) \leq \sum_{s = 1}^{n - N_r + 1} \lambda_s\Wscr_p^p(\mu_{Y|x}, \mur^{(s)}),
\end{align}
which implies that if $\Wscr_p(\mu_{Y|x}, \mur^{(s)}) < t$ for each $s =1,...,n - N_r + 1$, then $\Wscr_p(\mu_{Y|x}, \mur) < t$.

Next, we investigate each  $\Wscr_p(\mu_{Y|x}, \mur^{(s)})$ conditioned on the value of $N_r$.
For $\mur^{(1)}$, we have
$$\frac{K((x - x_i)/h_n)}{\sum_{j\in[N_r]} K((x-x_j)/h_n)} \leq \frac{b_R}{k b_r}, \quad \forall i\in [N_r],$$
which is similar to Case (2) of the boxed kernel.
For $d_y = 1$,
following \eqref{eq:mur-kenel-box-case1}, we could obtain
\begin{align}\label{eq:general-kernel-mu1}
\P(\Wscr_p(\mu_{Y|x}, \mur^{(1)}) > t | N_r) \leq C_1^{(1)} \exp(-C_2^{(1)} N_r t^2).
\end{align}
For $s\ge2$, each $\mur^{(s)}$ is a uniform distribution over $N_r + s - 1$ points. Using \Cref{lemma:empirical-was-concentration}, we have
\begin{align}\label{eq:general-kernel-mu2}
\P(\Wscr_p(\mu_{Y|x}, \mur^{(s)}) > t | N_r)\leq C_1^{(s)} \exp(-C_2^{(s)} (N_r + s - 1) t^2) \leq C_1^{(s)} \exp(-C_2^{(s)} N_r t^2).
\end{align}

By combining \eqref{eq:general-kernel-mu1} with \eqref{eq:general-kernel-mu2}, we have:
\begin{align}
\P(\Wscr_p(\mu_{Y|x}, \mur) > t)
& \leq \sum_{s=1}^{n - N_r + 1} \P(\Wscr_p(\mu_{Y|x}, \mur^{(s)}) > t, N_r > \tilde{C}_0 n \mu_X(x) h_n^{d_x}/4) + \P(N_r \le \tilde{C}_0 n \mu_X(x) h_n^{d_x}/4) \nonumber \\
& \le \sum_{s=1}^{n - N_r + 1} C_1^{(s)} \exp(-\tilde{C}_0 C_2^{(s)} /4 \cdot n \mu_X(x) h_n^{d_x} t^2) + \exp(-\tilde{C}_0/8 \cdot n\mu_X(x) h_n^{d_x}) \nonumber \\
& \le C_1^* \exp(-C_2^* n \mu_X(x) h_n^{d_x} t^2), \label{eq:mur-kenel-general-case1}
\end{align}
where we set $C_1^* = \sum_{s=1}^{n - N_r + 1} C_1^{(s)} + 1$ and $C_2^* = \min\{-C_2^{(s)} \tilde{C}_0/4 \text{ for }  s=1,...,n - N_r + 1, ~\tilde{C}_0/8 \}$.

Concluding all three cases (expressions \eqref{eq:mur-kenel-naive-case1}, \eqref{eq:mur-kenel-box-case1}, and \eqref{eq:mur-kenel-general-case1}),
for $d_y=1$,
we establish a probability bound $\P(\Wscr_p(\mu_{Y|x}, \mur) > t) \le C_1^* \exp(-C_2^* n \mu_X(x) h_n^{d_x} t^2)$ ,
where $C_1^*$ and $C_2^*$ are constants independent of $x$, $n$, and $t$.
The results for other values of $d_y$ can be derived similarly.

\smallskip
\noindent
$\bullet~$ \textbf{Part 2.} For $\Wscr_p(\mur, \munp)$,
since $\mur  = \sum_{i = 1}^n w_i \delta_{\{f(x)+\eta_i\}}$ and $\munp  = \sum_{i = 1}^n w_i \delta_{y_i}$
share the same probability mass $w_i = \frac{K((x - x_i)/h_n)}{\sum_{j = 1}^n K((x-x_j)/h_n)}$,
we construct a probability coupling $\xi = \sum_{i = 1}^n w_i \delta_{\{(f(x) + \eta_i, y_i)\}}$ to bound the Wasserstein distance by definition as follows:
\begin{align}
\Wscr_p(\mur, \munp) & \le \left( \E_{(Y_1, Y_2)\sim \xi}\Big[\|Y_1 - Y_2\|_1^p\Big] \right)^{\frac{1}{p}}
= \bigg(\sum_{i=1}^n w_i \|f(x)-f(x_i)\|_1^p \bigg)^{\frac{1}{p}} \nonumber \\
& \le L h_{n}^{\beta} \left(\frac{\sum_{i=1}^n K((x-x_i)/h_n) \|(x-x_i)/h_n\|_2^{\beta p} }{\sum_{i = 1}^n K((x-x_i)/h_n)}  \right)^{\frac{1}{p}}.  \label{eq:mur-munp-step1}
\end{align}
Here, the inequality holds due the $(L, \beta)$-H\"older condition on function $f$ from \Cref{asp:regular_dist}.
Next, we bound \eqref{eq:mur-munp-step1} by analyzing its numerator and denominator terms separately.

\smallskip
\underline{The numerator of \eqref{eq:mur-munp-step1}}:
We define a set $\Iscr_{\tilde R} = \{i \in [n]: \|x - x_i\|_2 \leq \tilde R h_n\}$ with a radius $\tilde{R}>0$,
and then we decompose the term $\sum_{i=1}^n K((x-x_i)/h_n) \|(x-x_i)/h_n\|_2^{\beta p}$ as follows:
\begin{align}
\sum_{i \in \Iscr_{\tilde R}} K((x-x_i)/h_n) \|(x-x_i)/h_n\|_2^{\beta p} + \sum_{i \notin \Iscr_{\tilde R}} K((x-x_i)/h_n) \|(x-x_i)/h_n\|_2^{\beta p}. \label{eq:mur-np-step3}
\end{align}

For the first sum in~\eqref{eq:mur-np-step3},
we define $N_{\tilde{R}} := \sum_{i = 1}^n \I{\|x - x_i\|_2 \leq \tilde{R} h_n}$, so that it consists of $|\Iscr_{\tilde R}| = N_{\tilde{R}}$ terms.
Following \eqref{eq:lemma-number-near-sample-3} in \Cref{lemma:lower-bd-Nnr},
we can set $\tilde{C}_{\tilde R} = \frac{\pi^{d_x/2}}{\Gamma(d_x / 2 + 1)}\tilde R^{d_x}$ and obtain
\begin{align}
\label{eq:mur-np-step3-bias1-prob}
\P(N_{\tilde R} \geq 4 \tilde C_{\tilde R} n \mu_X(x) h_n^{d_x}) \leq \exp(-\tilde C_{\tilde R} n \mu_X(x) h_n^{d_x}/8).
\end{align}
Given $i \in \Iscr_{\tilde R}$, we have $\|(x-x_i)/h_n\|_2 \le \tilde R$ and $K((x-x_i)/h_n) \le b_R$.
With \eqref{eq:mur-np-step3-bias1-prob}, this implies:
\begin{align}
\label{eq:mur-np-step3-bias1}
\begin{split}
 & \P\bigg(\sum_{i \in \Iscr_{\tilde R}} K((x-x_i)/h_n) \|(x-x_i)/h_n\|_2^{\beta p} \geq b_R \tilde R^{\beta p} \cdot 4  \tilde C_{\tilde R} n\mu_X(x)h_n^{d_x}\bigg)  \\
\leq~ & \P(N_{\tilde R} \geq 4 \tilde C_{\tilde R} n \mu_X(x) h_n^{d_x}) \leq \exp(-\tilde C_{\tilde R} n \mu_X(x) h_n^{d_x}/8).
\end{split}
\end{align}

For the second sum in~\eqref{eq:mur-np-step3}, we analyze the three cases of the kernel functions separately:

Case (2): A box kernel function that satisfies \Cref{asp:kernel2}. (This subsumes Case (1) of the naive kernel function.)
We set $\tilde R = R$, so that $K((x-x_i)/h_n) = 0$ for $i \notin \Iscr_{\tilde R}$. Hence, we have
\begin{align}
\sum_{i \notin \Iscr_{\tilde R}} K((x-x_i)/h_n) \|(x-x_i)/h_n\|_2^{\beta p} = 0.
\label{eq:mur-np-step3-bias2-1}
\end{align}

Case (3) A general kernel function that satisfies \Cref{asp:kernel}.
We have
\begin{align}
\sum_{i \notin \Iscr_{\tilde R}} K((x-x_i)/h_n) \|(x-x_i)/h_n\|_2^{\beta p}
& \leq n \sup_{\|\tau\|_2 >\tilde{R}} K(\tau)\|\tau\|_2^{\beta p}
\leq n C_K \sup_{\tilde{r} > \tilde R} e^{-\tilde{r}} \tilde{r}^{\beta p}
\leq n C_K \exp(-\tilde R) \tilde R^{\beta p}\nonumber \\
& = n C_K h_n^{2 d_x}(- 2d_x \log h_n)^{\beta p} \leq C_K n \mu_X(x) h_n^{d_x}.\label{eq:mur-np-step3-bias2-2}
\end{align}
by setting $\tilde R = - 2d_x \log h_n = \Theta(d_x \psi \log n) > \beta p$.
Here, the third inequality follows from the fact that $f(\tilde{r}) = e^{-\tilde{r}} \tilde{r}^{\beta p}$ is decreasing for $\tilde{r} \in [\beta p, \infty)$.
The last inequality holds because $h_n^{d_x}=O(\mu_X(x))$.

By combining~\eqref{eq:mur-np-step3-bias1}, \eqref{eq:mur-np-step3-bias2-1}, and \eqref{eq:mur-np-step3-bias2-2},
we derive a bound for the numerator term:
\begin{equation}\label{eq:bias-target-numerator}
\hspace{-1.5cm} \P\bigg(\sum_{i \in [n]} K((x-x_i)/h_n) \|(x-x_i)/h_n\|_2^{\beta p} \geq (4 b_R \tilde C_{\tilde R} \tilde R^{\beta p} + C_K) n\mu_X(x)h_n^{d_x}\bigg) \leq \exp(-\tilde C_{\tilde R} n \mu_X(x) h_n^{d_x}/8),
\end{equation}
where $\tilde{C}_{\tilde R} = \frac{\pi^{d_x/2}}{\Gamma(d_x / 2 + 1)}\tilde R^{d_x}$ is independent of $n$ and $x$.
Note that, given the choice $\tilde R = \Theta(d_x \psi \log n)$ in Case (3), we suppress the dependence of $\tilde{C}_{\tilde R}$ on $\log n$, as it is of smaller order compared to $n h_n^{d_x}$.

\smallskip
\underline{The denominator of \eqref{eq:mur-munp-step1}}:
We have $\sum_{i=1}^n K((x-x_i)/h_n) \ge b_r\sum_{i = 1}^n \I{\|x - x_i\|_2 \leq r h_n} =  b_r N_r$.
By applying \eqref{eq:lemma-number-near-sample-2} from \Cref{lemma:lower-bd-Nnr}, we obtain
\begin{align}
\hspace{-1.8cm} \P\bigg(\sum_{i = 1}^n K((x - x_i)/h_n) \le b_r \cdot C_0 n \mu_X(x) h_n^{d_x}/4\bigg)
\le \P(N_r \le \tilde C_0 n \mu_X(x) h_n^{d_x}/4)
\le \exp(-\tilde C_0 n \mu_X(x) h_n^{d_x}/8). \label{eq:mur-np-step2}
\end{align}

Using the cases where the numerator takes large values, as given in \eqref{eq:bias-target-numerator},
and the denominator takes small values, as given in \eqref{eq:mur-np-step2}, we can derive the following  probability bound:
\begin{align}
\P \left(\frac{\sum_{i=1}^n K((x-x_i)/h_n) \|(x-x_i)/h_n\|_2^{\beta p} }{\sum_{i = 1}^n K((x-x_i)/h_n)} > C_0 \right) \leq 2\exp(-C_3^* n \mu_X(x) h_n^{d_x}),
\end{align}
where we set $C_0 = \frac{4(4 b_R\tilde C_{\tilde R} \tilde R^{\beta p} + C_K)}{\tilde C_0 b_r}$ and $C_3^* = \min\{\tilde C_{\tilde R}/8, ~\tilde C_0/8\}$.
Note that $C_0$ and $C_3^*$ are constants independent of $x$, $n$, and $t$.
We then revisit \eqref{eq:mur-munp-step1} and obtain:
\begin{align}\label{eq:bias-target2}
\P \big(\Wscr_p(\mur,\munp)> C_0 Lh_n^{\beta}\big) \leq 2\exp(-C_3^* n \mu_X(x) h_n^{d_x}).
\end{align}

\smallskip
\noindent
$\bullet~$ \textbf{Part 3.}
Finally, by combining the analysis for $\Wscr_p(\mu_{Y|x}, \mur)$ (results \eqref{eq:mur-kenel-naive-case1}, \eqref{eq:mur-kenel-box-case1}, and \eqref{eq:mur-kenel-general-case1}) and $\Wscr_p(\mur, \munp)$ (result \eqref{eq:bias-target2}),
we prove the first case of \eqref{eq:convergence_rate_np} in \Cref{thm:concentration-np} as follows:
\begin{align}
\P \big(\Wscr_p(\mu_{Y|x},\munp)>t + C_0 Lh_n^{\beta}\big) & \le \P\big(\Wscr_p(\mu_{Y|x},\mur)>t\big)+ \P \big(\Wscr_p(\mur,\munp)> C_0 Lh_n^{\beta}\big) \nonumber \\
&\le C_1^* \exp(-C_2^* n \mu_X(x) h_n^{d_x} t^2) + 2 \exp(-C_3^* n \mu_X(x) h_n^{d_x}) \nonumber\\
& \le C_1 \exp(-C_2 n \mu_X(x) h_n^{d_x} t^2),
\end{align}
where we set $C_1 = C_1^* + 2$ and $C_2 = \min\{C_2^*, ~C_3^*\}$. Note that $C_1$ and $C_2$ are constants independent of $x, n$, and $t$. The other two cases can be derived in a similar manner.
\hfill \halmos

\subsection{Measure Concentration Result of $k$NN}\label{app:kNN}
For the nonparametric estimator $\munp$,
we can incorporate $k$NN method as an alternative within our {\iwdro} framework.
Specifically, we sort the training samples $\{(x_i, y_i)\}_{i = 1}^{n}$ in non-decreasing order based on the distance $\|x - x_i\|_2$,
and define the $k$NN estimator as
\begin{equation}\label{eq:knn}
\munp \coloneqq \sum_{i = 1}^{k_n}\frac{1}{k_n}\delta_{y_i}.
\end{equation}
We assume that the number of selected neighbors, denoted by $k_n$, satisfies the following conditions.
\begin{assumption}\label{asp:knn}
As $n \to \infty$, it holds that $k_n \to \infty$ and $\frac{k_n}{n} \to 0$.
\end{assumption}
For example, choosing $k_n = \lfloor n^{\delta} \rfloor$ with some $0 < \delta < 1$ to satisfies Assumption~\ref{asp:knn}.

Similar to the measure concentration analysis performed for the NW kernel estimator,
we can derive a nonasymptotic rate for this $k$NN estimator using the Wasserstein distance.
\begin{theorem}[Measure Concentration of the $k$NN estimator]\label{prop:ws-concentration-knn}
Under Assumptions \ref{asp:regular_dist} and \ref{asp:knn},
for any $ x\in \Xscr$ and $t\in (0,1]$,
the $k$NN estimator $\munp$ defined in \eqref{eq:knn} satisfies:
\begin{equation}\label{eq:convergence_rate_knn}
\P \Bigg(\Wscr_p( \mu_{Y | x}, \munp)>t + C_0 L \bigg(\frac{k_n}{n \mu_X(x)}\bigg)^{\frac{\beta}{d_x}}\Bigg)\le
\begin{cases} C_1 \exp(-C_2 k_n t^2),  & \text{if } p > d_y / 2, \\
C_1 \exp(-C_2 k_n (t/\log (2 + 1/t))^2), & \text{if } p = d_y / 2, \\
C_1 \exp(-C_2 k_n t^{d_y/p}),  & \text{if } p \in [1, d_y / 2),\\
\end{cases}
\end{equation}
where $C_0$, $C_1$, and $C_2$ are some constants independent of $x$, $n$, and $t$.
\end{theorem}
\emph{Proof: }
Here, we prove the case of $d_y = 1$.
We define $\mur = \sum_{i=1}^{k_n} \frac{1}{k_n} \delta_{\{f(x) + \eta_i\}}$,
which is a uniform distribution over $k_n$ neighbors.
Since $k_n$ is predetermined and independent of $\eta_i$, following from \Cref{lemma:empirical-was-concentration}, we have:
\begin{equation}\label{eq:mur-knn-step1}
\P\big(\Wscr_p(\mu_{Y | x}, \mur) > t\big) \leq \tilde{C}_1 \exp(-\tilde{C}_2 k_n t^2).
\end{equation}

For $\Wscr_p(\mur, \munp)$,
we construct a probability coupling $\xi = \sum_{i = 1}^{k_n} \frac{1}{k_n} \delta_{\{(f(x) + \eta_i, y_i)\}}$ to bound the Wasserstein distance by definition as below
\begin{align}
\Wscr_p(\mur, \munp) & \le \left( \E_{(Y_1, Y_2)\sim \xi}\Big[\|Y_1 - Y_2\|_1^p\Big] \right)^{1/p}
= \bigg(\sum_{i=1}^{k_n} \frac{1}{k_n} \|f(x)-f(x_i)\|_1^p \bigg)^{1/p} \nonumber \\
& \le L \bigg(\sum_{i=1}^{k_n} \frac{1}{k_n} \|x-x_i\|_2^{\beta p} \bigg)^{1/p}  \le L\|x - x_{k_n}\|_2^{\beta}, \label{eq:mur-knn-step2}
\end{align}
where the last inequality follows from the non-decreasing order of $\|x - x_i\|_2$.

We define $N_{\zeta} := \sum_{i = 1}^{n} \I{\|x - x_i\|_2 \leq \zeta}$ with a radius $\zeta > 0$.
Using a similar tail estimation for $N_r$ as in \eqref{eq:lemma-number-near-sample-2} of \Cref{lemma:lower-bd-Nnr},
we can set $\tilde{C}_0 = \pi^{\frac{d_x}{2}} / \Gamma(\frac{d_x}{2} + 1)$ and derive the following bound:
\begin{align}
\P(N_{\zeta} \le \tilde{C}_0 n \mu_X(x) \zeta^{d_x}/4) \le \exp(-\tilde{C}_0 n\mu_X(x) \zeta^{d_x} / 8). \nonumber
\end{align}
By choosing $\zeta = (\frac{4 k_n}{\tilde{C}_0 n \mu_X(x)})^{\frac{1}{d_x}}$
so that $\tilde{C}_0 n \mu_X(x) \zeta^{d_x}/4 = k_n$, we obtain $\P(N_{\zeta} \le k_n) \le \exp(-k_n / 2).$
On the other hand, conditioned on $N_{\zeta} > k_n$, we have $\|x - x_{k_n}\|_2 \le \zeta$ by definition.
This implies $\P (\Wscr_p(\mur, \munp) > L\zeta^{\beta} | N_{\zeta} > k_n) = 0$ from \eqref{eq:mur-knn-step2}.
By combing these two bounds, we have
\begin{align}
\P (\Wscr_p(\mur, \munp) > L\zeta^{\beta}) \le \P (\Wscr_p(\mur, \munp) > L\zeta^{\beta} | N_{\zeta} > k_n) + \P(N_{\zeta} \le k_n)
\le \exp(-k_n / 2).  \label{eq:mur-knn-step3-1}
\end{align}
Next, we set $C_0 = (4/\tilde{C}_0)^{\frac{\beta}{d_x}}$ and substitute the choice $\zeta = (\frac{4 k_n}{\tilde{C}_0 n \mu_X(x)})^{\frac{1}{d_x}}$ in \eqref{eq:mur-knn-step3-1}, and then achieve
\begin{align}
\P \Bigg(\Wscr_p(\mur, \munp)> C_0 L \bigg(\frac{k_n}{n \mu_X(x)}\bigg)^{\frac{\beta}{d_x}}\Bigg) \le \exp(-k_n / 2).  \label{eq:mur-knn-step3}
\end{align}

Finally, by combining \eqref{eq:mur-knn-step1} with \eqref{eq:mur-knn-step3}, we obtain
\begin{align}
\P \Bigg(\Wscr_p( \mu_{Y | x}, \munp)>t + C_0 L \bigg(\frac{k_n}{n \mu_X(x)}\bigg)^{\frac{\beta}{d_x}}\Bigg) 
\le~ & \P(\Wscr_p(\mu_{Y | x}, \mur) > t) + \P \Bigg(\Wscr_p(\mur, \munp)> C_0 L \bigg(\frac{k_n}{n \mu_X(x)}\bigg)^{\frac{\beta}{d_x}}\Bigg) \nonumber \\
\le~ & \tilde{C}_1 \exp(-\tilde{C}_2 k_n t^2) + \exp(-k_n / 2) \nonumber \\
\le~ & C_1 \exp(-C_2 k_n t^2),  \nonumber
\end{align}
where we set $C_1 = \tilde{C}_1 + 1$ and $C_2 = \min\{\tilde{C}_2, ~1/2 \}$. Note that $C_1$ and $C_2$ are constants independent of $x$, $n$, and $t$.
An analogous analysis can be carried out for other values of $d_y$ by applying the corresponding cases in \eqref{eq:convergence_rate_knn}.
$\hfill \square$

\subsection{Proof of Example~\ref{prop:p-estimator-example}.}\label{app:p-estimator-proof}
We impose the following conditions for the case of $d_y = 1$.
For $d_y > 1$, we assume that these conditions hold for each coordinate of $y$.
\begin{assumption}\label{asp:add-asp}
In Example~\ref{prop:p-estimator-example},
the linear model $f_{\theta}$, the noise $\eta$, and the nonlinear term $g(x)$
satisfy conditions 1 - 4 in \cite{hsu2012random}.
Moreover, for all $x\in \Xscr$, there exist some constant $C$ such that $\|x\|_2 \leq C$ and $\|g(x)\|_2 \leq C$.
\end{assumption}

To upper bound $\Wscr_p(\mu_{Y|x},\mup)$,
we introduce an auxiliary distribution $\muapx^{(1)} \coloneqq \theta^{*\top}x  + g(\tilde{x}) + \eta$,
where $\tilde{x} \sim \mu_{X}$ follows the covariate distribution and is independent of both the realization of $x$ and $\eta$ in the training samples.
We use $\muapx^{(1)}$ as an intermediate measure to 
apply the triangle inequality
$$\Wscr_p(\mu_{Y|x}, \mup) \leq \Wscr_p(\mu_{Y|x}, \muapx^{(1)}) + \Wscr_p(\muapx^{(1)}, \mup).$$

By definition, we have $\mu_{Y|x} = \theta^{*\top}x  + g(x) + \eta$ and $\eapx = 2 \sup_{x \in \Xscr} \|g(x)\|_1$. Hence, we obtain
\begin{equation}\label{eq:example3-part1}
\Wscr_p(\mu_{Y|x}, \muapx^{(1)}) = \Wscr_p(\theta^{*\top}x  + g(x) + \eta, \theta^{*\top} x + g(\tilde{x})+ \eta) = \|g(x) - g(\tilde x)\|_1 \le \eapx.
\end{equation}

Next, we introduce another intermediate measure $\muapx^{(2)} = \frac{1}{n}\sum_{i = 1}^n\delta_{\{\theta^{*\top} x + g(x_i) + \eta_i\}}$ and apply
\begin{equation}
\Wscr_p(\muapx^{(1)}, \mup)\leq \Wscr_p(\muapx^{(1)}, \muapx^{(2)}) + \Wscr_p(\muapx^{(2)}, \mup). \nonumber
\end{equation}

For $\Wscr_p(\muapx^{(1)}, \muapx^{(2)})$,
$\{g(x_i) + \eta_i\}_{i = 1}^n$ can be viewed as i.i.d. samples from the distribution $g(\tilde{x}) + \eta$,
i.e., $\muapx^{(2)}$ can be viewed as the empirical distribution for $\muapx^{(1)}$.
Following \Cref{lemma:empirical-was-concentration}, we have
\begin{equation}\label{eq:example3-part2}
\P\big(\Wscr_p(\muapx^{(1)}, \muapx^{(2)}) \geq t_1\big) \leq \tilde{C}_1 \exp(-\tilde{C}_2 n t_1^2).
\end{equation}

For $\Wscr_p(\muapx^{(2)}, \mup)$,
given that $\mup = \frac{1}{n}\sum_{i = 1}^n \delta_{\{\hat\theta^{\top} x + y_i - \hat\theta^{\top} x_i \}}$ and $\muapx^{(2)}$
share the same probability mass,
we construct a probability coupling $\xi = \sum_{i = 1}^n \frac{1}{n} \delta_{\{(\theta^{*\top} x + g(x_i) + \eta_i, ~\hat\theta^{\top} x + y_i - \hat\theta^{\top} x_i)\}}$ to bound the Wasserstein distance by definition as below
\begin{align}
\Wscr_p(\muapx^{(2)}, \mup) & \le \left( \E_{(Y_1, Y_2)\sim \xi}\Big[\|Y_1 - Y_2\|_1^p\Big] \right)^{\frac{1}{p}}
= \bigg(\sum_{i=1}^n \frac{1}{n} \|(\theta^{*\top} x + g(x_i) + \eta_i) - (\hat\theta^{\top} x + y_i - \hat\theta^{\top} x_i) \|_1^p \bigg)^{\frac{1}{p}} \nonumber \\
& = \bigg(\frac{1}{n} \sum_{i=1}^n \big|(\theta^*-\hat\theta)^{\top} x - (\theta^*-\hat\theta)^{\top} x_i) \big|^p \bigg)^{\frac{1}{p}} \nonumber \\
& \le \bigg(\frac{2^{p - 1}}{n}\sum_{i = 1}^n \big|(\theta^*-\hat\theta)^{\top} x \big|^p + \frac{2^{p - 1}}{n}\sum_{i = 1}^n \big| (\theta^*-\hat\theta)^{\top} x_i) \big|^p \bigg)^{\frac{1}{p}}. \nonumber
\end{align}
We define $\Sigma = \E_{X \sim \mu_X}[X X^{\top}]$,
and let $\lambda_{\min}(\Sigma)$ denote the smallest eigenvalue of $\Sigma$.
Additionally, we denote the vector norm $\|x\|_{\Sigma} = \sqrt{x^{\top}\Sigma x}$.
Then, we have
\begin{align}
\Wscr_p(\muapx^{(2)}, \mup) \le \bigg(\frac{2^{p}}{n}\sum_{i = 1}^n \|\theta^*-\hat\theta\|_{\Sigma}^p\|x_i\|_{\Sigma^{-1}}^p\bigg)^{\frac{1}{p}}
\le \frac{2C}{\sqrt{\lambda_{\min}(\Sigma)}}\|\theta^*-\hat\theta\|_{\Sigma}, \nonumber
\end{align}
where the second inequality follows from
$\|x\|_{\Sigma^{-1}}= \sqrt{x^{\top}\Sigma^{-1}x} \leq \frac{\|x\|_2}{\sqrt{\lambda_{\min}(\Sigma)}} \leq  \frac{C}{\sqrt{\lambda_{\min}(\Sigma)}}$
for all $x \in \Xscr$.
Furthermore, by Remark 12 from \cite{hsu2012random},
with probability at least $1- 3 \exp(-n \tilde{t}^2)$,
we have $\|\theta^*-\hat\theta\|_{\Sigma} \leq \sqrt{2 \max\{d_x, \lambda_{\min}(\Sigma)^{-1}\}} C \tilde{t}$, which leads to
\begin{align}\label{eq:example3-part3}
\P\big(\Wscr_p(\muapx^{(2)}, \mup) \ge t_2\big) \leq 3 \exp(-\tilde{C}_3 n t_2^2),
\text{ where } \tilde{C}_3 = \frac{\lambda_{\min}(\Sigma)}{8C^4\max\{d_x, \lambda_{\min}(\Sigma)^{-1}\}}.
\end{align}

Finally, by choosing $t_1 = t/2$ in \eqref{eq:example3-part2} and $t_2 = t/2$ in \eqref{eq:example3-part3}, we have
\begin{align*}
\P\big(\Wscr_p(\mu_{Y|x}, \mup) \geq t + \eapx\big)
& \leq \P\big(\Wscr_p(\mu_{Y|x}, \muapx^{(1)}) \geq \eapx\big) + \P(\Wscr_p\big(\muapx^{(1)}, \muapx^{(2)}\big) \ge t/2) + \P\big(\Wscr_p(\muapx^{(2)}, \mup) \ge t/2\big) \nonumber \\
& \leq 0 + \tilde{C}_1 \exp(-\tilde{C}_2 n t^2/4) + 3 \exp(-\tilde{C}_3 n t^2/4) \le C_4 \exp(-C_5n t^2),
\end{align*}
where we set $C_4 = \tilde{C}_1 + 3$ and $C_5 = \min\{\tilde{C}_2/4, \tilde{C}_3/4 \}$. Note that $C_4$ and $C_5$ are constants independent of $x$, $n$, and $t$.
$\hfill\square$

\section{Proofs in Section~\ref{sec:approximation-formulation}}
\subsection{Proof of Theorem~\ref{thm:coverage-iwdro-apx}.}

We first show  $\setaiw \subseteq \setame$ by definition.
For any distribution $\mu \in \setaiw$, we have $\Wscr_p(\mu, \munp) \leq \varepsilonnp$ and $\Wscr_p(\mu, \mup) \leq \varepsilonp$.
By the convexity of Wasserstein distance given in Lemma~\ref{lemma:was-convexity}, we have
\begin{align*}
\Wscr_p(\mu, \mume) & \leq \big(\kappa_x\Wscr_p^p(\mu, \munp) + (1-\kappa_x) \Wscr_p^p(\mu,\mup)\big)^{\frac{1}{p}}
\leq \big(\kappa_x \varepsilonnp^p + (1-\kappa_x) \varepsilonp^p\big)^{\frac{1}{p}} = \varepsilonme,
\end{align*}
and thus $\mu \in \setame$.
Under Assumptions~\ref{asp:regular_dist} and~\ref{asp:kernel},
\Cref{thm:choice-size} shows $\P(\mu_{Y|x} \in \setaiw) \geq 1- \alpha$,
which directly implies $\P(\mu_{Y|x} \in \setame) \geq 1- \alpha$.
$\hfill \square$

\subsection{Proof of~\Cref{prop:optimality-gap-iwdro-approximate}.}

We proceed under the event $\mu_{Y|x} \in \setame$, which occurs with probability at least $1-\alpha$.
We first show \eqref{eq:optimality-iwdro-approximate-bound}.
Since $\ziw$ is optimal for $\sup_{\mu \in \setaiw}\E_{Y\sim \mu}[c(z, Y)]$ and
$\setaiw \subseteq \setame$, it immediately follows that
$$\sup_{\mu \in \setaiw} \E_{Y\sim\mu}[c(\ziw, Y)] \le \sup_{\mu \in \setaiw} \E_{Y\sim\mu}[c(\zme, Y)] \le \sup_{\mu \in \setame} \E_{Y\sim\mu}[c(\zme, Y)].$$

Next, we establish~\eqref{eq:optimality-iwdro-approximate-gap}.
Referring to the analysis in \eqref{eq:single-ball-dro-gap-2} from \Cref{prop:optimality-gap-iwdro}, we have
$$\E_{\mu_{Y|x}}[c(\zme, Y)] - \E_{\mu_{Y|x}}[c(z^{*}, Y)] \le 2 \mathcal{L}_{z^*} \varepsilonme.$$
Furthermore, we can derive
$\varepsilonme  = \big(\kappa_x \varepsilonnp^p + (1-\kappa_x) \varepsilonp^p\big)^{\frac{1}{p}} \leq \Cme \min\{\varepsilonnp, \varepsilonp\} /2$,
where $\Cme > 2$ is a constant independent of $n$.
This is equivalent to show that $\varepsilonme$ and $\min\{\varepsilonnp, \varepsilonp\}$ are of the same order with respect to $n$ under $\tilde\Theta(\cdot)$.

We prove the case for $p < d_y / 2$, and analogous analysis can be performed for other values of $p$.
Recall that \Cref{thm:choice-size} provides
$\varepsilonnp=  \Theta\Para{\Para{\frac{\rho_x}{n}}^{\rnp}}$ with $\rnp=\frac{\beta}{\beta d_y/p + d_x}$
and $\varepsilonp= \Theta\Para{\Para{\frac{1}{n}}^{\rp}} + \eapx$ with $\rnp = \frac{p}{d_y}$.
We also have $h_n = \Theta\big((n\mu_X(x))^{-\frac{1}{\beta d_y/p + d_x}} \big)$.
By applying \Cref{lemma:lower-bd-Nnr} and $\rho_x = \Theta(\frac{1}{\mu_X(x)})$,
we have $N_r = \sum_{i = 1}^n \I{\|x - x_i\|_2 \leq rh_n} = \Theta(n\mu_X(x) h_n^{d_x})= \Theta((n/\rho_x)^{\frac{\beta d_y/ p}{\beta d_y / p + d_x}})$.
Then, by definition \eqref{eq:mixture-weight} and $\rme = -\frac{p^2}{d_y}$, we obtain
$$\kappa_x = \max\bigg\{1 - \Theta\bigg(\Big(\frac{\rho_x}{n}\Big)^{\frac{\beta p}{\beta d_y / p + d_x}}\bigg), ~0\bigg\}.$$

We analyze $\kappa_x$  based on the size of $\rho_x$.
Specifically, we consider $\rho_x = \Theta(n^{\gamma})$ for $\gamma \in [0, \infty)$.
We have $\varepsilonp = \Theta(n^{-\frac{p}{d_y}}) + \eapx = \tilde\Theta(1)$,
given that the model mispecification error $\eapx$ is assumed to be of constant error $\tilde\Theta(1)$.
(In the following, we denote $f(n) = o(g(n))$ and $g(n) = \omega(f(n))$
if, for any $k>0$, there exists some constant $n_0$ such that $f(n) \leq k g(n)$ for all $n \geq n_0$.)

\begin{itemize}

\item If $\gamma \in [0, 1)$, then $\rho_x = o(n)$ and $\Theta\bigg(\Big(\frac{\rho_x}{n}\Big)^{\frac{\beta p}{\beta d_y / p + d_x}}\bigg) = o(1)$.
We have $\kappa_x = 1 - o(1)$ and $\varepsilonnp = o(1)$.
Hence, $\min\{\varepsilonnp, \varepsilonp\} = \varepsilonnp = \Theta\bigg(\Big(\frac{\rho_x}{n}\Big)^{\frac{p}{\beta d_y / p + d_x}}\bigg)$
and $\varepsilonme = \Big((1 - o(1)) \varepsilonnp^p + o(1) \varepsilonp^p\Big)^{\frac{1}{p}} \approx \varepsilonnp$.

\item If $\gamma =1$, then $\rho_x = \Theta(n)$ and $\Theta\bigg(\Big(\frac{\rho_x}{n}\Big)^{\frac{\beta p}{\beta d_y / p + d_x}}\bigg) = \Theta(1)$.
We have $\kappa_x \in (0, 1)$ and $\varepsilonnp = \Theta(1)$.
Hence, all of $\varepsilonnp$, $\varepsilonp$,
and $\varepsilonme = \Big(\kappa_x \varepsilonnp^p + (1-\kappa_x) \varepsilonp^p\Big)^{\frac{1}{p}} = \Theta(1)$
are of constant sizes.

\item If $\gamma >1$, then $\rho_x = \omega(n)$ and $\Theta\bigg(\Big(\frac{\rho_x}{n}\Big)^{\frac{\beta p}{\beta d_y / p + d_x}}\bigg) = \omega(1)$.
We have $\kappa_x  = \max\{1-\omega(1),0\}$ and $\varepsilonnp = \omega(1)$.
Hence, $\min\{\varepsilonnp, \varepsilonp\} = \varepsilonp = \tilde\Theta(1)$
and $\varepsilonme = \Big(\kappa_x \varepsilonnp^p + (1-\kappa_x) \varepsilonp^p\Big)^{\frac{1}{p}} \approx \varepsilonp$.
$\hfill \square$
\end{itemize}

\section{Additional Details for \Cref{sec:numerical}}\label{app:numerical}

\subsection{Income Prediction}
\label{app:numerical-prediction}
\paragraph{Hyperparameter Selection.}
Given an experiment instance, for each model, we conduct a 4-fold cross-validation to choose the radius $\varepsilon$ that minimizes the average objective value within the training dataset.
For \textsf{NP-DRO}, we set $h_n = 20 n^{-\frac{1}{d_x + 2}}$ and consider $\varepsilon = k / (\sum_{i \in [n]}K((x - x_i)/h_n))$ with $k \in \{5, 10, 20, 40\}$.
For \textsf{P-DRO}, we consider $\varepsilon \in \{1,2,5,10\}$.


\subsection{Portfolio Allocation}
\label{app:numerical-portfolio}

\paragraph{Reformulation.} Program \eqref{eq:mean-cvar-problem2} has a piecewise linear cost function as follows:
\begin{align}
\min_{z\in \Zscr, r} \sup_{\mu \in \setaiw} ~\E_{Y\sim\mu} \bigg[\max\Big\{-\big(1 + 1/\phi\big) z^{\top}Y + \big(1 - 1/\phi\big) r, ~ -z^{\top}Y + r\Big\}\bigg]. \nonumber
\end{align}
According to \Cref{coro:linear}, it can be simplified to a linear program:
\begin{equation}\label{eq:cvar-mean-reformulation}
\begin{aligned}
\min \quad &\lambda_1 \varepsilon_1 +\lambda_2 \varepsilon_2  + \sum_{i = 1}^{n}w_{1,i} u_i +  \sum_{j = 1}^n w_{2,j} v_j\\
\text{s.t.}\quad & u_i + v_j \geq (1-1/\phi) r - (\alpha_{i,1}^{\top}y_{1,i} + \beta_{j,1}^{\top} y_{2,j}), \quad \forall i, j \in [n];\\
& u_i + v_j \geq r - (\alpha_{i,2}^{\top}y_{1,i} + \beta_{j,2}^{\top} y_{2,j}), \quad \forall i, j \in [n];\\
&\alpha_{i,1} + \beta_{j, 1} = \Para{1 + 1/\phi}z, \quad \forall i, j \in [n];\\
&\alpha_{i,2} + \beta_{j,2} =  z, \quad \forall i, j \in [n];\\
&\|\alpha_{i,1}\|_{\infty}, \|\alpha_{i,2}\|_{\infty} \leq \lambda_1, ~ \|\beta_{j,1}\|_{\infty}, \|\beta_{j,2}\|_{\infty} \leq \lambda_2, \quad \forall i,j \in [n];\\
& \lambda_1, \lambda_2 \geq 0; \quad u_i, v_i \in \R, ~~\forall i \in [n]; \\
& \sum_{i \in [d_z]} z_i = 1, ~z_i \ge 0.
\end{aligned}
\end{equation}

\paragraph{Hyperparameter Selection.}
In the real-data study, we determine the radius of the ambiguity set using 4-fold cross-validation, selecting the configuration that maximizes the empirical Sharpe Ratio over the training samples.
For {\iwdro}, \textsf{NP-DRO}, and \textsf{P-DRO}, we adopt the same hyperparameter structures as in the synthetic-data study.
For \textsf{OTCMV} and \textsf{RNW}, we calibrate the radius directly.
The hyperparameter configurations are provided in \Cref{tab:model-exp-config}.

\begin{table}[htb]
\SingleSpacedXI
\small
\renewcommand{\arraystretch}{1.1}
\renewcommand{\tabcolsep}{5mm}
\caption{Model hyperparameter configuration in the real-data study of portfolio allocation}
\label{tab:model-exp-config}
\begin{center}
\begin{tabular}{l l}
\toprule
Model & Hyperparameter tuning set\\
\midrule
\textsf{IW-DRO} &$k_1 \in \{0.4, 0.8, 1.2\}, ~k_2 \in \{0.002, 0.005, 0.1\}$\\
\textsf{NP-DRO} & $k \in \{0.2, 0.4, 0.6, 0.8, 1, 1.2\}$\\
\textsf{P-DRO} & $\varepsilon \in \{0.5, 1, 2\}$ \\
\textsf{OTCMV} & $\varepsilon \in \{0.1, 0.2, 0.5\}, ~a \in \{1.1, 1.2, 1.5\}$ \\
\textsf{RNW} & $\varepsilon \in \{0.01, 0.05, 0.1, 0.2, 0.5, 1\}$ \\
\bottomrule
\end{tabular}
\end{center}
\footnotesize
\end{table}


\end{document}